\documentclass[preprint,12pt,3p]{elsarticle}
\usepackage[T1]{fontenc}
\usepackage[utf8]{inputenc}
\usepackage{amsmath,amssymb,amsthm,mathtools}
\usepackage{enumitem}
\usepackage{booktabs}
\usepackage{multirow}
\usepackage{verbatim}
\usepackage{graphicx}
\usepackage{subcaption}
\usepackage{xcolor}
\usepackage{algorithm}
\usepackage{algpseudocode}
\usepackage{lineno}
\usepackage[expansion=false]{microtype}
\usepackage{hyperref}
\hypersetup{
  colorlinks=true,
  linkcolor=red,
  filecolor=green,
  urlcolor=green,
  anchorcolor=green,
  citecolor=cyan,
  pdftitle={Structure-Preserving Nonstandard Finite Difference Methods for Port-Hamiltonian Systems},
}
\allowdisplaybreaks

\DeclareMathOperator{\diag}{diag}
\DeclareMathOperator{\blkdiag}{blkdiag} 

\DeclareMathOperator{\spec}{spec} 
\DeclareMathOperator{\sym}{sym}

\DeclareMathOperator{\skewpart}{skew} 
\DeclareMathOperator{\ran}{ran} 
\DeclareMathOperator{\ke}{ker} 

\DeclareMathOperator{\Real}{Re}
\DeclareMathOperator{\Imag}{Im}

\DeclareMathOperator{\tanhc}{tanhc} 

\newcommand{\R}{\mathbb{R}}
\newcommand{\C}{\mathbb{C}}
\newcommand{\Z}{\mathbb{Z}}

\newcommand{\dgrad}{\overline{\nabla}} 
\newcommand{\Jb}{\overline{J}}
\newcommand{\Rb}{\overline{R}}
\newcommand{\Gb}{\overline{G}}

\newcommand{\nph}{n+\frac{1}{2}} 

\newcommand{\iu}{\mathrm{i}} 
\newcommand{\dd}{\mathrm{d}}

\newcommand{\pHNSFD}{\mbox{PH-NSFD}}

\theoremstyle{definition}
\newtheorem{definition}{Definition}[section]
\newtheorem{assumption}[definition]{Assumption}
\newtheorem{example}[definition]{Example}

\theoremstyle{remark}
\newtheorem{remark}[definition]{Remark}

\theoremstyle{plain}
\newtheorem{proposition}[definition]{Proposition}
\newtheorem{theorem}[definition]{Theorem}
\newtheorem{corollary}[definition]{Corollary}
\newtheorem{lemma}[definition]{Lemma}

\usepackage[figcolor=white]{todonotes}

\biboptions{sort&compress}
\begin{document}
\begin{frontmatter}
\title{Structure-Preserving Nonstandard Finite Difference Methods \\for Port-Hamiltonian Systems}

\author[BUW]{Matthias Ehrhardt\corref{Corr}}
\cortext[Corr]{Corresponding author}
\ead{ehrhardt@uni-wuppertal.de}

\author[FPT]{Manh Tuan Hoang}
\ead{tuanhm16@fe.edu.vn}

\address[BUW]{University of Wuppertal, Applied and Computational Mathematics,\\ 
Gau\ss{}stra\ss{}e 20, 42119 Wuppertal, Germany}

\address[FPT]{Department of Mathematics, FPT University, Hoa Lac Hi-Tech Park, \\
Km29 Thang Long Blvd, Hanoi, Viet Nam}

\begin{abstract}
We propose a structure-preserving, nonstandard finite difference (NSFD) methodology for numerically integrating linear and nonlinear port-Hamiltonian (pH) systems.
This approach incorporates Mickens's nonlocal discretization principles, discrete-gradient representations of the Hamiltonian, and structure-preserving approximations of the interconnection and dissipation operators.
Under symmetry and positive semidefiniteness conditions on the discrete structure operators, we establish an exact discrete power balance and discrete passivity for every solved step.
Equilibria of the continuous system are exact equilibria of the step equation of the scalar scheme. 
At nondegenerate energy minima, the linearization of the scheme is a Cayley transform that reproduces the linear stability type for every admissible positive denominator. 
Broader stability statements are made only under explicitly stated solvability and spectral conditions.
The NSFD denominator function provides an additional degree of freedom that can be selected to reproduce characteristic frequencies, dissipation rates, or the linearized flow at an equilibrium while preserving the pH algebraic structure.


The central idea is to distinguish structural preservation, which is achieved through discrete gradients, from dynamic consistency, which is achieved through nonstandard denominator functions.
We derive sufficient conditions for energy, passivity, invariance, and equilibrium preservation, and we construct denominator functions from the spectrum of the linearized pH vector field.
Since the discrete-gradient time-stepping is of the midpoint type, the locally exact scalar denominators belong to the hyperbolic tangent family
$\phi(h) = \frac{2}{\lambda}\tanh\bigl( \frac{\lambda h}{2} \bigr)$, rather than the exponential family  $\phi(h) = (\mathrm{e}^{\lambda h}-1)/\lambda$,  which is commonly used for Euler-type NSFD schemes. 
The latter reduces the order of the method from two to one.
For linear systems, exact matrix denominators lead to step-dependent pH representations whose structure operators converge to continuous ones of second order.

Numerical experiments on conservative, dissipative, nonlinear, and forced benchmark systems confirm the theory and quantify the improvement over standard discrete-gradient integrators.

\end{abstract}

\begin{keyword}
port-Hamiltonian systems \sep nonstandard finite difference \sep discrete gradients \sep structure-preserving integration \sep passivity \sep dynamic consistency \sep locally exact integrators 

\textit{2020 Mathematics Subject Classification:} 65L05, 65L12, 65P10, 37M15, 93C10, 93D05.
\end{keyword}

\journal{ESAIM Mathematical Modelling and Numerical Analysis} 


\end{frontmatter}

\section{Introduction}\label{Sec1}

\textit{Port-Hamiltonian (pH) systems} provide a unified framework for modeling interconnected physical systems in terms of energy storage, energy routing, dissipation and external power exchange \cite{DuindamEtAl2009, MehrmannUnger2023, vanderSchaft2017, vanderSchaftJeltsema2014}.
A finite-dimensional nonlinear port-Hamiltonian system can be written in the input--state--output 
form
\begin{equation}\label{eq:continuous_ph}
    \dot{x} = \bigl( J(x)-R(x) \bigr)\nabla H(x) + G(x)u,
    \qquad y = G(x)^\top \nabla H(x),
\end{equation}
where $x\in\R^n$ denotes the state,
$H\colon\R^n\to\R$ is the Hamiltonian (the stored energy),
$J(x)$ is the skew-symmetric interconnection matrix,
$R(x)$ is the positive semidefinite dissipation matrix,
and $u,y\in\R^m$ are the collocated input and output variables.
The defining structural assumptions are
\begin{equation}\label{eq:ph_structure}
    J(x) = -J(x)^\top, \qquad R(x)=R(x)^\top,
    \qquad R(x)\succeq 0 .
\end{equation}
As an immediate consequence one obtains the \textit{power balance}
\begin{equation}  \label{eq:power_balance}
    \frac{\dd}{\dd t}H(x(t))
    =    -\nabla H(x(t))^\top R(x(t))\nabla H(x(t))   + y(t)^\top u(t)
    \le y(t)^\top u(t) .
\end{equation}
In particular the closed system ($u\equiv 0$) is dissipative, and the open system is passive with respect to the \textit{supply rate} $y^\top u$ \cite{vanderSchaft2017, Willems1972}.
Port-Hamiltonian modeling has become a standard formalism in multiphysics simulation, electrical circuits \cite{FalaizeHelie2016}, mechanics, thermodynamics, power networks and fluid dynamics, see \cite{BeattieEtAl2018, DuindamEtAl2009, MehrmannUnger2023, vanderSchaftJeltsema2014} and the references therein.
Its two decisive features are that the structure \eqref{eq:ph_structure} is invariant under power-preserving interconnection and under Galerkin projection, and that stability and passivity statements are available directly from the algebraic structure, without solving the equations.
Moreover, 
pH systems can be coupled and decoupled in a structure-preserving way \cite{Ehrhardt2026}, such that one consider a large monolithic system or a group of smaller systems for each physical process.

%
%
It is essential to preserve these properties after discretization. 
Artificial energy production yields incorrect long-time behavior, even for small local errors. 
Artificial dissipation predicts spurious asymptotic stability. 
A single non-passive component can destabilize a simulated interconnection.

Geometric numerical integration \cite{HairerLubichWanner2006, LeimkuhlerReich2004, SanzSernaCalvo1994} provides two complementary families of methods for such problems.
\textit{Symplectic integrators} preserve the symplectic two-form of a canonical Hamiltonian system and, by backward error analysis, conserve a modified Hamiltonian over exponentially long times.
\textit{Discrete-gradient} (energy-preserving) methods \cite{CelledoniEtAl2012, Gonzalez1996, ItohAbe1988, McLachlanQuispelRobidoux1999,
McLachlanQuispel2014, QuispelMcLaren2008} replace $\nabla H$ by a nonlocal, two-argument object $\dgrad H(x,z)$ satisfying a discrete chain rule, and thereby conserve the energy exactly rather than approximately.  
Because the port-Hamiltonian power balance \eqref{eq:power_balance} is an energy statement rather than a symplectic one, and because dissipation destroys symplecticity anyway, discrete gradients are the natural point of departure for pH systems.
Related constructions preserve Lyapunov functions \cite{GrimmQuispel2005}, several first integrals simultaneously
\cite{DahlbyOwrenYaguchi2011, NortonQuispel2014}, or the GENERIC structure of thermodynamic systems \cite{Oettinger2018}.

For pH systems specifically, structure-preserving time discretizations have been developed along several lines:
Aoues et al.\ \cite{AouesEtAl2017} study losslessness, passivity and composability of discrete-time approximations;
Celledoni and H{\o}iseth \cite{CelledoniHoiseth2017} construct energy-preserving and passivity-consistent discretizations based on discrete gradients and averaged vector fields;
Kotyczka and Lef\`evre \cite{KotyczkaLefevre2018, KotyczkaLefevre2019} define discrete-time pH systems through symplectic collocation and a discrete-time Dirac structure;
Mehrmann and Morandin \cite{MehrmannMorandin2019} and Kinon, Morandin and Schulze \cite{KinonMorandinSchulze2026} treat pH descriptor systems;
Schulze \cite{Schulze2023} introduces discrete gradient \textit{pairs} for systems with a state-dependent mass matrix;
D'Ambrosio and Di Donato \cite{Dambrosio2026} study discrete-gradient $\theta$-methods;
Giesselmann, Karsai and Tscherpel \cite{GiesselmannKarsaiTscherpel2025} give an energy-consistent Petrov--Galerkin time discretization;
and Kinon et al.\ \cite{KinonEtAl2023} apply the framework to nonlinear elastodynamics.
Structure-preserving spatial semi-discretizations of distributed pH systems are reviewed in \cite{Kotyczka2019}, see also \cite{CardosoRibeiroMatignonLefevre2021,  Egger2019,
KotyczkaMaschkeLefevre2018, MoullaLefevreMaschke2012, TrenchantEtAl2018}.

Structural preservation alone does not determine dynamical accuracy.
For example, when using the standard denominator $\varphi(h)=h$, the AVF discrete-gradient method applied to a harmonic oscillator reduces to the implicit midpoint method. 
Although it preserves the quadratic energy up to round-off errors, it introduces an $O(h^2)$ phase error over one period. Additionally, the unwrapped phase lag 
increases linearly over time.

A conceptually different route to qualitative correctness is provided by the \textit{nonstandard finite difference} (NSFD) 
approach pioneered by Mickens \cite{Mickens1994, 
Mickens2002, Mickens2020}.
Its central tenet is that a numerical approximation should reproduce the essential qualitative properties of the underlying differential equation rather than merely approximate derivatives locally. 
These Mickens' rules are described in detail in Section~\ref{sec:nsfd}.
The theoretical underpinning of the method (consistency, elementary stability, and the notion of \textit{dynamic consistency}) was developed by Anguelov and Lubuma \cite{AnguelovLubuma2001, AnguelovLubuma2003} and by Mickens \cite{Mickens2005dc}, cf.\ also the surveys in \cite{HoangMatthias2026,  
Patidar2016}.
Systematic constructions of denominator functions under positivity constraints are given in \cite{Mickens2007denominator, Mickens2013}, and second-order variants in \cite{Alalhareth2025, AnguelovEtAl2022, HoangMatthias2024, HoangMatthias2026, KojouharovEtAl2021}, and high-order mass- and positivity-preserving variants for production--destruction systems in \cite{Conte2026}. 
NSFD schemes have been 
successful in mathematical biology, epidemiology, finance and transport problems \cite{EhrhardtMickens2013, Maamar2024, MickensWashington2020, Morais2024, Treibert2022, ZinihiEhrhardtSidiAmmi2026}, where positivity, boundedness and correct equilibrium structure are decisive.

Mickens also considered conservative oscillators and one-dimensional discrete Hamiltonian systems \cite{Mickens1996energy, Mickens2005hamilton, Mickens2020}, 
and Songolo and Bid\'egaray-Fesquet \cite{SongoloBidegaray2024} recently reformulated Mickens' rules for systems of ODEs, 
showing that the resulting schemes are naturally interpreted as exponential integrators acting on the linear part. 
A closely related and, for the present work, particularly important line is Cie\'sli\'nski's theory of \textit{locally exact} integrators \cite{Cieslinski2013camwa, Cieslinski2013pla}, 
in which a scheme is modified so as to reproduce the linearization of the vector field exactly at every point, while retaining energy conservation.
Related linearization-preserving integrators are studied in \cite{McLachlanQuispelTse2009}.
The same circle of ideas underlies exponential integrators \cite{HochbruckOstermann2010}, exponential fitting \cite{IxaruVandenBerghe2004} and trigonometrically fitted methods \cite{Gautschi1961, HochbruckLubich1999}.


Here we will combine these two structural viewpoints: NSFD supplies a mechanism for structure-preservation and for matching dynamical time scales, while the pH formalism supplies the algebraic invariants that must survive the discretization step.
We propose a class of \textit{port-Hamiltonian nonstandard finite difference methods} (\pHNSFD{}) 
in which different components of the discrete scheme are assigned distinct 
structural roles:
\begin{equation*}
\begin{array}{rcl}
\text{discrete gradient }\dgrad H
&\longrightarrow &\text{exact discrete chain rule and energy balance},\\[1mm]
\text{skew-symmetric interconnection }\Jb_n
&\longrightarrow & \text{lossless power exchange},\\[1mm]
\text{positive semidefinite dissipation }\Rb_n
&\longrightarrow & \text{discrete passivity}, \\[1mm]
\text{nonlocal Casimir representation }\dgrad C
&\longrightarrow& \text{invariant preservation},\\[1mm]
\text{denominator function }\phi(h)
&\longrightarrow&\text{dynamic consistency and spectral matching}.
\end{array}
\end{equation*}
%
The last line distinguishes the proposed approach from standard discrete-gradient methods: the denominator is used to reproduce characteristic frequencies, dissipation rates and, at a designated equilibrium, the exact linearized flow.
This freedom is \textit{structurally free of charge}, since replacing $h$ by $\phi(h)$ does not touch any algebraic identity on which the structural theorems rest, but it is not free with respect to accuracy.

For autonomous problems with a fixed step size, the scalar scheme is the underlying discrete-gradient update with step
length $\varphi(h)$ and its states are assigned to the physical times $t_n=nh$ (see Lemma~\ref{lem:reduction}). 
The denominator therefore calibrates the time scale of the update. 
In general, a shared scalar denominator cannot fit distinct modes independently. 
However, the matrix extensions allow for separate modal fitting.

The main contributions of this work are the following.
We develop a structure-preserving NSFD framework for nonlinear pH systems, with well-posedness of the implicit step and a consistency and order analysis.
We distinguish three levels of qualitative statements:
(i) an \textit{exact discrete power balance} and \textit{discrete passivity}, which hold for every solved step; 
(ii) stability results under stated assumptions, namely a discrete Lyapunov and LaSalle-type analysis and the preservation of the linear stability type at nondegenerate energy minima for every admissible positive denominator; and 
(iii) exact matrix fitting for linear systems on an explicitly specified nonresonant range of step sizes.


A second group of results concerns the construction of denominator functions.
Since the discrete-gradient time stepping is of \textit{midpoint type}, the locally exact scalar denominators are $\phi(h)=\frac{2}{\lambda}\tanh\bigl(\frac{\lambda h}{2}\bigr)$, whereas the exponential denominator $(\mathrm{e}^{\lambda h}-1)/\lambda$ of Euler-type NSFD schemes generically reduces the order from two to one.
For damped oscillatory modes, we provide a scalar denominator that exactly matches the rotation angle.
For matrix denominators, we introduce a class of \textit{step-dependent discrete pH methods}.
The exact midpoint-type matrix denominator and locally exact energy-preserving schemes for canonical systems are known \cite{Cieslinski2013camwa, Cieslinski2013pla}.
We derive their pH representation with second-order consistent operators and nonnegative dissipation, equilibrium compatibility conditions and step-size domains.
Additionally, we provide a step-dependent pH representation with positive semidefinite dissipation for a rational approximation that coincides with the two-stage Gauss method for linear problems.
%
%
Numerical experiments confirm the theory for conservative, dissipative, nonlinear, and forced benchmarks. These experiments include a long-time integration in which the \pHNSFD{} scheme's final-time error is 10 orders of magnitude smaller than the implicit midpoint rule's error.

The paper is organized as follows.
Sections~\ref{sec:ph}--\ref{sec:dgrad} collect the required material on pH systems, NSFD schemes and discrete gradients.
The \pHNSFD{} scheme is defined and analyzed in Section~\ref{sec:scheme}, and its structural properties are proved in Section~\ref{sec:properties}.
Sections~\ref{sec:denominator} and~\ref{sec:locallyexact} are devoted to scalar and matrix denominator functions, Section~\ref{sec:models} to model problems, and Section~\ref{sec:numerics} to numerical experiments.

\section{Port-Hamiltonian Systems}\label{sec:ph}
We consider \textit{port-Hamiltonian (pH) systems} of the form \eqref{eq:continuous_ph}, i.e.
\begin{equation} \label{eq:ph_system}
    \dot{x} = \bigl( J(x)-R(x) \bigr)\nabla H(x) + G(x)u,
    \qquad x(0)=x_0,
\end{equation}
with collocated output
\begin{equation} \label{eq:ph_output}
    y  = G(x)^\top \nabla H(x).
\end{equation}
In the sequel we have to assume the following properties.
\begin{assumption}
\label{ass:main}
Let $\Omega\subseteq\R^n$ be open and connected.  
We assume throughout that
\begin{enumerate}[label=(A\arabic*),leftmargin=2.6em]
  \item $H\in C^2(\Omega;\R)$ and $\nabla H$ is locally Lipschitz on $\Omega$;
  \item $J,R\colon\Omega\to\R^{n\times n}$ and $G\colon\Omega\to\R^{n\times m}$ are locally Lipschitz;
  \item for all $x\in\Omega$,  we have for interconnection matrix $J$ and dissipation matrix $R$
        \begin{equation}\label{eq:structure}
            J(x)^\top = -J(x), \qquad 
            R(x)^\top = R(x), \qquad R(x)\succeq0 ;
        \end{equation}
  \item $u\in L^{\infty}_{\mathrm{loc}}([0,\infty);\R^m)$.
\end{enumerate}
If $H$ is bounded from below on $\Omega$ and has compact sublevel sets, we call $H$ \textit{coercive}.
\end{assumption}

Assumption~\ref{ass:main} suffices for all algebraic results below. 
For $u\in L^\infty_{\mathrm{loc}}$, solutions are understood in the Carath\'eodory sense, identities hold almost everywhere, and the discrete input is the cell average, since point values of $u$ are not defined.
For the order, linearization and quadrature results we additionally assume, where stated, that $H\in C^3$, $J,R,G\in C^2$ locally, and $u\in C^2$ on the time interval considered.

\begin{lemma}[Power balance]\label{lem:powerbalance}
Every classical solution of \eqref{eq:ph_system}--\eqref{eq:ph_output} satisfies
\begin{equation}\label{eq:cont_pb}
    \dot{H} = -\nabla H^\top R(x)\nabla H + y^\top u
    = -\mathcal{D}(x) + y^\top u,
    \quad
    \mathcal{D}(x) := \nabla H(x)^\top R(x)\nabla H(x)\ge 0 .
\end{equation}
In particular $\dot{H}\le y^\top u$, and $\dot{H} \le 0$ for $u\equiv0$.
If $R\equiv0$ and $u\equiv0$, then the Hamiltonian (the 'energy') is conserved, i.e.\ $H(x(t))=H(x_0)$ for all~$t$.
\end{lemma}
\begin{proof}
Differentiating $H(x(t))$ and inserting \eqref{eq:ph_system} yields
\begin{equation*}
   \dot{H} = \nabla H^\top \dot{x} = \nabla H^\top J \nabla H - \nabla H^\top R \nabla H+\nabla H^\top G u.
\end{equation*}
The first term vanishes because the interconnection matrix $J$ is skew-symmetric, so we have $\nabla H^\top J \nabla H = \tfrac{1}{2} \nabla H^\top (J+J^\top) \nabla H=0$, and the third term equals $y^\top u$ by \eqref{eq:ph_output}.
\end{proof}

Integrating the first identity in \eqref{eq:cont_pb} gives the
(integral form of the) \textit{dissipation inequality}
\begin{equation}\label{eq:cont_passivity}
    H(x(t_1)) - H(x(t_0)) \le \int_{t_0}^{t_1} y(s)^\top u(s)\,\dd s ,
\end{equation}
so that $H$ is a \textit{storage function} (energy) for the \textit{supply rate} $y^\top u$ and \eqref{eq:ph_system} is \textit{passive}, whenever $H$ is bounded from below \cite{vanderSchaft2017, Willems1972}. 
It is precisely these structures \eqref{eq:cont_pb}--\eqref{eq:cont_passivity} that we shall reproduce \textit{exactly}, not merely to leading order, at the discrete level.

We note that \emph{linear pH systems} take the form
\begin{equation}\label{eq:linear_ph}
    \dot{x} = (J-R) Q x + Bu, \qquad y = B^\top Q x,
    \qquad H(x)=\tfrac{1}{2} x^\top Q x,
\end{equation}
with constant $J=-J^\top$, $R=R^\top\succeq0$, $Q=Q^\top\succ0$.
They arise both as models in their own right and as linearizations of \eqref{eq:ph_system} at an equilibrium. 
The following elementary spectral result for linear systems will be used repeatedly.

\begin{lemma}[Spectrum of a linear pH vector field]\label{lem:spectrum}
Let $A:=(J-R)Q$ with $J=-J^\top$, $R=R^\top\succeq0$, $Q=Q^\top\succ0$.  
Then
\begin{equation}
    \spec(A)\subseteq\overline{\C}_-,
\end{equation}
where $\mathbb C^-:=\{\lambda\in\mathbb C\colon \Real\lambda\le0\}$ denotes as usual the closed left half-plane.
$\lambda\in\spec(A)$ with $\Real\lambda=0$ implies that the corresponding eigenvector $v$ satisfies $R Q v = 0$. 
If $R=0$, then $\spec(A)\subseteq \iu\R$ and $A$ is diagonalizable with purely imaginary eigenvalues.
\end{lemma}
\begin{proof}
$A$ is similar to $\widetilde{A} := Q^{1/2}AQ^{-1/2} = Q^{1/2}(J-R)Q^{1/2}$, whose symmetric part is $-Q^{1/2}RQ^{1/2}\preceq0$.  
Hence for any unit eigenvector $w$ of $\widetilde{A}$ with eigenvalue $\lambda$, $\Real\lambda = \Real(w^{*}\widetilde{A}w)=w^{*}\sym(\widetilde{A})w\le0$.
Equality forces $Q^{1/2}RQ^{1/2}w=0$, i.e.\ $RQv=0$ for $v=Q^{-1/2}w$.
If $R=0$ then $\widetilde{A}$ is real skew-symmetric, hence normal with purely imaginary spectrum.
\end{proof}

\begin{lemma}[Contractivity in the energy norm]\label{lem:contract}
For $A$ as in Lemma~\ref{lem:spectrum}, we have that $\mathrm{e}^{tA^\top}Q\,\mathrm{e}^{tA}\preceq Q$ for all $t\ge0$.
\end{lemma}
\begin{proof}
$\frac{\dd}{\dd t}\bigl(\mathrm{e}^{tA^\top}Q\mathrm{e}^{tA}\bigr)
 = \mathrm{e}^{tA^\top}\bigl(A^\top Q+QA\bigr)\mathrm{e}^{tA}
 = -2\,\mathrm{e}^{tA^\top}QRQ\,\mathrm{e}^{tA}\preceq0$,
since $A^\top Q+QA=Q(J-R)^\top Q+Q(J-R)Q=-2QRQ$.
\end{proof}

%
Geometrically, \eqref{eq:ph_system} is the coordinate expression of a Dirac structure \cite{vanderSchaftJeltsema2014}.
For our purposes the relevant consequence is the algebraic identity $\nabla H^\top J\nabla H=0$, due to the skew-symmetry of the interconnection matrix $J$.
A discrete counterpart is noted in Remark~\ref{rem:dirac}.

\begin{definition}[Casimir]\label{def:casimir}
A function $C\in C^1(\Omega;\R)$ is a \textit{Casimir} of the interconnection structure, if
\begin{equation}\label{eq:casimir_cont}
    \nabla C(x)^\top J(x) =0 \qquad\text{for all } x\in\Omega .
\end{equation}
If in addition $\nabla C^\top R=0$ and $\nabla C^\top G=0$, then $C$ is conserved along all solutions of \eqref{eq:ph_system}, for every Hamiltonian $H$ and every input $u$.
\end{definition}

If only \eqref{eq:casimir_cont} holds, then $\dot{C}=-\nabla C^\top R\nabla H+\nabla C^\top Gu$, and $C$ is in general not conserved.
Casimirs are conserved because of the structure, independently of $H$, which makes them valuable in control \cite{vanderSchaftJeltsema2014} and their discrete preservation a requirement separate from energy conservation.

\begin{definition}[Forced equilibrium]\label{def:equilibrium}
A pair $(x^\ast,u^\ast)\in\R^n\times \R^n$ is an \textit{equilibrium} of \eqref{eq:ph_system} if
\begin{equation}\label{eq:equilibrium}
    \bigl( J(x^\ast)-R(x^\ast) \bigr)\nabla H(x^\ast) + G(x^\ast)u^\ast = 0 .
\end{equation}
\end{definition}

For $u^\ast = 0$ every critical point of $H$ is an \emph{equilibrium}, and for a constant, non-zero external port input $u^\ast\in\R^n$ it is a \emph{forced equilibrium}.
If $x^\ast$ is a strict local minimum of $H$ then $H-H(x^\ast)$ is a \emph{local Lyapunov function}. 
Asymptotic stability follows under a \emph{zero-state detectability} 
condition \cite[Section~3.2]{vanderSchaft2017}.

In the sequel we briefly collect the examples used throughout the paper.

\begin{example}[Damped harmonic oscillator {\cite[Example~1.1.]{Moench2024}}]
\label{ex:damped}
With $x=(q,p)^\top$,
\begin{equation}\label{eq:damped_oscillator}
    \dot{q} = p, \qquad 
    \dot{p} = -\omega^2 q-\gamma p, \qquad 
    H(q,p) = \tfrac{1}{2}\omega^2q^2 + \tfrac{1}{2} p^2 ,
\end{equation}
\begin{equation}\label{eq:damped_JRG}
    J = \begin{pmatrix}0&1\\-1&0\end{pmatrix},\quad
    R = \begin{pmatrix}0&0\\0&\gamma\end{pmatrix},\quad
    G = \begin{pmatrix}0\\1\end{pmatrix},\quad
    \nabla H(q,p)=\begin{pmatrix}\omega^2q\\p\end{pmatrix}.
\end{equation}
This is the workhorse example: $\gamma=0$ gives the conservative oscillator,
$\gamma>0$ the dissipative one, and $u\neq0$ the forced one.
\end{example}

\begin{example}[Nonlinear pendulum {\cite[Example~15.2.]{vanderSchaftJeltsema2014}, \cite[Section~8.4]{Strogatz}}]\label{ex:pendulum}
We consider $\dot{q} = p$, $\dot{p} = -\sin q$, with $J$ as above, $R=0$ and
\begin{equation}\label{eq:pendulum}
    H(q,p) = \tfrac{1}{2} p^2+1-\cos q .
\end{equation}
The linearization at the stable equilibrium $(0,0)$ has frequency $\omega_0=1$.
It tests the nonlinear discrete-gradient construction, the preservation of the energy landscape (separatrix at $H=2$) \cite[Sect.~8.4]{Strogatz}, and the benefit of locally exact denominators away from the linear regime.
\end{example}


\begin{example}[Rigid body / Casimir test {\cite[Example~8.1.]{vanderSchaftJeltsema2014}}]\label{ex:rigidbody}
The free rigid body $\dot{x} = J(x)\nabla H(x)$ with $J(x) = \widehat{x}$ (the cross-product matrix), $H(x)=\tfrac{1}{2}\sum_j x_j^2/I_j$, possesses the Casimir $C(x)=\tfrac{1}{2}\|x\|^2$, since $\nabla C^\top J=x^\top\widehat{x}=0$.  
Adding $R=\rho\,\widehat{x}\,\widehat{x}^\top\succeq0$, $\rho\ge0$, models a damping which still annihilates $\nabla C$. Then $C$ remains conserved, whereas $H$ decreases in general.
This is the reference problem for Section~\ref{sec:casimir}.
\end{example}

\section{The Nonstandard Finite Difference Principle}\label{sec:nsfd}
In the NSFD framework, introduced by Mickens \cite{Mickens1994, Mickens2020}, the standard first-order difference quotient $\dot{x}(t_n)\approx(x_{n+1}-x_n)/h$ is replaced by $\dot{x}(t_n) \approx (x_{n+1}-x_n)/\phi(h)$,
where the \textit{denominator function} $\phi$ satisfies the conditions
\begin{equation}   \label{eq:phi_conditions}
    \phi\in C^1(0,h_{\max}),\qquad
    \phi(h)>0,\qquad
    \phi(h) = h+\mathcal{O}(h^2)\quad (h\to0^+) .
\end{equation}

Mickens' general NSFD ansatz also permits a numerator function $\psi(h)$ in the difference quotient $(x_{n+1}-\psi(h)x_n)/\phi(h)$.
Within the present framework this is \emph{not} an additional degree of freedom: consistency already requires $\psi(h)=1+o(h)$, and exact preservation of a nonzero equilibrium forces $\psi\equiv1$ (Proposition~\ref{prop:rigidity} in \ref{app:numerator}).
Therefore, the denominator $\phi(h)$ is the relevant scalar NSFD degree of freedom.

\subsection{The Design Rules for NSFD Schemes}
Mickens' design rules \cite{AnguelovLubuma2001, Mickens2002, Patidar2016} require that the order of the discrete derivative equals that of the continuous one, that $h$ is replaced by a nontrivial denominator $\phi(h)$ as in \eqref{eq:phi_conditions}, that nonlinear terms are approximated nonlocally, and that essential properties of the exact solution (e.g.\ positivity, convexitiy, conservation laws, equilibria) are inherited without spurious solutions or instabilities.
This is formalized by the notion of \textit{dynamic consistency}.

\begin{definition}[Dynamic consistency \cite{AnguelovLubuma2001, Mickens2005dc}]\label{def:dc}
Let $\mathcal{P}$ be a property of the differential equation (e.g.\ positivity, boundedness, monotone energy decay, stability type of the equilibria, oscillation frequency, decay rate).  
A discretization is \textit{dynamically consistent with respect to $\mathcal{P}$} if the discrete system possesses $\mathcal{P}$ for all step sizes $h>0$ in the admissible range.
\end{definition}

\begin{definition}[Elementary stability \cite{AnguelovLubuma2001}]\label{def:elemstab}
A one-step method $x_{n+1}=\Psi_h(x_n)$ for the ODE $\dot{x} = f(x)$ is \textit{elementarily stable} if, for every $h>0$, the map $\Psi_h$ is well defined, its fixed points $x^\ast$ are exactly the equilibria of the right-hand side $f$, and each of them has the same linear stability type for $\Psi_h$ as for the flow of~$f$.
\end{definition}

For implicit schemes the well-definedness of $\Psi_h$ (and, if necessary, the choice of solution branch) is part of the statement, and even the exact time-$h$ flow map can have fixed points that are not equilibria (Proposition~\ref{prop:nonresonance}).
The same happens for A-stable methods: the stability function of the three-stage Gauss method equals $1$ at $z=\pm\iu\sqrt{60}$, so for the harmonic oscillator with $h\omega=\sqrt{60}$ every state is a fixed point. 
Restrictions of this kind are stated explicitly below.

\subsection{Exact Schemes and the Origin of the Denominator Function}
The origin of the denominator function lies in \textit{exact finite difference schemes}.
The prototype is the linear scalar 
ODE $\dot{x} = \lambda\, x$, whose exact
solution satisfies $x(t_{n+1}) = \mathrm{e}^{\lambda h}x(t_n)$. 
This is reproduced \textit{exactly} by the explicit NSFD scheme
\begin{equation}\label{eq:exact_euler}
    \frac{x_{n+1}-x_n}{\phi_{\mathrm{E}}(h)}=\lambda x_n,
\quad\text{with}\quad
    \phi_{\mathrm{E}}(h)=\frac{\mathrm{e}^{\lambda h}-1}{\lambda},
\end{equation}
which for $\lambda=-\alpha<0$ (Dahlquist test equation) gives the classical Mickens denominator
$\phi_{\mathrm{E}}(h) = (1-\mathrm{e}^{-\alpha h})/\alpha$.
The essential point is that $\phi_{\mathrm{E}}$ is \emph{exact for the forward Euler scheme}.
For the midpoint-type discrete-gradient scheme used here, the exact denominator is $\frac{2}{\lambda}\tanh\bigl(\frac{\lambda h}{2}\bigr)$, and $\phi_{\mathrm{E}}$ generically reduces the order from two to one (Section~\ref{sec:denominator}).

Exact difference schemes exist for several linear constant-coefficient systems \cite{DangHoang2018, Hoang2023, Mickens1994, Roeger2008, SongoloBidegaray2024} and some nonlinear scalar equations \cite{Mickens2020}.
The NSFD methodology retains as much of this \emph{exactness} as possible when no exact scheme is available.
Mickens and Ramadhani~\cite{MickensRamadhani1994} chose the denominators of nonlinear NSFD schemes such that they are exact for the linearized equations, which preserves the linear stability properties for all step sizes.

The second NSFD ingredient is the \textit{nonlocal approximation} $f(x_n)\approx\widetilde{f}(x_n,x_{n+1})$ with $\widetilde{f}(x,x) = f(x)$, chosen such that a desired qualitative property is inherited \cite{AnguelovLubuma2003}.
For pH systems the canonical choice is dictated by the energy, i.e.\ by a discrete chain rule, which leads to discrete gradients.

\section{Discrete Gradients and Discrete Structure Operators} \label{sec:dgrad}

Let us start with the introduction of discrete gradients.
\begin{definition}[Discrete gradient]\label{def:dgrad}
A continuous map $\dgrad H\colon\Omega\times\Omega\to\R^n$ is called a \textit{discrete gradient} of $H$, if
\begin{align}
    H(z)-H(x) &= \dgrad H(x,z)^\top (z-x) &&\text{for all } x,z\in\Omega,
    \label{eq:discrete_gradient}\\
    \dgrad H(x,x) &= \nabla H(x) &&\text{for all } x\in\Omega .
    \label{eq:consistency_gradient}
\end{align}
It is called \textit{symmetric} if $\dgrad H(x,z) = \dgrad H(z,x)$ for all $x,z\in\Omega$.
\end{definition}

Identity \eqref{eq:discrete_gradient} is a \textit{discrete chain rule}, exact along every step independently of the step size. This is what allows us to insert $\phi(h)$ later without affecting any conservation statement.

\begin{example}[Standard discrete gradients]\label{ex:dgrads}\mbox{}
\begin{enumerate}[label=(\alph*),leftmargin=2.2em]
  \item \textit{Average vector field (AVF) / mean value} \cite{QuispelMcLaren2008,CelledoniEtAl2012}:
    \begin{equation}\label{eq:avf}
        \dgrad H^{\mathrm{AVF}}(x,z) = \int_0^1 \nabla H\bigl(x+s(z-x)\bigr)\,\dd s .
    \end{equation}
  \item \textit{Gonzalez / midpoint} \cite{Gonzalez1996}:
    \begin{equation}\label{eq:gonzalez}
        \dgrad H^{\mathrm{G}}(x,z) = \nabla H\Bigl(\frac{x+z}{2}\Bigr)
        +\frac{H(z)-H(x)-\nabla H\bigl(\frac{x+z}{2}\bigr)^\top(z-x)}{\|z-x\|^2}\,(z-x),
    \end{equation}
    for $z\neq x$, extended by $\nabla H(x)$ for $z=x$.
  \item \textit{Itoh--Abe / coordinate increment} \cite{ItohAbe1988}:
    the $j$-th component is the difference quotient of $H$ in the $j$-th
    coordinate direction, taken after the first $j-1$ coordinates are updated.
\end{enumerate}
The AVF and Gonzalez gradients are symmetric. The Itoh--Abe gradient is not, but
it is derivative-free and leads to a scheme that can be solved componentwise.
The AVF gradient requires the segment $[x,z]$, and the Gonzalez gradient and the midpoint evaluations \eqref{eq:midpointops} below require the midpoint, to lie in $\Omega$.
Hence ``for all $x,z\in\Omega$'' is read as ``for all admissible pairs $(x,z)\in\mathcal{A}:=\{(x,z)\in\Omega^2:[x,z]\subset\Omega\}$'' (note that $\mathcal{A}=\Omega^2$ for convex $\Omega$). 
For the Itoh--Abe discrete gradient, admissible pairs must
additionally have all intermediate coordinate-update points
and coordinate segments used in its definition contained in
$\Omega$. 
This additional requirement is automatic when $\Omega=\mathbb R^n$.
\end{example}

\begin{lemma}[Quadratic Hamiltonians]\label{lem:quadratic}
Let $H(x)=\frac{1}{2} x^\top Qx + b^\top x + c$ with $Q=Q^\top$. Then
\begin{equation}\label{eq:quadratic_dg}
    \dgrad H^{\mathrm{AVF}}(x,z)=\dgrad H^{\mathrm{G}}(x,z)=Q\,\frac{x+z}{2}+b .
\end{equation}
For the Itoh--Abe gradient with coordinate updates in the order $1,2,\dots,n$, we write \linebreak $Q=D+L+L^\top$ with $D$ diagonal and $L$ strictly lower triangular. Then we have
\begin{equation}\label{eq:itohabe_quadratic}
    \dgrad H^{\mathrm{IA}}(x,z)=Q\,\frac{x+z}{2}+b+\frac{1}{2}\bigl(L-L^\top\bigr)(z-x).
\end{equation}
In particular, all three discrete gradients coincide with the midpoint gradient if $Q$ is diagonal.
\end{lemma}
\begin{proof}
Since $\nabla H$ is affine, the integral in \eqref{eq:avf} equals the average of the endpoint values.
For \eqref{eq:gonzalez}, the identity $H(z)-H(x)=\bigl(Q\frac{x+z}{2}+b\bigr)^\top(z-x)$ holds exactly for quadratic $H$, so the numerator of the correction term vanishes.
For the Itoh--Abe gradient $\dgrad H^{\mathrm{IA}}$, let $w^{(j)}=(z_1,\dots,z_j,x_{j+1},\dots,x_n)^\top$.
Since $H$ is quadratic, $H(w^{(j)})-H(w^{(j-1)}) = (z_j-x_j) \bigl( Q\frac{w^{(j-1)}+w^{(j)}}{2} +b \bigr)_j$, which gives the $j$-th component of \eqref{eq:itohabe_quadratic}.
\end{proof}

%
The additional skew-symmetric term in \eqref{eq:itohabe_quadratic} does not contribute to the discrete chain rule, but it is nonzero in general.
Continuity and the discrete chain rule do not determine the rate at which $\dgrad H(x,z)$ approaches $\nabla H$ near the diagonal, 
since a continuous term $w(x,z)$ of low regularity with $w(x,z)^\top(z-x)=0$ and $w(x,x)=0$ can always be added.
For order results we therefore impose the approximation property explicitly.

\begin{definition}[Second-order midpoint approximation]\label{def:dg2}
A discrete gradient $\dgrad H$ satisfies 
approximation property (DG2) on a compact set $K\subset\Omega$ if there is a constant $C_K$ such that
\begin{equation}\label{eq:DG2}
    \Bigl\|\dgrad H(x,z)-\nabla H\Bigl(\frac{x+z}{2}\Bigr)\Bigr\|\le C_K\|z-x\|^2
    \qquad\text{for all } x,z\in K .\tag{DG2}
\end{equation}
Likewise, admissible operators $\Jb$, $\Rb$, $\Gb$ (Definition~\ref{def:discreteops} below) satisfy the condition (DG2) on $K$ if $\Jb(x,z)=J\bigl(\frac{x+z}{2}\bigr)+\mathcal{O}(\|z-x\|^2)$ uniformly on $K$, and similarly for $\Rb$ and $\Gb$.
\end{definition}

\begin{lemma}[Approximation order]\label{lem:dgradorder}
Let $H$ be of class $C^3$ on a neighborhood of a compact convex set $K\subset\Omega$.
Then $\dgrad H^{\mathrm{AVF}}$ and $\dgrad H^{\mathrm{G}}$ satisfy the condition \eqref{eq:DG2} on $K$.
More generally, \eqref{eq:DG2} holds for every symmetric discrete gradient of class $C^{1,1}$ near the diagonal.
If $\dgrad H$ is merely locally Lipschitz, then $\dgrad H(x,z)=\nabla H(x)+\mathcal{O}(\|z-x\|)$ on compact sets.
\end{lemma}

\begin{proof}
Let $m=\frac{x+z}{2}$ and $\delta=z-x$.
Taylor expansion about $m$ gives, for the AVF gradient, a linear term proportional to $\int_0^1(\theta-\frac{1}{2})\,\dd\theta=0$, and, for the Gonzalez gradient, $H(z)-H(x)-\nabla H(m)^\top\delta=\mathcal{O}(\|\delta\|^3)$, so that its correction term is $\mathcal{O}(\|\delta\|^2)$.
For a symmetric $C^{1,1}$ discrete gradient, $\delta\mapsto\dgrad H(m-\frac{1}{2}\delta,m+\frac{1}{2}\delta)$ is even with value $\nabla H(m)$ at $\delta=0$, so its derivative vanishes there.
The last statement follows from consistency and Lipschitz continuity.
\end{proof}

The midpoint evaluations \eqref{eq:midpointops} below satisfy the condition (DG2) trivially, and the endpoint averages $\frac{1}{2}( J(x)+J(z) )$ do so if $J$ is of class $C^2$.
Symmetry and continuity of the operators alone do not ensure a second-order estimate.
The power-balance results of Section~\ref{sec:properties} do not require (DG2) and therefore retain a broader scope than the order results.

In floating-point arithmetic, the Gonzalez gradient $\dgrad H^{\mathrm{G}}$ suffers from cancellation, when $\|z-x\| \approx \sqrt{\varepsilon_{\mathrm{mach}}}$.
The AVF gradient avoids this whenever its integral can be evaluated exactly, e.g.\ for a polynomial $H$ of degree $d$ by an $N$-point Gauss--Legendre rule with $N\ge\lceil d/2\rceil$ (2 nodes for a quartic $H$) \cite{CelledoniEtAl2012}.
Inexact quadrature causes a chain-rule defect, see \eqref{eq:residual}.

Alongside the discrete gradient we require nonlocal representations of $J$, $R$ and $G$.

\begin{definition}[Admissible discrete structure operators]\label{def:discreteops}
Maps $\Jb,\Rb\colon \Omega\times\Omega\to\R^{n\times n}$, $\Gb\colon\Omega\times\Omega\to\R^{n\times m}$ are called \textit{admissible}, 
if they are continuous and satisfy 
the \textit{structural conditions} 
\begin{equation}\label{eq:discrete_J}
    \Jb(x,z)^\top=-\Jb(x,z),\qquad
    \Rb(x,z)^\top=\Rb(x,z),\qquad
    \Rb(x,z)\succeq0,\quad\forall x,z\in\Omega,
\end{equation}
together with the \textit{consistency conditions}
\begin{equation}\label{eq:discrete_consistency}
    \Jb(x,x)=J(x),\qquad \Rb(x,x)=R(x),\qquad \Gb(x,x)=G(x).
\end{equation}
They are called \textit{symmetric} if $\Jb(x,z)=\Jb(z,x)$, $\Rb(x,z)=\Rb(z,x)$, $\Gb(x,z)=\Gb(z,x)$.
\end{definition}

The simplest admissible choices are the \textit{midpoint evaluations}
\begin{equation}\label{eq:midpointops}
    \Jb(x,z)=J\Bigl(\frac{x+z}{2}\Bigr),\quad
    \Rb(x,z)=R\Bigl(\frac{x+z}{2}\Bigr),\quad
    \Gb(x,z)=G\Bigl(\frac{x+z}{2}\Bigr),
\end{equation}
and the \textit{symmetrized} evaluations $\Jb(x,z) = \tfrac{1}{2} \bigl( J(x)+J(z) \bigr)$, and similarly for $R$ and $G$.
Both are symmetric in the sense of Definition~\ref{def:discreteops}, 
and both inherit \eqref{eq:discrete_J} pointwise from \eqref{eq:structure} because skew-symmetry, symmetry and positive semidefiniteness are preserved under convex combinations.
This is a small but important observation: unlike energy conservation, the structural conditions \eqref{eq:discrete_J} come for free.
%
Non-symmetric admissible choices such as $\Jb(x,z)=J(x)$ are also allowed, but lead to first-order schemes.

\section{The PH-NSFD Scheme}\label{sec:scheme}
Here, we introduce the \pHNSFD{} scheme and analyze its basic mathematical features.
\begin{definition}[\pHNSFD{} scheme]\label{def:phnsfd}
Let $\phi$ be any denominator function satisfying \eqref{eq:phi_conditions}. Let $\dgrad H$ be a discrete
gradient of $H$ and $\Jb$, $\Rb$, $\Gb$ be admissible in the sense of Definition~\ref{def:discreteops}. 
The \emph{port-Hamiltonian nonstandard finite difference (PH-NSFD)} scheme is defined by
\begin{equation} \label{eq:ph_nsfd}
    \frac{x_{n+1}-x_n}{\phi(h)}
    = \bigl(\Jb_n - \Rb_n \bigr) \dgrad H_n + \Gb_n u_{\nph}, \qquad
    y_{\nph} = \Gb_n^\top \dgrad H_n ,
\end{equation}
where $\dgrad H_n=\dgrad H(x_n,x_{n+1})$, $\Jb_n=\Jb(x_n,x_{n+1})$, $\Rb_n=\Rb(x_n,x_{n+1})$, $\Gb_n=\Gb(x_n,x_{n+1})$ and $u_{\nph}$ is a
representative value of the input on $[t_n,t_{n+1}]$, e.g.\ $u_{\nph}=u(t_n+h/2)$ or $u_{\nph}=\frac{1}{h}\int_{t_n}^{t_{n+1}}u(s)\,\dd s$.
The scheme is called \textit{symmetric}, if $\dgrad H$, $\Jb$, $\Rb$, $\Gb$ are symmetric and $\phi$ extends to an odd function, $\phi(-h)=-\phi(h)$, of class $C^3$ near $h=0$.
\end{definition}

We refer to the schemes of Definition~\ref{def:phnsfd} as the \textit{scalar class}.
Its defining feature, which will be essential for the preservation of equilibria, is the exact diagonal consistency \eqref{eq:consistency_gradient}, \eqref{eq:discrete_consistency} of all discrete operators, for every step size.
Step-dependent matrix extensions, for which this property is not automatic, are treated separately in Section~\ref{sec:locallyexact}.
Throughout, a statement that holds for every step refers to every $h$ in the stated range for which $\phi(h)$ is finite and positive and the implicit step has a solution in the domain under consideration. 
If the step equation has several solutions, the chosen branch is specified.

Setting $\phi(h)\equiv h$ recovers the discrete-gradient methods of \cite{Gonzalez1996, McLachlanQuispelRobidoux1999} in the pH form used by
\cite{CelledoniHoiseth2017, KinonMorandinSchulze2026, Schulze2023}.
If, in addition, $H$ is quadratic, the discrete operators are the constant continuous ones, $\Jb=J$, $\Rb=R$, $\Gb=G$, and the AVF or Gonzalez gradient is used, the scheme reduces to the implicit midpoint rule.
For the Itoh--Abe gradient this holds only for a diagonal quadratic form (Lemma~\ref{lem:quadratic}). 
The scalar class 
contains these methods as special cases and extends them by one scalar degree of freedom, namely the denominator function $\phi(h)$.

\begin{lemma}[Reduction lemma]\label{lem:reduction}
The \pHNSFD{} scheme \eqref{eq:ph_nsfd} with step size $h>0$ and denominator function $\phi(h)$ coincides with the discrete-gradient scheme ($\phi\equiv\mathrm{id}$)
applied with the modified step size $\widetilde{h}:=\phi(h)$ and input value $u_{\nph}$.
\end{lemma}
\begin{proof}
Immediate: \eqref{eq:ph_nsfd} reads $x_{n+1} = x_n + \widetilde{h} \bigl[ (\Jb_n-\Rb_n)\dgrad H_n + \Gb_nu_{\nph} \bigr]$.
\end{proof}

%
%
%
%
Consequently, all algebraic structural results of the discrete-gradient method that hold for every step size (energy balance, passivity, Casimir conservation, equilibrium preservation) transfer to the \pHNSFD{} scheme with $h$ replaced by $\phi(h)$, and so do existence and uniqueness results (Theorem~\ref{thm:wellposed}).
Accuracy, however, does not transfer: the scheme advances the solution by an amount appropriate to a step of length $\phi(h)$, but interprets it at $t_n+h$, and the discrepancy $h-\phi(h)$ enters the error directly (Theorem~\ref{thm:order}).

Next, we address the well-posedness of the implicit step.
\begin{theorem}[Solvability of the nonlinear step]\label{thm:wellposed}
Let $x_n\in\Omega$, and let $K\subseteq\Omega$ be nonempty
and closed with $\{x_n\}\times K\subseteq\mathcal A$.
Fix the step input $u=u_{n+1/2}$ with
$\|u_{n+1/2}\|\le U$, and let
\begin{equation}
    \mathcal{F}(x,z) := \bigl( \Jb(x,z)-\Rb(x,z) \bigr)\dgrad H(x,z) +\Gb(x,z)u ,
    \qquad
    \Theta(z) := x_n +\phi(h)\,\mathcal{F}(x_n,z).
\end{equation}
Assume that $z\mapsto\mathcal{F}(x_n,z)$ is Lipschitz continuous on $K$ with constant $L_K$, uniformly for $\|u\|\le U$, and that
\begin{equation}\label{eq:contraction}
    \Theta(K)\subseteq K
    \qquad\text{and}\qquad
    \phi(h)\,L_K<1 .
\end{equation}
Then the step equation \eqref{eq:ph_nsfd} has exactly one solution $x_{n+1}\in K$, and the fixed-point iteration $z^{(k+1)}=\Theta(z^{(k)})$ converges to it from every starting point $z^{(0)}\in K$.
\end{theorem}
\begin{proof}
The step equation is equivalent to the fixed-point equation $z=\Theta(z)$.
As a closed subset of $\R^n$, $K$ is complete.
By \eqref{eq:contraction}, $\Theta$ maps $K$ into itself and satisfies $\|\Theta(z_1)-\Theta(z_2)\| \le \phi(h)L_K\|z_1-z_2\|$ with $\phi(h)L_K<1$.
The Banach fixed-point theorem applies.
\end{proof}

Both conditions in \eqref{eq:contraction} are needed, and uniqueness is asserted only within $K$.
Two practical versions are the following.

\begin{corollary}\label{cor:wellposed}
\begin{enumerate}[label=(\alph*),leftmargin=2.2em]
  \item \textnormal{(Local ball condition.)}
  Let $\overline B(x_n,r) \subseteq \Omega$ be a closed ball with $\{x_n\}\times\overline B(x_n,r)\subseteq\mathcal A$.
Assume that $\|F(x_n,z)\| \le M$ and that 
$F(x_n,\cdot)$ is Lipschitz with constant $L$ on this ball.
If $\varphi(h)M \le r$ and $\varphi(h)L<1$, then the step has exactly one solution in $\overline B(x_n,r)$.
  \item \textnormal{(Uniform conditions.)} If the same set $K$ satisfies \eqref{eq:contraction} for all $h$ in a range $\mathcal{H}$ and $\sup_{h\in\mathcal{H}}\phi(h)<1/L_K$, the conclusion of Theorem~\ref{thm:wellposed} holds for every $h\in\mathcal{H}$. If $\Omega=\R^n$ and $\mathcal{F}(x_n,\cdot)$ is globally Lipschitz, one may take $K=\R^n$, and the self-map 
  property is automatic.
\end{enumerate}
\end{corollary}

\begin{proof}
 For (a), $\|\Theta(z)-x_n\| = \phi(h) \|\mathcal{F}(x_n,z)\| \le r$ shows $\Theta(\overline{B}(x_n,r))\subseteq\overline{B}(x_n,r)$; apply Theorem~\ref{thm:wellposed}. Part (b) is immediate.
\end{proof}

%
Bounded denominators, such as $\frac{2}{\alpha} \tanh\bigl( \frac{\alpha h}{2} \bigr) < 2/\alpha$, cap the effective step size in the sense of Mickens' renormalized time step \cite{Mickens2002}.
If the cap satisfies the conditions of Corollary~\ref{cor:wellposed}(b), solvability holds uniformly in $h$.
Boundedness alone does not guarantee this.
For $\phi(h)=h$, the contractivity bound $hL<1$ is sufficient, but not necessary, for solvability.

%
In practice \eqref{eq:ph_nsfd} can, for example, be solved by a (damped) Newton iteration, initialized by a suitable predictor and subject to an iteration limit. 
The solvers and initial guesses used in the numerical experiments are specified in Section~\ref{sec:numerics}.
With $g(z)=\dgrad H(x_n,z)$ and $S(z)=\Jb(x_n,z)-\Rb(x_n,z)$, the residual and its Jacobian are
\begin{equation*}
    r(z) = z - x_n - \phi(h) \bigl[ S(z)g(z)+\Gb(x_n,z)u \bigr], \quad
    D_z r(z) = I - \phi(h)\,D_z\bigl[ S(z)g(z)+\Gb(x_n,z)u \bigr],
\end{equation*}
where the full derivative includes the derivatives of $S$, $g$ and $\Gb$.
For a symmetric consistent right-hand side, $D_z\mathcal{F}(x,x) = \frac{1}{2} Df(x)$, which motivates the frozen Jacobian $I-\frac{\phi(h)}{2} Df(x_n)$.
Newton updates are computed by linear solves.
Since all results concern solved steps, a step is accepted only if the residual criterion holds at the returned iterate (re-evaluated after the last update), otherwise a solver failure is recorded.
Outputs and diagnostics are evaluated at the accepted state, where, for an exact discrete gradient, the residual \eqref{eq:residual} equals $g^\top r/\phi(h)$.


We now discuss consistency and order of accuracy, i.e.\ we want to quantify the price and the benefit of a nontrivial denominator function $\phi(h)$.
\begin{definition}\label{ExtendedDF}
The denominator $\phi$ is said to be of \textit{order $q$} if $\phi(h) = h + \mathcal{O}(h^{q+1})$ as $h \to 0^+$.
\end{definition}
Note that \eqref{eq:phi_conditions} only requires $q\ge1$.
We first treat autonomous problems, i.e.\ $u\equiv0$ or $u$ constant.
Time-dependent inputs $u=u(t)$ are considered later in Corollary~\ref{cor:order_input}.

\begin{theorem}[Order of the \pHNSFD{} scheme]\label{thm:order}
Consider \eqref{eq:ph_system} with $u\equiv0$ or $u$ constant on a finite interval $[0,T]$, with a sufficiently smooth vector field $f$ and solution $x(t)$.
Let $\Psi_s$ denote the map of the underlying discrete-gradient scheme ($\phi \equiv \mathrm{id}$) with step size $s$, assume that $\Psi_s$ has order $p\in\{1,2\}$ and is Lipschitz stable on a neighborhood of the trajectory, and let $\phi$ be of order $q \ge 1$.
Then, for $x_n=x(t_n)$, the local error satisfies
\begin{equation}\label{eq:localerror}
    \Psi_{\phi(h)}(x_n) - x(t_{n+1}) = \bigl(\phi(h)-h\bigr) f(x_n) + \mathcal{O}\bigl( h^{p+1}+h^{q+2} \bigr),
\end{equation}
and the global error satisfies
\begin{equation}\label{eq:orderformula}
    \max_{nh\le T} \|x_n-x(t_n)\| = \mathcal{O}\bigl( h^{\min\{p,q\}} \bigr).
\end{equation}
Thus the scheme has order \textit{at least} $\min\{p,q\}$.
In particular, we have the following:
\begin{enumerate}[label=(\alph*),leftmargin=2.2em]
  \item if the underlying scheme is of order $p=2$ (for instance a symmetric discrete gradient and symmetric operators satisfying (DG2)) and $q\ge2$, the scheme has order $2$;
  \item if $p=2$ and $\phi=\phi_{\mathrm{E}}$ is the exponential denominator \eqref{eq:exp_denominator}, $\phi_{\mathrm{E}}(h)=h-\frac{\alpha}{2}h^2+\mathcal{O}(h^3)$, the local error contains the term $-\frac{\alpha}{2}h^2f(x_n)$; the method is therefore of first order along every trajectory on which $f$ does not vanish identically.
\end{enumerate}
\end{theorem}
\begin{proof}
Let $s=\phi(h)$ and $\delta=s-h=\mathcal{O}(h^{q+1})$.
Let $\mathcal{E}_t$ denote the exact flow.
By assumption $\Psi_s(x_n)=\mathcal{E}_s(x_n)+\mathcal{O}(s^{p+1})$, and $s=\mathcal{O}(h)$.
Expanding the exact flow about time instance $h$,
\begin{equation*}
  \mathcal{E}_{s}(x_n)-\mathcal{E}_h(x_n)
  = \int_h^{h+\delta} f\bigl(\mathcal{E}_\tau(x_n)\bigr)\,\dd\tau
  = \delta f\bigl(\mathcal{E}_h(x_n)\bigr) +\mathcal{O}(\delta^2)
  = \delta f(x_n) +\mathcal{O}\bigl(h|\delta|+\delta^2\bigr),
\end{equation*}
and $h|\delta|+\delta^2 = \mathcal{O}(h^{q+2})$.
Together with $\mathcal{E}_h(x_n) = x(t_{n+1})$ this gives \eqref{eq:localerror}.
The local error is therefore $\mathcal{O}(h^{p+1}+h^{q+1})$, and Lipschitz stability together with the discrete Gronwall inequality yields \eqref{eq:orderformula}, which is part (a).
The order $p=2$ of the symmetric underlying scheme follows from its self-adjointness \cite[Thm.~II.3.2]{HairerLubichWanner2006}, see Remark~\ref{rem:symmetry}.
For part (b), $\delta = -\frac{\alpha}{2}h^2 +\mathcal{O}(h^3)$ and $p=2$ give a local error $-\frac{\alpha}{2}h^2f(x_n) +\mathcal{O}(h^3)$, which accumulates to a global error of exact order one unless $f$ vanishes along the trajectory.
\end{proof}

Theorem~\ref{thm:order} gives a lower bound: a denominator may also cancel the leading error term and even yield an exact scheme for linear problems (Section~\ref{sec:denominator}).

\begin{corollary}[Time-dependent inputs] \label{cor:order_input}
Let $H\in C^3$ and $J,R,G\in C^2$ near the trajectory, let $u\in C^2([0,T];\R^m)$, let $u_{\nph}$ be the midpoint value $u(t_n+\frac{h}{2})$ or the cell average $\frac{1}{h}\int_{t_n}^{t_{n+1}} u(s)\,\dd s$.
Let $\dgrad H$ and $\Jb$, $\Rb$, $\Gb$ be symmetric and satisfy (DG2), and let $\phi(h)=h+\mathcal{O}(h^3)$.
If the steps are computed on the solution branch with $x_{n+1}-x_n = \mathcal{O}(h)$ in a compact neighborhood of the trajectory (Theorem~\ref{thm:wellposed}), then the \pHNSFD{} scheme is (globally) convergent of order $2$ on $[0,T]$.
\end{corollary}
\begin{proof}
By (DG2), the right-hand side equals $f\bigl(\frac{x_n+x_{n+1}}{2},t_n+\frac{h}{2}\bigr) +\mathcal{O}(h^2)$, and $\phi(h)-h = \mathcal{O}(h^3)$.
Hence, the step is an $\mathcal{O}(h^3)$ perturbation of the implicit midpoint rule in physical time $h$, which has a local error of order $\mathcal{O}(h^3)$ and is Lipschitz stable for small $h$.
\end{proof}

\begin{remark}[Symmetry]\label{rem:symmetry}
If $\dgrad H$, $\Jb$, $\Rb$, $\Gb$ are symmetric in their two arguments and the step equation selects a unique inverse branch, the underlying discrete-gradient map satisfies $\Psi_{-s} = \Psi_s^{-1}$, i.e.\ it is self-adjoint and hence of even order.
The complete method $x_n\mapsto\Psi_{\phi(h)}(x_n)$ is \textit{self-adjoint} if, in addition, $\phi$ is odd.
Together with sufficient smoothness, oddness implies $\phi(h)=h+\mathcal{O}(h^3)$.
%
Oddness alone is not sufficient (e.g.\ $\phi(h)=h+h|h|$).
When a nonodd denominator such as $\phi_{\mathrm{E}}$ is combined with a symmetric discrete gradient, only the underlying discrete-gradient method is \textit{symmetric}, the complete method is not.
\end{remark}

The functions $h \mapsto \frac{2}{\omega}\tan\bigl(\frac{\omega h}{2}\bigr)$ and $h \mapsto \frac{2}{\alpha} \tanh\bigl(\frac{\alpha h}{2}\bigr)$ are odd, analytic near $h=0$ and of order $q=2$. 
The former is positive on the first admissible interval $0 < h < \pi/\omega$, the only one used here. 
The function $h \mapsto (1-\mathrm{e}^{-\alpha h})/\alpha$ is neither odd nor of order $q\ge2$.

Dang and Hoang~\cite{DangHoang2020} introduced a class of bounded denominator functions satisfying $\phi(h) = h +\mathcal{O}(h^{q+1})$ as $h \to 0^+$ for arbitrary $q$, which ensure that the resulting nonstandard explicit Runge--Kutta (ERK) methods preserve the order of accuracy of the underlying standard Runge--Kutta methods, see also \cite{HoangMatthias2026NSFD} for second-order nonstandard methods for non-autonomous systems.
These denominator functions are, however, not odd in general.

A family that is positive for every $h>0$, bounded, odd and of arbitrarily high order is given, for $m\in\mathbb{N}$ and $\tau>0$, by
\begin{equation}\label{eq:phi_family_bounded}
\begin{split}
    \phi_m(h) &= \frac{h}{\bigl(1 + \tau h^{2m}\bigr)^{1/(2m)}},\qquad
    0<\phi_m(h) < \tau^{-1/(2m)}\quad (h>0),\\
    \phi_m(h) &= h - \frac{\tau}{2 m}h^{2 m+1}+\mathcal{O}\bigl(h^{4 m+1}\bigr).
\end{split}
\end{equation}
Its cap $\tau^{-1/(2m)}$ can be properly chosen in order to satisfy the contraction and invariant-domain conditions of Corollary~\ref{cor:wellposed}(b), or the cap required in Theorem~\ref{thm:pendulum} below.


\section{Discrete Structural Properties of the \pHNSFD{} Scheme} \label{sec:properties}
\subsection{Discrete Power Balance} \label{sec:powerbalance}
The principal structural property of the proposed \pHNSFD{} method follows directly from the discrete-gradient identity \eqref{eq:discrete_gradient}.
\begin{theorem}[Exact discrete power balance] \label{thm:powerbalance}
Let \eqref{eq:ph_nsfd} be a \pHNSFD{} scheme in the sense of Definition~\ref{def:phnsfd}.  
Then for every $n$ and every solved step with $\phi(h)>0$,
\begin{equation}\label{eq:discrete_power_balance}
    H(x_{n+1})-H(x_n)
    = - \phi(h)\,\dgrad H_n^\top \Rb_n \dgrad H_n
      + \phi(h)\, y_{\nph}^\top u_{\nph}
    = - \phi(h)\,\mathcal{D}_n + \phi(h)\,y_{\nph}^\top u_{\nph},
\end{equation}
with the \textit{discrete dissipation rate} $\mathcal{D}_n := \dgrad H_n^\top \Rb_n \dgrad H_n \ge 0$.
\end{theorem}
\begin{proof}
By the defining property \eqref{eq:discrete_gradient} of the discrete gradient, we have
\begin{equation*}
    H(x_{n+1}) - H(x_n) = \dgrad H_n^\top (x_{n+1}-x_n).
\end{equation*}
Inserting \eqref{eq:ph_nsfd}, yields
\begin{equation*}
    \begin{split}
    H(x_{n+1})-H(x_n)
    &=    \phi(h)\,\dgrad H_n^\top \Jb_n \dgrad H_n
    -    \phi(h)\,\dgrad H_n^\top\Rb_n\dgrad H_n
    +    \phi(h)\,\dgrad H_n^\top \Gb_n u_{\nph}.
 \end{split}
\end{equation*}
Since $\Jb_n$ is skew-symmetric, $\dgrad H_n^\top\Jb_n\dgrad H_n =\tfrac{1}{2}\dgrad H_n^\top(\Jb_n+\Jb_n^\top)\dgrad H_n=0$.
Moreover $\dgrad H_n^\top\Gb_n=y_{\nph}^\top$ by \eqref{eq:ph_nsfd}, and $\mathcal{D}_n\ge0$ because $\Rb_n\succeq0$.
\end{proof}

The identity \eqref{eq:discrete_power_balance} is algebraic:
it holds for every solved
step with $\varphi(h)>0$ and does not require the condition (DG2).

\begin{corollary}[Discrete dissipation]\label{cor:dissipation}
For a closed \pHNSFD{} system, $u_{\nph}=0$, we have
\begin{equation}  \label{eq:discrete_dissipation}
    H(x_{n+1})-H(x_n) = -\phi(h)\,\dgrad H_n^\top\Rb_n\dgrad H_n \le0 ,
\end{equation}
hence, the 
discrete energy is monotonically non-increasing for every solved step.
\end{corollary}

\begin{corollary}[Exact conservation]\label{cor:conservation}
If, in addition, $\Rb_n=0$, then $H(x_{n+1})=H(x_n)$ for every solved step.
Hence, $H(x_n)=H(x_0)$ for all $n$.
\end{corollary}
Thus, the \pHNSFD{}  method preserves the conservative Hamiltonian exactly.
%

The exact discrete power balance 
\eqref{eq:discrete_power_balance} can also be interpreted as an approximation of the integrated continuous balance. However, this interpretation requires additional assumptions.

\begin{proposition}[Consistency of the discrete dissipation rate] \label{prop:dissipation_consistency}
Let $x(t)$ be a sufficiently smooth solution, let $\dgrad H$ and $\Rb$ satisfy (DG2) on a compact neighborhood of the trajectory, and let $\mathcal{D}$ be locally Lipschitz.
Evaluate $\mathcal{D}_n$ at the exact endpoints $x(t_n)$, $x(t_{n+1})$, or at numerical states whose global error is $\mathcal{O}(h^2)$.
Then, the discrete dissipation $\mathcal{D}_n$  satisfies
\begin{equation}\label{eq:dissipation_midpoint}
    \mathcal{D}_n = \mathcal{D}\Bigl( x\bigl( t_n+\tfrac{h}{2} \bigr) \Bigr) + \mathcal{O}(h^2).
\end{equation}
If, in addition, $\phi(h) = h+\mathcal{O}(h^3)$ and $t \mapsto \mathcal{D}(x(t))$ is of class $C^2$ (e.g.\ $H\in C^3$, $R\in C^2$), then 
\begin{equation}\label{eq:dissipation_quadrature}
    \phi(h)\,\mathcal{D}_n = \int_{t_n}^{t_{n+1}} \mathcal{D}\bigl(x(t)\bigr)\,\dd t +\mathcal{O}(h^3).
\end{equation}
Thus, the dissipative term in \eqref{eq:discrete_power_balance} is a second-order accurate approximation of the corresponding integral in \eqref{eq:cont_pb}.
\end{proposition}
\begin{proof}
Since $\frac{1}{2}\bigl(x(t_n)+x(t_{n+1})\bigr)=x\bigl(t_n+\frac{h}{2}\bigr)+\mathcal{O}(h^2)$, (DG2) and Lipschitz continuity give \eqref{eq:dissipation_midpoint}.
Numerical states with global error $\mathcal{O}(h^2)$ change $\mathcal{D}_n$ by $\mathcal{O}(h^2)$ only.
The midpoint quadrature rule and $\phi(h)=h+\mathcal{O}(h^3)$ yield \eqref{eq:dissipation_quadrature}.
\end{proof}

If instead $\varphi(h) = h + ch^2 +O(h^3)$ with $c\ne0$, then
we get \begin{equation*}
  \varphi(h)\mathcal D_n = \int_{t_n}^{t_{n+1}}\mathcal D(x(t))\,dt
  + ch^2\mathcal D\bigl(x(t_n+\tfrac{h}{2})\bigr) +O(h^3).
\end{equation*}
Thus, this approximation of the dissipated energy is generally only first-order accurate.

If $H$ is coercive and $u\equiv0$, Corollary~\ref{cor:dissipation} implies that $x_n$ remains in the compact set $\{H\le H(x_0)\}$, i.e.\ the numerical solution is bounded for every step size for which the steps are solved.

\subsection{Discrete Passivity and Stability}\label{sec:passivity}
\begin{theorem}[Unconditional discrete passivity]\label{thm:passivity}
Assume $\Rb_n\succeq0$.
Then the \pHNSFD{} scheme satisfies the discrete passivity inequality
\begin{equation} \label{eq:discrete_passivity}
    H(x_{n+1})-H(x_n) \le \phi(h)\, y_{\nph}^\top u_{\nph}
\end{equation}
and, by summation,
\begin{equation}\label{eq:discrete_passivity_integrated}
    H(x_N)-H(x_0) \le  \phi(h)\sum_{n=0}^{N-1} y_{\nph}^\top u_{\nph}
    \qquad\text{for all } N\in\mathbb{N} .
\end{equation}
If $H$ is bounded below, $H\ge0$ say, the discrete system with input
$u_{\nph}$, output $y_{\nph}$ and storage function $H$ is passive for every
solved step.
\end{theorem}

The proof is immediate from Theorem~\ref{thm:powerbalance}. 
Passivity follows from the algebraic properties of the discrete operators for every solved step. 
In contrast, for $\dot{x}=JQx$ with $J^\top = -J$ and $Q\succ0$, explicit Euler increases the quadratic energy at every non-equilibrium state.
Algebraic stability of a Runge--Kutta method does not replace the discrete chain rule for nonquadratic $H$: for $H=\frac{1}{4}q^4+\frac{1}{2}p^2$ and $h=1$, the implicit midpoint step $(\frac{1}{2},\frac{3}{2})\mapsto(\frac{3}{2},\frac{1}{2})$ increases $H$ by $\frac{1}{4}$ without any supply. 
Estimate \eqref{eq:discrete_passivity_integrated} is the exact discrete analogue of \eqref{eq:cont_passivity}.
For a state-dependent denominator $\phi_n=\phi(h;x_n)$ the supply sum carries the weights $\phi_n$.
The supply in \eqref{eq:discrete_passivity} is weighted by $\phi(h)$ (in physical time, the output is $\frac{\phi(h)}{h}y_{\nph}$): 
coupling two lossless subsystems with $\phi_1\neq\phi_2$ by $u_1=-y_2$, $u_2=y_1$ gives $\Delta(H_1+H_2)=(\phi_2-\phi_1)y_1^\top y_2$, so interconnected components should share a denominator or be coupled through the scaled outputs.

%
Analogously, let $(x^\ast,u^\ast)$ be a forced equilibrium, $H$ convex, and $\Jb=J$, $\Rb=R$, $\Gb=G$ constant (constancy of $J$, $R$, $G$ alone does not suffice). 
With $\tilde{g} = \dgrad H_n-\nabla H(x^\ast)$, the step reads $x_{n+1}-x_n = \phi(h)\bigl[(J-R)\tilde{g}+G(u_{\nph}-u^\ast)\bigr]$ and $y_{\nph}-y^\ast = G^\top \tilde{g}$, so the shifted supply rate $(y_{\nph}-y^\ast)^\top (u_{\nph}-u^\ast)$ satisfies \eqref{eq:discrete_passivity} with the \emph{Bregman divergence} $H(x)-H(x^\ast)-\nabla H(x^\ast)^\top(x-x^\ast)\ge0$ as storage function.

\begin{theorem}[Discrete Lyapunov stability]\label{thm:lyap}
Let $u\equiv0$, let $x^\ast$ be a strict local minimum of $H$ with $\nabla H(x^\ast)=0$, and assume Assumption~\ref{ass:main}.  
Let $h$ be such that the \pHNSFD{} step defines a continuous map $\Psi_h$ on a neighborhood $U$ of $x^\ast$ with $\Psi_h(x^\ast)=x^\ast$.
Then $x^\ast$ is a stable fixed point of $\Psi_h$.
\end{theorem}

\begin{proof}
By Corollary~\ref{cor:dissipation}, $V = H-H(x^\ast)$ is non-increasing along the discrete flow.
For small $c>0$, the connected component $U_c$ of $\{V\le c\}$ containing $x^\ast$ is a compact subset of $U$ that shrinks to $\{x^\ast\}$ as $c\to0$, and $\Psi_h(U_c)$ is connected, contains $x^\ast$ and lies in $\{V\le c\}$.
Hence, $\Psi_h(U_c)\subseteq U_c$.
\end{proof}

The continuity assumption is not a technicality: since $V$ is only non-increasing, a discontinuous selection of solutions of the implicit step could jump to another component of the sublevel set, see Section~\ref{ssec:model_pendulum} for an example.

\begin{theorem}[Discrete LaSalle] \label{thm:lasalle}
In addition to the hypotheses of Theorem~\ref{thm:lyap}, assume that $H$ is coercive, let $c>H(x^\ast)$, let $K$ be the connected component of $\{H\le c\}$ containing $x^\ast$, which is then compact and contains a neighborhood of $x^\ast$, and assume that $\Psi_h$ is defined and continuous on $K$. 
Define $\mathcal{D}_h(x) := \dgrad H(x,\Psi_h(x))^\top \Rb(x,\Psi_h(x))\, \dgrad H(x,\Psi_h(x))$ and assume the \textit{zero-state detectability} condition that the largest $\Psi_h$-invariant subset of $\{x\in K\colon \mathcal{D}_h(x) = 0\}$ is $\{x^\ast\}$.
Then $x^\ast$ is asymptotically stable for $\Psi_h$, and its basin of attraction contains $K$.
\end{theorem}
\begin{proof}
As in the proof of Theorem~\ref{thm:lyap}, $\Psi_h(K)\subseteq K$.
Since $H(x_{n+1})-H(x_n)=-\phi(h)\mathcal{D}_h(x_n)$ and $H$ is bounded below on $K$, $\mathcal{D}_h(x_n)\to0$, and LaSalle's invariance principle \cite[Ch.~4]{LaSalle1976} together with detectability gives the claim.
\end{proof}

%

Theorem~\ref{thm:lyap} is a nonlinear but qualitative statement.  
For the \textit{linear stability type} of the equilibria a sharper result is available.

\begin{theorem}[Linear stability preservation at nondegenerate energy minima]\label{thm:elemstab}
Let $u\equiv0$, let $x^\ast$ satisfy $\nabla H(x^\ast)=0$ and $Q_\ast := \nabla^2H(x^\ast)\succ0$, and let
$A_\ast = \bigl(J(x^\ast)-R(x^\ast)\bigr)Q_\ast$.
We consider a \pHNSFD{} scheme with a symmetric discrete gradient satisfying (DG2) and admissible operators $\Jb$, $\Rb$ such that the two-argument right-hand side is of class $C^1$ near $(x^\ast,x^\ast)$.
Then, for every $h$ in the domain of $\phi$ with $s:=\phi(h)>0$, the step equation possesses a unique local solution branch $\Psi_h$ near $x^\ast$ with $\Psi_h(x^\ast)=x^\ast$, and
\begin{equation}\label{eq:cayley}
    D\Psi_h(x^\ast)=\mathcal{C}_{s}(A_\ast)
    :=\Bigl(I-\tfrac{s}{2}A_\ast\Bigr)^{-1}
      \Bigl(I+\tfrac{s}{2}A_\ast\Bigr).
\end{equation}
Moreover, we have the following:
\begin{enumerate}[label=(\roman*),leftmargin=2.2em]
  \item an eigenvalue $\lambda$ of $A_\ast$ with $\Real\lambda<0$ is mapped to an eigenvalue of $\mathcal{C}_s(A_\ast)$ of modulus $<1$, and an eigenvalue with $\Real\lambda=0$ to one of modulus $1$;
  \item the eigenvalues of $A_\ast$ on the imaginary axis are semisimple, and so are the corresponding eigenvalues of $\mathcal{C}_s(A_\ast)$ on the unit circle;
  \item $\mathcal{C}_s(A_\ast)$ is contractive in the $Q_\ast$-norm,
  \begin{equation}\label{eq:cayley_contractive}
      \mathcal{C}_s(A_\ast)^\top Q_\ast\,\mathcal{C}_s(A_\ast)-Q_\ast
      = -2s\Bigl(I-\tfrac{s}{2}A_\ast\Bigr)^{-\top} Q_\ast R(x^\ast)Q_\ast \Bigl(I-\tfrac{s}{2}A_\ast\Bigr)^{-1}\preceq0 .
  \end{equation}
\end{enumerate}
Consequently, the linear stability type of $x^\ast$ for the flow and for $\Psi_h$ coincide, for every admissible positive denominator value.
\end{theorem}
\begin{proof}
By property (DG2), $\dgrad H(x^\ast+\xi_n,x^\ast+\xi_{n+1}) = Q_\ast\frac{\xi_n+\xi_{n+1}}{2}+o(\|\xi_n\|+\|\xi_{n+1}\|)$, and the derivatives of $\Jb$, $\Rb$ multiply $\nabla H(x^\ast) = 0$.
Hence, the step residual has the partial derivatives $I-\frac{s}{2}A_\ast$ and $-I-\frac{s}{2}A_\ast$ at $(x^\ast,x^\ast)$.
The former is invertible by Lemma~\ref{lem:spectrum}, and the implicit function theorem gives the local branch and \eqref{eq:cayley}.
For an eigenvalue $\lambda$ of $A_\ast$, the corresponding eigenvalue $r_s(\lambda) = (1+s\lambda/2)/(1-s\lambda/2)$ satisfies
\begin{equation*}
    |r_s(\lambda)|^2-1 = \frac{2s\Real\lambda}{|1-s\lambda/2|^2},
\end{equation*}
which proves (i).
Since $\mathrm{e}^{tA_\ast}$ is contractive in the $Q_\ast$-norm by Lemma~\ref{lem:contract}, it is bounded for $t\ge0$, so the eigenvalues of $A_\ast$ on the imaginary axis are semisimple. 
As $r_s$ is analytic near $\spec(A_\ast)$ with $r_s'(\lambda)=s/(1-s\lambda/2)^2 \neq 0$, $\mathcal{C}_s(A_\ast) = r_s(A_\ast)$ has Jordan blocks of the same sizes, which proves (ii).
For (iii), let $W=I-\frac{s}{2}A_\ast$ and $V = I+\frac{s}{2}A_\ast$,
then $V^\top Q_\ast V - W^\top Q_\ast W = s(A_\ast^\top Q_\ast+Q_\ast A_\ast) = -2 sQ_\ast R(x^\ast)Q_\ast$, and congruence with $W^{-1}$ gives \eqref{eq:cayley_contractive}.
\end{proof}

Theorem~\ref{thm:elemstab} concerns the local branch near $x^\ast$ and is not a global solvability statement.
Combined with Proposition~\ref{prop:equilibrium} below and Theorem~\ref{thm:wellposed}, it yields elementary stability for systems whose equilibria are nondegenerate energy minima, on every range of step sizes on which the step map is well defined and its branch retains the diagonal solutions. 
For equilibria that are not energy minima, the same Cayley formula holds with the full Jacobian $A_\ast = Df(x^\ast)$ for symmetric consistent right-hand sides, provided $I-\frac{s}{2}A_\ast$ is invertible.
This is guaranteed for every $h$ by a bounded denominator with $\sup_h\phi(h)<2/\Lambda$ if the positive real eigenvalues of $A_\ast$ are bounded by $\Lambda$.
\subsection{Casimir Preservation}\label{sec:casimir}
Suppose the continuous system possesses a Casimir $C$ as in
Definition~\ref{def:casimir}.
Exact preservation of $C$ by the scheme is only meaningful if $C$ is an invariant of the full system, i.e.\ if in addition $\nabla C^\top R=0$ and $\nabla C^\top G=0$.
If only $\nabla C^\top J=0$ holds, the appropriate discrete statement is the Casimir balance \eqref{eq:casimirbalance} below.
Let $\dgrad C$ be a discrete gradient of $C$, then
\begin{equation}\label{eq:dgradC}
    C(x_{n+1})-C(x_n) = \dgrad C_n^\top (x_{n+1}-x_n), \qquad 
    \dgrad C(x,x)=\nabla C(x).
\end{equation}

\begin{definition}[Discrete Casimir conditions]\label{def:casimirconditions}
We say that the \pHNSFD{} scheme satisfies the \textit{discrete Casimir
conditions} for $C$ if for all $x,z\in\Omega$
\begin{equation}\label{eq:casimirconditions}
    \dgrad C(x,z)^\top\,\Jb(x,z)=0,\qquad
    \dgrad C(x,z)^\top\,\Rb(x,z)=0,\qquad
    \dgrad C(x,z)^\top\,\Gb(x,z)=0 .
\end{equation}
\end{definition}

\begin{theorem}[Exact Casimir preservation]\label{thm:casimir}
Under \eqref{eq:casimirconditions}, the \pHNSFD{} scheme preserves $C$ exactly:
\begin{equation}\label{eq:44new}
    C(x_{n+1})=C(x_n)\qquad \qquad\text{for every solved step with }\varphi(h)>0.
\end{equation}
If only the first condition in \eqref{eq:casimirconditions} holds, one obtains the exact \emph{discrete Casimir balance}
\begin{equation}\label{eq:casimirbalance}
    C(x_{n+1})-C(x_n)
    = \phi(h)\Bigl(-\dgrad C_n^\top\Rb_n\dgrad H_n
       +\dgrad C_n^\top\Gb_nu_{\nph}\Bigr),
\end{equation}
which is the discrete analogue of $\dot{C} = -\nabla C^\top R\nabla H + \nabla C^\top Gu$.
\end{theorem}
\begin{proof}
By \eqref{eq:dgradC} and \eqref{eq:ph_nsfd},
\begin{equation*}
    C(x_{n+1})-C(x_n) = \phi(h)\, \dgrad C_n^\top
    \bigl[ (\Jb_n-\Rb_n) \dgrad H_n +\Gb_nu_{\nph} \bigr].
\end{equation*}
The three conditions in \eqref{eq:casimirconditions} annihilate the three terms
and give \eqref{eq:44new}. 
If only the first condition holds, the $J$-term vanishes and the remaining terms give \eqref{eq:casimirbalance}.
\end{proof}

The discrete Casimir conditions \eqref{eq:casimirconditions} couples the choice of $\dgrad C$ to the choice of $\Jb$, $\Rb$, $\Gb$. It is not automatic.  
The following constructive routes are available.
\begin{enumerate}[label=(\alph*),leftmargin=2.2em]
  \item \textit{Affine Casimirs.}  
  If $C(x)=c^\top x$, we choose $\dgrad C\equiv c$ (which is the AVF, Gonzalez and Itoh--Abe gradient of $C$).
    Then \eqref{eq:casimirconditions} reduces to the pointwise linear conditions
    $c^\top\Jb(x,z)=0$, $c^\top\Rb(x,z)=0$, $c^\top\Gb(x,z)=0$, which hold for the midpoint or symmetrized operators \eqref{eq:midpointops} whenever $c^\top J=c^\top R=0$ and $c^\top G=0$ on $\Omega$.
    Most Casimirs arising from Kirchhoff-type interconnection constraints are of this form.
  \item \textit{Projection on a regular region.}  For nonlinear $C$ and a symmetric discrete gradient $\dgrad C$ with $\|\dgrad C\|\ge\eta>0$, the projected operators $\Pi\Jb\Pi$, $\Pi\Rb\Pi$, $\Pi\Gb$ with $\Pi=I-\dgrad C\,\dgrad C^\top/(\dgrad C^\top\dgrad C)$ satisfy \eqref{eq:casimirconditions} and retain skew-symmetry and positive semidefiniteness.
    If the continuous kernel conditions hold, they also retain the diagonal consistency \eqref{eq:discrete_consistency}. 
    If (DG2) holds for $\dgrad C$ \textit{and} for $\Jb$, $\Rb$, $\Gb$, the projected operators satisfy (DG2), and second order is retained under the assumptions of Corollary~\ref{cor:order_input}. 
    At critical points of $C$ the projector is singular.
    One then works on a regular region or uses the structure of the model, as in Example~\ref{ex:rigid_discrete}.
  \item \textit{Simultaneous conservation.}  Alternatively one may use the
    discrete gradients of \cite{DahlbyOwrenYaguchi2011,NortonQuispel2014},
    which are designed to conserve several first integrals at once.
    In the pH setting the pair $(H,C)$ is exactly of this type.
\end{enumerate}

\begin{example}[Rigid body]\label{ex:rigid_discrete}
For Example~\ref{ex:rigidbody} with $C(x)=\frac{1}{2}\|x\|^2$, set $m=\frac{x+z}{2}$ and choose
\begin{equation*}
    \dgrad C(x,z)=m,\qquad
    \Jb(x,z)=\widehat{m},\qquad
    \Rb(x,z)=\rho\,\widehat{m}\,\widehat{m}^\top,\qquad
    \Gb=0,\qquad \rho\ge0,
\end{equation*}
where $\widehat{m}$ is the cross-product matrix of $m$, and let $\dgrad H(x,z)=\diag(I_1^{-1},I_2^{-1},I_3^{-1})\,m$ be the midpoint gradient of the diagonal quadratic Hamiltonian (Lemma~\ref{lem:quadratic}).
Since $m^\top\widehat{m}=0$, both $\Jb$ and $\Rb$ annihilate $\dgrad C$, so \eqref{eq:casimirconditions} holds.
%
By Theorems~\ref{thm:casimir} and~\ref{thm:powerbalance},
\begin{equation*}
    C(x_{n+1})=C(x_n),
    \qquad
    H(x_{n+1})-H(x_n)=-\phi(h)\,\rho\,\bigl\|\widehat{m}^\top\dgrad H_n\bigr\|^2\le0 .
\end{equation*}
For $\rho=0$ both $H$ and $C$ are conserved for every solved step and every $\phi$.
For $\rho>0$, $C$ is preserved exactly, whereas $H$ decreases.
\end{example}

\begin{remark}[A discrete Dirac structure]\label{rem:dirac}
For each step, the graph
\begin{equation*}
    \mathcal{D}_n = \bigl\{ ((f_s,y),(e,u))\colon\,f_s = -\Jb_n e -\Gb_n u,\, y = \Gb_n^\top e\bigr\}
\end{equation*}
of a skew-symmetric linear map is isotropic, $e^\top f_s + u^\top y = 0$, and has dimension $n+m$, i.e.\ it is a \emph{(discrete) Dirac structure} \cite{KotyczkaLefevre2019, TalasilaEtAl2006}.
With $e=\dgrad H_n$ the step reads $(x_{n+1}-x_n) / \phi(h) = -f_s -\Rb_n \dgrad H_n$, so that \eqref{eq:discrete_power_balance} is the power pairing of $\mathcal{D}_n$ with $\Rb_n$ attached as a resistive relation.
The Casimir conditions are additional kernel conditions.
\end{remark}

\subsection{Equilibrium Preservation}\label{sec:equilibrium}
\begin{proposition}[Exact equilibria]\label{prop:equilibrium}
Let $(x^\ast,u^\ast)$ satisfy \eqref{eq:equilibrium} and let the discrete operators satisfy the consistency conditions \eqref{eq:discrete_consistency}.  
Then $x_n = x_{n+1} = x^{\ast}$, $u_{\nph} = u^\ast$ solves 
\eqref{eq:ph_nsfd} for every $h$ with $\phi(h)>0$, that is, every equilibrium of the continuous pH system is a diagonal solution of the implicit step relation.
Conversely, if $\phi(h)>0$ and $\Jb$, $\Rb$, $\Gb$ are consistent, every fixed point of a \pHNSFD{} map with constant input $u^\ast$ is an equilibrium of \eqref{eq:ph_system}.
\end{proposition}
\begin{proof}
Setting $x_n = x_{n+1} = x^\ast$ gives $\dgrad H(x^\ast,x^\ast) = \nabla H(x^\ast)$ by \eqref{eq:consistency_gradient} and $\Jb = J(x^\ast)$, $\Rb = R(x^\ast)$, $\Gb=G(x^\ast)$ by \eqref{eq:discrete_consistency}, so the right-hand side of \eqref{eq:ph_nsfd} equals the left-hand side of \eqref{eq:equilibrium}, which vanishes.
The left-hand side of \eqref{eq:ph_nsfd} vanishes trivially.
For the converse, a fixed point $x_{n+1} = x_n = \bar{x}$ forces the right-hand side of \eqref{eq:ph_nsfd} to vanish, and by consistency this is $\bigl( J(\bar{x})-R(\bar{x}) \bigr) \nabla H(\bar{x}) +G(\bar{x})u^\ast = 0$.
\end{proof}

This holds for every denominator function $\phi(h)$, but relies on the exact diagonal consistency of the scalar class.
For the matrix extensions of Section~\ref{sec:locallyexact} it must be verified separately.
The Euler methods also preserve the location of equilibria, but not their linear stability type.
If the step equation has several solutions, a selected step map need not choose the diagonal one: 
for $H(x) = \frac{1}{4} (x^2-1)^2$, $\dot{x} = -H'(x)$, the AVF step from $x_n = 0$ with $\phi(h) = 4$ reduces to $z^3-z = 0$, and the smooth branch through the regular solution $(0,1)$ maps the equilibrium $0$ to $1$.
Equality of fixed points and equilibria thus requires a branch retaining the diagonal solutions, e.g.\ by uniqueness (Theorem~\ref{thm:wellposed}) or as in Theorem~\ref{thm:elemstab}.


\section{Construction of Nonstandard Denominator Functions}\label{sec:denominator}
The denominator function is the principal additional degree of freedom introduced by the NSFD framework.  
By Section~\ref{sec:properties} it costs nothing structurally.
By Theorem~\ref{thm:order} it must be chosen with care if the order is to be preserved.  
%
%
%
%
We now give a systematic construction, for which we require the denominator function $\phi$ to be positive, of order $q\ge p$ for an underlying scheme of order $p$, preferably bounded (Corollary~\ref{cor:wellposed}(b)), and \textit{spectrally matched}, i.e.\ reproducing a designated characteristic time scale exactly.

\subsection{The midpoint scheme and its exact denominator}
Consider the scalar linear test problem $\dot{\xi} = \lambda\xi$.
By Theorem~\ref{thm:elemstab} the \pHNSFD{} scheme \eqref{eq:ph_nsfd} acts on it as
\begin{equation}\label{eq:scalarcayley}
    \xi_{n+1}=\frac{1+\frac{1}{2}\phi(h)\lambda}{1-\frac{1}{2}\phi(h)\lambda}\,\xi_n .
\end{equation}
Requiring exactness, $\xi_{n+1}=\mathrm{e}^{\lambda h}\xi_n$, and solving for $\phi$
gives, with $z = \phi(h)\lambda/2$,
\begin{equation*}
\frac{1+z}{1-z} = \mathrm{e}^{\lambda h}
 \quad \Longleftrightarrow\quad
 z = \frac{\mathrm{e}^{\lambda h}-1}{\mathrm{e}^{\lambda h}+1} = \tanh\Bigl(\frac{\lambda h}{2}\Bigr),
 \end{equation*}
i.e.\ we have the following result.

\begin{proposition}[Exact scalar denominator for the midpoint scheme]
\label{prop:exactscalar}
Let $\lambda\in\C\setminus\{0\}$ and step-size $h>0$ with $\mathrm{e}^{\lambda h}\neq-1$.
For $\dot{\xi} = \lambda\xi$, the \pHNSFD{} scheme is exact if and only if
\begin{equation}\label{eq:exact_scalar_phi}
    \phi(h)=\phi_\lambda(h):=\frac{2}{\lambda}\tanh\Bigl(\frac{\lambda h}{2}\Bigr).
\end{equation}
Moreover $\phi_\lambda$ is odd, of order $q=2$, with the expansion
\begin{equation}\label{eq:phiexpansion}
    \phi_\lambda(h)=h-\frac{\lambda^2 h^3}{12}
    +\frac{\lambda^4 h^5}{120}-\cdots ,
\end{equation}
and it is real whenever $\lambda\in\R$ or $\lambda\in\iu\R$.
\end{proposition}

Formula \eqref{eq:exact_scalar_phi} is a formal relation for the test equation $\dot{\xi} = \lambda\xi$.
It defines a member of the real scalar class only on step-size ranges where it is real, finite and positive. 

Two important special cases of \eqref{eq:exact_scalar_phi} are obtained by specializing $\lambda$.

\subsubsection*{Dissipative modes}
For a real decaying mode $\lambda=-\alpha$, $\alpha>0$, we obtain
\begin{equation}\label{eq:tanh_denominator}
    \phi(h) = \frac{2}{\alpha}\tanh\Bigl(\frac{\alpha h}{2}\Bigr)
    = h-\frac{\alpha^2 h^3}{12}+\mathcal{O}(h^5),\quad 
    \text{with}\quad 0<\phi(h)<\frac{2}{\alpha}.
\end{equation}
This should be compared with the classical Mickens \textit{exponential denominator function}
\begin{equation} \label{eq:exp_denominator}
    \phi_{\mathrm{E}}(h) =  \frac{1-\mathrm{e}^{-\alpha h}}{\alpha}
    = h-\frac{\alpha}{2}h^2+\mathcal{O}(h^3),\quad 
    \text{with}\quad 0<\phi_{\mathrm{E}}(h)<\frac{1}{\alpha},
\end{equation}
which is exact for the \textit{forward Euler scheme} \eqref{eq:exact_euler}.

\begin{proposition}[$\tanh$ versus $\exp$]\label{prop:tanhvsexp}
Let the underlying discrete-gradient scheme be symmetric ($p=2$).  Then:
\begin{enumerate}[label=(\alph*),leftmargin=2.2em]
  \item with $\phi=\phi_{\mathrm{E}}$ from \eqref{eq:exp_denominator}, the \pHNSFD{} scheme is generally of order one, i.e.\ the exponential denominator \textit{reduces} the order of the method from two to one;
  
  \item with $\phi$ from \eqref{eq:tanh_denominator}, the scheme retains order
    $2$ and is in addition \textit{exact} for every mode with eigenvalue $-\alpha$;
  
  \item both denominators are positive and bounded, with caps $1/\alpha$ and $2/\alpha$. 
  The $\tanh$ denominator is the larger of the two, $\phi_{\mathrm{E}}(h)<\phi(h)$ for $h>0$, and therefore closer to $h$. 
  Indeed, for $\alpha>0$ and $h>0$, we have $0<\varphi_E(h)<\varphi(h)<h$, where the last inequality follows from $\tanh s<s$ for $s>0$.

\end{enumerate}
\end{proposition}
\begin{proof}
(a) and (b) follow from Theorem~\ref{thm:order} together with the Taylor expansions in
\eqref{eq:tanh_denominator}, \eqref{eq:exp_denominator}, and from
Proposition~\ref{prop:exactscalar}.  
For (c), with $s=\alpha h/2>0$, $\alpha\phi_{\mathrm{E}}=1-\mathrm{e}^{-2s}=(\mathrm{e}^{2s}-1)/\mathrm{e}^{2s}$ and $\alpha\phi=2\tanh s=2(\mathrm{e}^{2s}-1)/(\mathrm{e}^{2s}+1)$.
Since $\mathrm{e}^{2s}-1>0$, the inequality $\alpha\phi_{\mathrm{E}}<\alpha\phi$ is equivalent to $\mathrm{e}^{2s}+1<2\mathrm{e}^{2s}$, i.e.\ to $\mathrm{e}^{2s}>1$, which holds.
\end{proof}

A fitted rate $\alpha$ fixes the cap. 
Enlarging $\alpha$ to enforce Corollary~\ref{cor:wellposed}(b) makes it a design parameter and sacrifices exactness, so the contraction conditions must be checked for the actual choice.
For another mode $-\beta$, the factor $(1-\beta\phi/2)/(1+\beta\phi/2)$ tends to $(\alpha-\beta)/(\alpha+\beta)$ as $h\to\infty$, e.g.\ $-99/101$ for $\beta=100\alpha$: fast modes keep their stability type but are only weakly damped.


\subsubsection*{Oscillatory modes}
For a purely imaginary pair $\lambda=\pm\iu\omega$, formula \eqref{eq:exact_scalar_phi} yields a real denominator $\phi(h)\in\R$, because $\tanh(\iu\theta) = \iu\tan\theta$:
\begin{equation} \label{eq:oscillatory_denominator}
    \phi(h) =  \frac{2}{\omega}\tan\Bigl(\frac{\omega h}{2}\Bigr)
    = h+\frac{\omega^2h^3}{12}+\mathcal{O}(h^5),
    \qquad 0<h<\frac{\pi}{\omega}.
\end{equation}
With this choice \eqref{eq:oscillatory_denominator} the midpoint-type discrete-gradient scheme reproduces the exact phase rotation of the harmonic oscillator: the amplification factor is $\mathrm{e}^{2\iu\arctan(\phi\omega/2)}=\mathrm{e}^{\iu\omega h}$.
We use it on the first admissible interval $0<h<\pi/\omega$, adjacent to $h=0$, on which it is exact for the modes $\pm\iu\omega$ (it is positive again on later intervals such as $(2\pi/\omega,3\pi/\omega)$, which are not used).
This restriction, i.e.\ more than two steps per period $2\pi/\omega$, is a \textit{resolution} requirement, not a stability restriction.
Note that the denominator $\phi$ is unbounded as $h\uparrow\pi/\omega$. 
This is an acceptable trade in the conservative case with quadratic $H$, where the implicit system is linear.
For nonlinear problems, the step-size range must be compatible with the solvability conditions of Theorem~\ref{thm:wellposed} (see Theorem~\ref{thm:pendulum} for the pendulum).

\subsubsection*{Damped oscillatory modes}
For a complex eigenvalue $\lambda = -a + \iu b$ with $a>0$ and $b>0$, the formula \eqref{eq:exact_scalar_phi} is complex 
and cannot be used as a real scalar denominator.
The simplest real choice is the \textit{oscillation-frequency-based} denominator $\frac{2}{b} \tan\bigl(\frac{bh}{2}\bigr)$.
It matches the rotation angle of the Cayley factor exactly only for $a=0$.
Indeed, for $s>0$ the Cayley factor $r(s) = (1+s\lambda/2) / (1-s\lambda/2)$ satisfies, after multiplication with the complex conjugate of its denominator,
\begin{equation}\label{eq:cayley_angle}
    \arg r(s) = \operatorname{atan2}\Bigl(bs,\ 1-\frac{a^2+b^2}{4}\,s^2\Bigr),
\end{equation}
which depends on both $a$ and $b$.
Exact phase matching is nevertheless possible with a real scalar denominator.

\begin{proposition}[Phase-matching denominator]\label{prop:phasematching}
Let $\lambda=-a+\iu b$ with $a\ge0$, $b>0$, and let $\theta=bh\in(0,\pi)$. Then
\begin{equation}\label{eq:phasematching}
    \phi_{\mathrm{ph}}(h) = \frac{2\sin\theta}{b\cos\theta + \sqrt{b^2+a^2\sin^2\theta}}
\end{equation}
is positive and satisfies $\arg r\bigl(\phi_{\mathrm{ph}}(h)\bigr)=\theta=bh$, i.e.\ the Cayley factor rotates by exactly the angle of the exact flow.
Moreover,
\begin{equation}\label{eq:phasematching_expansion}
    \phi_{\mathrm{ph}}(h)=h+\Bigl(\frac{b^2}{12}-\frac{a^2}{4}\Bigr)h^3+\mathcal{O}(h^5),
\end{equation}
and for $a=0$ the formula reduces to $\frac{2}{b}\tan\bigl(\frac{bh}{2}\bigr)$.
\end{proposition}
\begin{proof}
By \eqref{eq:cayley_angle}, $\arg r(s)=\theta$ with $\theta\in(0,\pi)$ holds if and only if $bs\cos\theta=\bigl(1-\frac{a^2+b^2}{4}s^2\bigr)\sin\theta$ with $bs>0$, i.e.\ if and only if $s$ solves
\begin{equation*}
    \frac{a^2+b^2}{4}\,\sin\theta\, s^2+b\cos\theta\, s-\sin\theta=0 .
\end{equation*}
Since $\sin\theta>0$, it has exactly one positive root, and rationalizing it with $b^2\cos^2\theta+(a^2+b^2)\sin^2\theta=b^2+a^2\sin^2\theta$ gives \eqref{eq:phasematching}, whose denominator is positive because $\sqrt{b^2+a^2\sin^2\theta}\ge b>-b\cos\theta$.
A Taylor expansion in the step-size $h$ yields \eqref{eq:phasematching_expansion}.
\end{proof}

The denominator \eqref{eq:phasematching} matches the rotation angle but in general not the damping factor $|\mathrm{e}^{\lambda h}|=\mathrm{e}^{-ah}$.
Simultaneous phase and damping matching requires the matrix construction of Section~\ref{sec:locallyexact}.
Extended to negative $h$, it is odd and of order $q=2$.
It is unbounded as $h\uparrow\pi/b$.

\subsection{Composite and state-dependent denominators}

\begin{table}[t]
\centering
\caption{Scalar denominator functions for \pHNSFD{} schemes.
Here $\alpha>0$ is a decay rate, $\omega>0$ a frequency and $\lambda = -a+\iu b$ ($a\ge0$, $b>0$) a damped oscillatory eigenvalue of the linearization.
``Order'' refers to $q$ with $\phi(h)=h+\mathcal{O}(h^{q+1})$.
By Theorem~\ref{thm:order} the scheme order is at least $\min\{2,q\}$ for a symmetric scheme. ``Range'' is the step-size range on which $\phi$ is positive.}
\label{tab:denominators}
\small
\begin{tabular}{@{}lccll@{}}
\toprule
$\phi(h)$ & order $q$ & range & bounded & exact for \\
\midrule
$h$ & $\infty$ & $h>0$ & no & nothing \\[1mm]
$\dfrac{1-\mathrm{e}^{-\alpha h}}{\alpha}$ & $1$ & $h>0$ & yes ($1/\alpha$)
   & forward Euler on $\dot{\xi}=-\alpha\xi$ \\[2mm]
$\dfrac{2}{\alpha}\tanh\Bigl(\dfrac{\alpha h}{2}\Bigr)$
   & $2$ & $h>0$ & yes ($2/\alpha$) & \pHNSFD{} on $\dot{\xi}=-\alpha\xi$ \\[3mm]
$\dfrac{2}{\omega}\tan\Bigl(\dfrac{\omega h}{2}\Bigr)$
   & $2$ & $0<h<\pi/\omega$ & no & \pHNSFD{} on the modes $\pm\iu\omega$ \\[3mm]
$\dfrac{h}{(1+\tau h^{2m})^{1/(2m)}}$ & $2m$ & $h>0$ & yes ($\tau^{-1/(2m)}$) & --- \\[4mm]
$\dfrac{\sin(\omega h)}{\omega}$ & $2$ & $0<h<\pi/\omega$ & yes ($1/\omega$)
   & --- \\[3mm]
$\phi_{\mathrm{ph}}(h)$ in \eqref{eq:phasematching} & $2$ & $0<h<\pi/b$ & no
   & rotation angle of the mode $-a+\iu b$ \\
\bottomrule
\end{tabular}
\end{table}

Table~\ref{tab:denominators} summarizes the choices.
For a general complex eigenvalue, \eqref{eq:exact_scalar_phi} should be used through a real $2\times2$ block (Section~\ref{ss:modal}) rather than through its real part, for which no positivity guarantee is available.
If several time scales are present, a single scalar $\phi$ can be exact for at most one of them.  
These three strategies can be used to treat this problem: 
\begin{enumerate}[label=(\roman*),leftmargin=2.2em]

  \item \textit{Dominant-mode matching}: choose $\lambda$ as the eigenvalue whose time scale one most wishes to resolve. This is typically the slowest oscillatory mode in a conservative problem, or the slowest decaying mode in a stiff dissipative one; accurate damping of unresolved fast modes is then not guaranteed.
  
  \item \textit{State-dependent denominators} $\phi_n=\phi(h;x_n)$ built from local spectral data of the full Jacobian $Df(x_n,u)$, the pH analogue of Cie\'sli\'nski's locally exact integrators \cite{Cieslinski2013camwa, Cieslinski2013pla}.
    The algebraic identities of Section~\ref{sec:properties} extend with stepwise weights $\phi_n$.
    The LaSalle argument additionally requires $\inf_{x\in K}\phi(h;x)>0$, second order requires $\phi(h;x)=h+\mathcal{O}(h^3)$ uniformly on compact sets, and time symmetry requires a denominator that is symmetric in $(x_n,x_{n+1})$ and odd in the step-size $h$.
  
  \item \textit{Matrix denominators}: replace the scalar denominator $\phi(h)$ by a matrix denominator $\Phi(h)$, so that different modes receive different effective time scales simultaneously. 
    This is the subject of Section~\ref{sec:locallyexact}, and
    it requires care, since a naive matrix denominator destroys the power
    balance, and a structured one need not preserve the equilibria.
\end{enumerate}

\section{Step-Dependent Matrix Extensions and Locally Exact Methods} \label{sec:locallyexact}
A scalar denominator  $\phi(h)$ can be exact for at most one characteristic time scale.
Matrix denominators remove this limitation, but they leave the scalar class of Definition~\ref{def:phnsfd}: the resulting structure operators depend on $h$ and coincide with $J$, $R$, $G$ on the diagonal only asymptotically.
Here we therefore define a separate class of step-dependent methods.

\subsection{Linearization and the step-dependent class}
For a constant input $u^\ast$, we write $f(x,u^\ast) = \bigl( J(x)-R(x) \bigr) \nabla H(x) +G(x)u^\ast$.
The full Jacobian $A(x)=D_xf(x,u^\ast)$ acts on $v\in\R^n$ as
\begin{equation}\label{eq:full_jacobian}
    A(x)v=\bigl(J(x)-R(x)\bigr)\nabla^2H(x)\,v+\bigl(DJ(x)[v]-DR(x)[v]\bigr)\nabla H(x)+DG(x)[v]\,u^\ast .
\end{equation}
At an unforced equilibrium $x^\ast$ with $\nabla H(x^\ast)=0$ the last two terms vanish, and with $\xi=x-x^\ast$, we get
\begin{equation}\label{eq:linearized}
    \dot{\xi} = A_\ast\xi+\mathcal{O}(\|\xi\|^2),
    \qquad A_\ast=\bigl(J(x^\ast)-R(x^\ast)\bigr)\nabla^2H(x^\ast),
\end{equation}
with \emph{exact linearized flow} $\xi(t+h) = \mathrm{e}^{hA_\ast} \xi(t)$.
The remainder in \eqref{eq:linearized} requires a locally Lipschitz derivative of the vector field, e.g.\ $H\in C^3$ and $J,R\in C^2$. 
At a forced equilibrium or a general point, the derivatives of $J$, $R$, $G$ must be retained.
Fitting to the reduced matrix $(J-R)\nabla^2H$ is only a \textit{frozen-structure} fitting.

\begin{definition}[Step-dependent discrete pH methods]\label{def:classH}
For each step size $h$ in a range $\mathcal{H}$, let $J_h,R_h\colon\Omega\times\Omega\to\R^{n\times n}$ and $G_h\colon\Omega\times\Omega\to\R^{n\times m}$ be continuous with
\begin{equation}\label{eq:classH_structure}
    J_h(x,z)^\top=-J_h(x,z),\qquad R_h(x,z)=R_h(x,z)^\top\succeq0 .
\end{equation}
With a discrete gradient $g=\dgrad H(x,z)$, the \textit{step-dependent discrete pH method} reads
\begin{equation}\label{eq:classH}
    z-x = h\bigl[ \bigl(J_h(x,z)-R_h(x,z)\bigr)g + G_h(x,z)\,u_{\nph}\bigr], \qquad
    y_h=G_h(x,z)^\top g ,
\end{equation}
with $x=x_n$, $z=x_{n+1}$.
Its structure operators are \textit{second-order consistent}, if 
\begin{equation}\label{eq:classH_consistency}
    J_h(x,z) = J(m) + \mathcal{O}\bigl(h^2+\|z-x\|^2\bigr),\quad
     \text{with}\; m=\frac{x+z}{2},
\end{equation}
and analogously for $R_h$ and $G_h$, uniformly on the compact set under consideration.
With (DG2) for $g$, smoothness, local solvability and stability, this gives order two.
(DG2) for $g$ cannot be dropped: for $H=\frac{1}{2}x^\top Qx$, $Q=\begin{psmallmatrix}2&1\\1&2\end{psmallmatrix}$, $J=\begin{psmallmatrix}0&1\\-1&0\end{psmallmatrix}$, $R=0$, the Itoh--Abe gradient \eqref{eq:itohabe_quadratic} turns the scheme with $J_h=J$ into $z-x=\frac{h}{1-h/2}JQ\frac{x+z}{2}$, which is of first order.
The method is \textit{equilibrium compatible} at an equilibrium $(x^\ast,u^\ast)$ of \eqref{eq:ph_system} if, for every $h\in\mathcal{H}$,
\begin{equation}\label{eq:classH_equilibrium}
    \bigl[J_h(x^\ast,x^\ast)-R_h(x^\ast,x^\ast)\bigr]\nabla H(x^\ast)+G_h(x^\ast,x^\ast)\,u^\ast=0 .
\end{equation}
\end{definition}

The scalar class is contained in this class, with $J_h=\frac{\phi(h)}{h}\Jb$ etc.
It is equilibrium compatible by exact diagonal consistency.
For general step-dependent methods, second-order consistency does not imply equilibrium compatibility at a fixed $h$.
The proof of Theorem~\ref{thm:powerbalance} gives the following.

\begin{proposition}[Power balance in the step-dependent class]\label{prop:classH_pb}
Every solved step of \eqref{eq:classH} satisfies
\begin{equation}\label{eq:classH_pb}
    H(x_{n+1})-H(x_n)=h\bigl[-g^\top R_h\,g+y_h^\top u_{\nph}\bigr]\le h\,y_h^\top u_{\nph}.
\end{equation}
\end{proposition}

\subsection{Exact matrix denominators}
For a linear system $\dot{\xi}=A\xi$, a matrix-valued denominator $\Phi(h)\in\R^{n\times n}$ replaces $\phi(h)I$.
Since a canonical matrix denominator may be singular, we write the schemes in increment form and avoid division by $\Phi(h)$.
Let
\begin{equation}\label{eq:matrix_functions}
    \varphi_1(\zeta) = \frac{\mathrm{e}^{\zeta}-1}{\zeta},\qquad
    \tanhc(\zeta) := \frac{\tanh\zeta}{\zeta}, \qquad
    \varphi_1(0)=\tanhc(0):=1,
\end{equation}
where $\varphi_1$ is entire, while $\tanhc$ is analytic at $0$ but has poles at $\zeta = \iu\pi(k+\frac{1}{2})$, $k\in\Z$.
Note that $\tanhc$ is defined analogously to the tanc function (tangent cardinal function), cf.\  \cite{Tanhc}.

\begin{proposition}[Exact matrix denominators]\label{prop:exactmatrix}
Let $A\in\R^{n\times n}$, $h>0$ and $E=\mathrm{e}^{hA}$.
Consider the two representations
\begin{equation}\label{eq:matrix_nsfd}
    \text{\emph{(Euler)}}\quad \xi_{n+1}-\xi_n=\Phi_{\mathrm{E}}(h)A\xi_n,
    \qquad
    \text{\emph{(midpoint)}}\quad \xi_{n+1}-\xi_n=\Phi_{\mathrm{M}}(h)A\,\frac{\xi_n+\xi_{n+1}}{2},
\end{equation}
with the canonical choices
\begin{equation}\label{eq:exact_matrix_phi}
    \Phi_{\mathrm{E}}(h)=h\,\varphi_1(hA),
    \qquad
    \Phi_{\mathrm{M}}(h)=h\tanhc\Bigl(\tfrac{h}{2}A\Bigr).
\end{equation}
$\Phi_{\mathrm{E}}$ is defined for every $A$, and $\Phi_{\mathrm{M}}$ is defined if $E+I$ is invertible.
On these domains both representations are exact, $\xi_{n+1}=E\xi_n$.
If $A$ is invertible, the exactness requirement determines $\Phi_{\mathrm{E}}$ and $\Phi_{\mathrm{M}}$ uniquely, namely $\Phi_{\mathrm{E}} = (E-I) A^{-1}$ and $\Phi_{\mathrm{M}} = 2 A^{-1}\tanh\bigl( \frac{h}{2}A \bigr)$.
If $A$ is singular, \eqref{eq:exact_matrix_phi} are canonical choices, but not the only ones.
Moreover, $\Phi_{\mathrm{M}}$ is odd in $h$, and $\Phi_{\mathrm{M}}(h) = hI - \frac{h^3}{12}A^2 +\mathcal{O}(h^5)$.
\end{proposition}
\begin{proof}
The Euler case follows from $h\varphi_1(hA)A=E-I$.
For the midpoint case, $Z=\frac{1}{2}\Phi_{\mathrm{M}}A = \tanh\bigl(\frac{h}{2}A \bigr) = (E-I)(E+I)^{-1}$, so $I-Z = 2(E+I)^{-1}$ and $(I-Z)^{-1}(I+Z) = E$.
Conversely, exactness requires $\Phi_{\mathrm{M}}A=2(E-I)(E+I)^{-1}$, which determines $\Phi_{\mathrm{M}}$ if $A$ is invertible.
For singular $A$, every exact denominator has the form $\Phi_{\mathrm{M}}+N$ (resp.\ $\Phi_{\mathrm{E}}+N$) with $NA=0$, i.e.\ $N$ vanishes on $\ran A$ and is arbitrary on a complement of $\ran A$.
This complement can be taken as $\ke A$ only if the eigenvalue $0$ is semisimple (for $A=\begin{psmallmatrix}0&1\\0&0\end{psmallmatrix}$, $\ke A=\ran A$). 
\end{proof}

The midpoint formula is the matrix form of Cie\'sli\'nski's locally exact construction \cite{Cieslinski2013camwa, Cieslinski2013pla}, and, as in Proposition~\ref{prop:tanhvsexp}, $\Phi_{\mathrm{M}}$ rather than $\Phi_{\mathrm{E}}$ is relevant for the midpoint-type scheme.
The quotient form $\Phi(h)^{-1}(\xi_{n+1}-\xi_n)$ additionally requires $\Phi(h)$ to be invertible, which excludes the resonant steps $h\lambda\in2\pi\iu(\Z\setminus\{0\})$.

\subsection{Structured matrix denominators}\label{Sec8.3}
A matrix denominator generally destroys the power balance: 
repeating the proof of Theorem~\ref{thm:powerbalance} with $x_{n+1}-x_n=\Phi\bigl[(\Jb_n-\Rb_n)\dgrad H_n+\Gb_nu_{\nph}\bigr]$ gives the \emph{discrete power balance}
\begin{equation}\label{eq:matrixpb}
    H(x_{n+1})-H(x_n) = \dgrad H_n^\top\Phi\Jb_n\dgrad H_n
    -\dgrad H_n^\top \Phi\Rb_n\dgrad H_n + \dgrad H_n^\top \Phi\Gb_nu_{\nph},
\end{equation}
and the first term no longer vanishes because $\Phi\Jb_n$ need not be
skew-symmetric.  Two remedies are available.

\begin{theorem}[Congruence form]\label{thm:matrixdenom}
Let $P=P(h)$ be symmetric positive definite for $h\in\mathcal{H}$, let $\Jb$, $\Rb$, $\Gb$ be admissible, and set
\begin{equation}\label{eq:congruence_scheme}
    J_h = P^{1/2}\Jb P^{1/2},\qquad
    R_h=P^{1/2}\Rb P^{1/2},\qquad
    G_h=P^{1/2}\Gb .
\end{equation}
\begin{enumerate}[label=(\alph*),leftmargin=2.2em]
  
  \item $J_h$ is skew-symmetric and $R_h$ is symmetric positive semidefinite, so \eqref{eq:classH} with the opera\-tors \eqref{eq:congruence_scheme} is a step-dependent discrete pH method and satisfies the power balance \eqref{eq:classH_pb} with $y_h=G_h^\top\dgrad H$.
  
  \item If $\Jb$, $\Rb$, $\Gb$ satisfy (DG2) and $P(h)=I+\mathcal{O}(h^2)$, the operators \eqref{eq:congruence_scheme} are second-order consistent.
  With (DG2) for $\dgrad H$, the \emph{method has consistency order two}.
  
  \item If $\dgrad H$, $\Jb$, $\Rb$, $\Gb$ are symmetric in their two arguments and $P(-h)=P(h)$, the method is \emph{symmetric under step reversal}.
  
  \item The method is \emph{equilibrium compatible} at $(x^\ast,u^\ast)$ if and only if
  \begin{equation}\label{eq:congruence_equilibrium}
      P^{1/2}\bigl( J(x^\ast) - R(x^\ast) \bigr) P^{1/2} \nabla H(x^\ast) +P^{1/2} G(x^\ast)u^\ast=0
      \qquad\text{for all } h\in\mathcal{H}.
  \end{equation}
  Condition \eqref{eq:congruence_equilibrium} holds whenever $u^*=0$ and $\nabla H(x^*) = 0$.
For general $P$, this condition does not necessarily hold at forced equilibria or at unforced equilibria with $\nabla H(x^*)\ne0$, including steady rotations of the rigid body in Example~\ref{ex:rigidbody}.
\end{enumerate}
If the denominator is retained in the form $(x_{n+1}-x_n)/\phi(h)$, all three operators in \eqref{eq:congruence_scheme} are multiplied by $\phi(h)/h$.
Second order then additionally requires $\phi(h) = h+\mathcal{O}(h^3)$.
\end{theorem}
\begin{proof}
(a) $\bigl(P^{1/2}\Jb P^{1/2}\bigr)^\top=P^{1/2}\Jb^\top P^{1/2}=-P^{1/2}\Jb P^{1/2}$, and $P^{1/2}\Rb P^{1/2}\succeq0$ is a congruence of a positive semidefinite matrix, Proposition~\ref{prop:classH_pb} applies.
(b) With $P^{1/2}=I+\mathcal{O}(h^2)$, \eqref{eq:classH_consistency} follows from (DG2).
(c) Interchanging $x$ and $z$ and replacing $h$ by $-h$ leaves \eqref{eq:classH} invariant.
(d) This is \eqref{eq:classH_equilibrium} with \eqref{eq:congruence_scheme}.
If $u^\ast=0$ and $\nabla H(x^\ast)=0$, the left-hand side vanishes.
\end{proof}

Thus congruence preserves the algebraic pH constraints, whereas consistency and equilibrium compatibility are separate requirements.
$P(h)\to I$ only ensures that a shifted discrete equilibrium converges to the correct one as $h\to0$.

\begin{example}[Shift of a forced equilibrium]\label{ex:congruence_shift}
Consider $\dot{x} = -x+u$, i.e.\ $H(x)=\frac{1}{2}x^2$, $J=0$, $R=1$, $G=1$, with the equilibrium $x^\ast=u$.
The congruence form with scalar $P>0$ gives the discrete equilibrium equation $-Px+\sqrt{P}\,u = 0$, hence $x_h^\ast = u/\sqrt{P}$.
For $P=(1+h^2)^2$ one obtains $x_h^\ast = u/(1+h^2)\neq u$.
At $h=1$, $u=x_n=1$, the step equation $x_{n+1}-1 = -4\cdot\frac{1+x_{n+1}}{2}+2$ gives $x_{n+1} = \frac{1}{3}$, so the original equilibrium is not preserved.
\end{example}

The second remedy preserves the full equilibrium set.

\begin{theorem}[Left multiplication]\label{thm:leftmultiplier}
Let $L_h$ be symmetric positive definite and commute with $\Jb(x,z)$ and $\Rb(x,z)$ for all pairs considered (for instance, $J$ and $R$ constant and $L_h$ a function of $J$ and $R$ that commutes with both).
Set $J_h = L_h\Jb$, $R_h = L_h\Rb$ and $G_h = L_h\Gb$.
Then $J_h$ is skew-symmetric, $R_h$ is symmetric positive semidefinite, the power balance \eqref{eq:classH_pb} holds with $y_h=G_h^\top\dgrad H$, and on the diagonal the modified vector field equals $L_hf(x,u)$.

Consequently, with constant input $u^\ast$, the diagonal solutions of the step equation are exactly the equilibria of \eqref{eq:ph_system}, and so are the fixed points of a step map whose branch retains these diagonal solutions. 
If $L_h=I+\mathcal{O}(h^2)$ and $\Jb$, $\Rb$, $\Gb$ satisfy (DG2), the operators are second-order consistent.
\end{theorem}
\begin{proof}
$(L_h\Jb)^\top=\Jb^\top L_h=-\Jb L_h=-L_h\Jb$.
Since $L_h^{1/2}$ is a polynomial in $L_h$, it commutes with $\Rb$, so $L_h\Rb=L_h^{1/2}\Rb L_h^{1/2}$ is symmetric positive semidefinite.
The power balance follows from Proposition~\ref{prop:classH_pb}.
On the diagonal, we have $(J_h-R_h) \nabla H + G_hu^\ast = L_h f(x,u^\ast)$, and $L_h$ is invertible.
\end{proof}

\subsection{The locally exact scheme for linear systems}
For linear pH systems the exact matrix denominator admits a step-dependent pH representation without any congruence.
\begin{theorem}[Locally exact scheme for linear pH systems]\label{thm:locallyexactisph}
Let $A=(J-R)Q$ with $J=-J^\top$, $R=R^\top\succeq0$, $Q\succ0$, let $h>0$ be such that $E+I$ is invertible, $E=\mathrm{e}^{hA}$, and set
\begin{equation*}
    M_h:=\Phi_{\mathrm{M}}(h)A=2\tanh\Bigl(\tfrac{h}{2}A\Bigr),\quad
    K_h:=\frac{1}{h}M_hQ^{-1},\quad
    J_h:=\skewpart(K_h), \quad
    R_h:=-\sym(K_h).
\end{equation*}
Then
\begin{enumerate}[label=(\roman*),leftmargin=2.2em]
  \item $J_h^\top=-J_h$ and $R_h\succeq0$;
  \item the scheme
  \begin{equation}\label{eq:locallyexact_scheme}
      x_{n+1}-x_n = h \bigl( J_h-R_h \bigr)Q\,\frac{x_n+x_{n+1}}{2}
  \end{equation}
    is a step-dependent discrete pH method in the sense of Definition~\ref{def:classH}, with constant operators;
  \item it is \textit{exact}: $x_{n+1}=Ex_n$;
  \item $J_h=J+\mathcal{O}(h^2)$ and $R_h=R+\mathcal{O}(h^2)$;
  \item it satisfies the \emph{exact power balance} $H(x_{n+1})-H(x_n) = -h\,v^\top QR_hQ\,v$, $v=\frac{x_n+x_{n+1}}{2}$, and, by exactness, $H(x_{n+1})-H(x_n) = -\int_{t_n}^{t_{n+1}}\mathcal{D}(x(t))\,\dd t$.
\end{enumerate}
\end{theorem}
\begin{proof}
Since $J_h-R_h=\skewpart(K_h)+\sym(K_h)=K_h$, the right-hand side of \eqref{eq:locallyexact_scheme} equals $M_h\frac{x_n+x_{n+1}}{2}$, which is the midpoint representation of Proposition~\ref{prop:exactmatrix}.
This gives (ii) with $Q\frac{x_n+x_{n+1}}{2} = \dgrad H(x_n,x_{n+1})$ (Lemma~\ref{lem:quadratic}) and (iii).
Skew-symmetry of $J_h$ is trivial.
For $R_h\succeq0$, note that $H$ is quadratic, so
\begin{equation*}
    H(x_{n+1})-H(x_n) = \Bigl(Q\frac{x_n+x_{n+1}}{2}\Bigr)^\top(x_{n+1}-x_n)
    =v^\top\sym(QM_h)\,v,
    \qquad v=\frac{1}{2}\bigl(I+E\bigr)x_n .
\end{equation*}
By (iii) and Lemma~\ref{lem:contract}, $H(x_{n+1})\le H(x_n)$ for every $x_n$.
Since $I+E$ is invertible, $v$ ranges over all of $\R^n$, whence $\sym(QM_h)\preceq0$.
Finally,
\begin{equation*}
    R_h=-\frac{1}{h}\sym\bigl(M_hQ^{-1}\bigr)=-\frac{1}{h}\,Q^{-1}\sym(QM_h)\,Q^{-1}\succeq0 .
\end{equation*}
For (iv), the expansion $M_h=hA-\frac{h^3}{12}A^3+\mathcal{O}(h^5)$ gives $K_h=AQ^{-1}+\mathcal{O}(h^2)=J-R+\mathcal{O}(h^2)$, and the skew-symmetric and symmetric parts converge accordingly.
Statement (v) is Proposition~\ref{prop:classH_pb} together with (iii) and the continuous power balance.
\end{proof}

The theorem is a statement in the step-dependent class: the operators $J_h$, $R_h$ converge to $J$, $R$ with an $\mathcal{O}(h^2)$ error, but they do not satisfy the exact diagonal consistency \eqref{eq:discrete_consistency} for a fixed $h$.
For $u=0$ every equilibrium $x^\ast\in\ke A$ is preserved, since $M_hx^\ast=0$.
Whether the scheme has \textit{additional fixed points} is a separate question.

\begin{proposition}[Nonresonance]\label{prop:nonresonance}
Let $E=\mathrm{e}^{hA}$ and assume
\begin{equation}\label{eq:nonresonance}
    \spec(A) \cap \Bigl\{\frac{2\pi\iu k}{h} \colon k\in\Z\setminus\{0\}\Bigr\}=\emptyset\,.
\end{equation}
Then $\ke(E-I)=\ke A$, so the fixed-point set of the exact time-$h$ map, and hence of \eqref{eq:locallyexact_scheme}, is $\ke A$.
For the affine map $x\mapsto Ex+h\varphi_1(hA)Bu$ with constant $u$, the fixed points are exactly the solutions of $Ax+Bu=0$.
The condition $\spec(A)\cap\{\iu k\pi/h\colon k\in\Z\setminus\{0\}\}=\emptyset$ implies both \eqref{eq:nonresonance} and the invertibility of $E+I$.
If $A$ is Hurwitz, both conditions hold for every $h>0$.
For a conservative system with frequencies $\omega_j$, the restriction $0<h<\pi/\max_j\omega_j$ is sufficient.
\end{proposition}
\begin{proof}
Write $E-I=hA\,\varphi_1(hA)$.
The zeros of $\varphi_1$ are $2\pi\iu k$, $k\in\Z\setminus\{0\}$, so \eqref{eq:nonresonance} makes $\varphi_1(hA)$ invertible.
Since $A$ and $\varphi_1(hA)$ commute, we have $\ke(E-I)=\ke A$.
For the affine map, a fixed point satisfies $h\varphi_1(hA)(Ax+Bu)=0$, i.e.\ $Ax+Bu=0$.
The matrix $E+I$ is singular if and only if $h\lambda\in\iu\pi(2\Z+1)$ for some $\lambda\in\spec(A)$, and extra fixed points require $h\lambda\in2\pi\iu(\Z\setminus\{0\})$.
Both are excluded by the stated condition, which allows zero eigenvalues.
If $A$ is Hurwitz, $h\lambda$ is never purely imaginary.
\end{proof}

The condition cannot be dropped: for the harmonic oscillator with $h=2\pi/\omega$ one has $\mathrm{e}^{hA}=I$, so every point is a fixed point of the sampled exact map, although $\mathrm{e}^{hA}+I=2I$ is invertible.

For linear systems with an input that is constant on each step, the analogous construction $x_{n+1}-x_n=\Phi_{\mathrm{M}}(h)\bigl(A\frac{x_n+x_{n+1}}{2}+Bu\bigr)$ is exact, $x_{n+1}=Ex_n+h\varphi_1(hA)Bu$, admits the same operators $J_h$, $R_h$ together with the modified port matrix $B_h=\Phi_{\mathrm{M}}(h)B/h$, and its fixed points are the solutions of $Ax+Bu=0$ if $\Phi_{\mathrm{M}}(h)$ is invertible.

\begin{remark}[Local accuracy near an equilibrium]\label{rem:localaccuracy}
If a method preserves an equilibrium $x^\ast$ and has the exact linearization there, the \emph{one-step error} $\mathcal{R}_h = \Psi_h-\mathcal{E}_h$ satisfies $\mathcal{R}_h(x^\ast) = 0$ and $D\mathcal{R}_h(x^\ast) = 0$, hence $\mathcal{R}_h(x^\ast+\xi) = \mathcal{O}(\|\xi\|^2)$ for fixed $h$, provided $D\mathcal{R}_h$ is locally Lipschitz (e.g.\ for a $C^2$ vector field and a $C^2$ step map). 
Local exactness thus removes the error component that is linear in the perturbation.
An estimate $\mathcal{O}(h^{p+1}\|\xi\|^2)$ requires the additional bound $\|D^2\mathcal{R}_h\|\le Ch^{p+1}$ near $x^\ast$.
This is consistent with \cite{Cieslinski2013camwa} and with the pendulum results in Section~\ref{ssec:exp4}.
\end{remark}

Evaluating $\Phi_{\mathrm{M}}(h)$ requires one matrix function evaluation per step size, e.g.\ by Krylov methods for large $n$ \cite{HochbruckOstermann2010}, compared with one scalar function evaluation for the scalar class.

\subsection{Blockwise and modal constructions}\label{ss:modal}
An exact blockwise representation is obtained from the spectral structure of $A_\ast$:
if $A_\ast = T\,\blkdiag(A_1,\dots,A_k)\,T^{-1}$, 
one sets $\Phi = T\,\blkdiag(\Phi_1,\dots,\Phi_k)\,T^{-1}$ with $\Phi_j = h\tanhc\bigl( \frac{h}{2} A_j \bigr)$.
For dissipative real modes $\lambda_j = -\alpha_j$, $\alpha_j>0$, and for oscillatory eigenvalues $\lambda_j = \pm \iu\omega_j$, the blocks are the scalar denominators
\begin{equation}\label{eq:blockdissipative}
    \phi_j(h) = \frac{2}{\alpha_j}\tanh\Bigl(\frac{\alpha_jh}{2}\Bigr)
    \qquad\text{and}\qquad
    \phi_j(h) = \frac{2}{\omega_j} \tan\Bigl( \frac{\omega_j h}{2} \Bigr),
\end{equation}
respectively, the latter for $0<h<\pi/\omega_j$. 
For a real $2\times2$ block $A_j = -a_jI_2 + b_jK_2$, where $K_2 = \bigl(\begin{smallmatrix}0&1\\-1&0\end{smallmatrix}\bigr)$,
the matrix function is represented as 
\begin{equation*}
  \Phi_j = \Real\bigl(\varphi_{-a_j + \iu b_j}(h) \bigr)I_2
  + \Imag\bigl( \varphi_{-a_j + \iu b_j}(h) \bigr)K_2,
\end{equation*}
with $\varphi_\lambda$ from the scalar denominator \eqref{eq:exact_scalar_phi}.

This representation assigns different effective time scales
to different modes. 
For a linear system satisfying the hypotheses of Theorem~\ref{thm:locallyexactisph}, the resulting denominator has a step-dependent pH representation that depends on the step, as proven in that theorem. 
However, for a frozen linearization of a nonlinear system, structural preservation requires a separate verification.
Theorems~\ref{thm:matrixdenom} and~\ref{thm:leftmultiplier} provide structured constructions under their respective hypotheses.

For conservative linear pH systems, the construction can be made fully explicit. 
It is simultaneously exactly energy-preserving and exactly frequency-matching for \textit{all} modes.
The following result provides an explicit modal realization of
the conservative case $R=0$ of Theorem~\ref{thm:locallyexactisph}, under the additional invertibility assumption on $J$.

\begin{theorem}[Modal locally exact conservative scheme]\label{thm:williamson}
Let $\dot{x} = JQx$ with $J=-J^\top$ invertible and $Q\succ0$, and let $\omega_1,\dots,\omega_{n/2}>0$ be the symplectic eigenvalues of $Q$ relative to $J$.  
By Williamson's theorem \cite{Williamson1936} 
there is an invertible $S\in\R^{n\times n}$ with
\begin{equation}
    S^{-1} JS^{-\top} = J_0 := \blkdiag\Bigl(\begin{psmallmatrix}0&1\\-1&0\end{psmallmatrix},
    \dots\Bigr), \qquad
    S^\top QS = \blkdiag\bigl( \omega_1 I_2,\dots,\omega_{n/2} I_2 \bigr) =: D .
\end{equation}
Let $0<h<\pi/\max_j\omega_j$ and define, in the coordinates $z=S^{-1}x$, the scheme
\begin{equation}\label{eq:modalscheme}
    z_{n+1}-z_n = J_0\,\Phi_D\,D\,\frac{z_n+z_{n+1}}{2},\qquad
    \Phi_D := \blkdiag\bigl( \phi_1(h)I_2,\dots,\phi_{n/2}(h)I_2 \bigr),
\end{equation}
with $\phi_j(h) = \frac{2}{\omega_j}\tan\bigl( \frac{\omega_jh}{2} \bigr)$.  Then
\begin{enumerate}[label=(\roman*),leftmargin=2.2em]
  \item $J_0\Phi_D$ is skew-symmetric, so \eqref{eq:modalscheme} is a step-dependent discrete pH method with $J_h = J_0\Phi_D/h = J_0 +\mathcal{O}(h^2)$ and $R_h=0$, and it conserves $\frac{1}{2}z^\top Dz = \frac{1}{2}x^\top Qx$ exactly;
  \item \eqref{eq:modalscheme} is \textit{exact}: $z_{n+1}=\mathrm{e}^{hJ_0D}z_n$, hence  $x_{n+1}=\mathrm{e}^{hJQ}x_n$;
  \item the scheme reproduces every characteristic frequency $\omega_j$ without phase error, and its only fixed point is the origin.
\end{enumerate}
\end{theorem}
\begin{proof}
With $x=Sz$, $\dot{z}=J_0Dz$ and $H=\frac{1}{2}z^\top Dz$.
$\Phi_D$ commutes with $J_0$, so $J_0\Phi_D$ is skew-symmetric, which gives (i) by Proposition~\ref{prop:classH_pb}.
Each $2\times2$ block is a harmonic oscillator, discretized with its exact denominator (Proposition~\ref{prop:exactscalar}), which gives (ii) and (iii), the latter together with Proposition~\ref{prop:nonresonance}.
\end{proof}

Theorem~\ref{thm:williamson} shows, in the linear conservative case, that energy preservation and spectral exactness are simultaneously achievable.
The extension to $R\neq0$ is Theorem~\ref{thm:locallyexactisph}.
Locally exact, energy-preserving schemes are already available for nonlinear canonical Hamiltonian systems \cite{Cieslinski2013pla}, and extending the present combination of pH structure, nonnegative discrete dissipation, equilibrium compatibility and local exactness to nonlinear systems with state-dependent structure and port operators remains a subject for further investigation. 

\subsection{A dissipative rational approximation}
When $S$ or $\tanh\bigl(\frac{h}{2}A\bigr)$ is unavailable, one may use rational approximations.
The $[1/2]$ Pad\'e approximant $\tanh(s)\approx s/(1+s^2/3)$ in $s$ gives
\begin{equation}\label{eq:rational_phi}
    \Phi(h)\approx h\Bigl(I+\frac{h^2}{12}A^2\Bigr)^{-1},
\end{equation}
which differs from $\Phi_{\mathrm{M}}(h)$ by $\mathcal{O}(h^5)$ and requires only one linear solve with a matrix of the same sparsity pattern as $A^2$.
For this approximation the sign conditions of Theorem~\ref{thm:locallyexactisph} can be proved.

\begin{proposition}[Dissipativity of the rational approximation]\label{prop:rational}
Let $A=(J-R)Q$ as in Theorem~\ref{thm:locallyexactisph}, $B:=A^\top Q +QA = -2QRQ\preceq0$, $C:=I+\frac{h^2}{12}A^2$, and assume that $C$ is invertible. Let $M:=hAC^{-1}$.
Then we have
\begin{equation}\label{eq:rational_identity}
    C^\top\bigl(QM+M^\top Q\bigr)C=h\Bigl(B+\frac{h^2}{12}A^\top BA\Bigr)\preceq0 ,
\end{equation}
so that $R_h:=-\frac{1}{h}\sym(MQ^{-1})=-\frac{1}{h}Q^{-1}\sym(QM)Q^{-1}\succeq0$.
The resulting scheme $x_{n+1}-x_n=M\frac{x_n+x_{n+1}}{2}$ has the step
\begin{equation}\label{eq:gauss_step}
    x_{n+1}=\Bigl(I-\frac{h}{2}A+\frac{h^2}{12}A^2\Bigr)^{-1}\Bigl(I+\frac{h}{2}A+\frac{h^2}{12}A^2\Bigr)x_n ,
\end{equation}
which is the linear step of the two-stage Gauss method, i.e.\ the $[2/2]$ Pad\'e approximation of the exponential, of order four on linear problems.
\end{proposition}

\begin{proof}
Since $MC=hA$, the left-hand side equals $h(C^\top QA+A^\top QC)=h\bigl(B+\frac{h^2}{12}((A^\top)^2QA+A^\top QA^2)\bigr)=h\bigl(B+\frac{h^2}{12}A^\top BA\bigr)$.
Since $B\preceq0$, the right-hand side of \eqref{eq:rational_identity} is negative semidefinite, and congruence with $C^{-1}$ gives $QM+M^\top Q\preceq0$.
Finally, since all matrices involved are rational functions of $A$ and commute, $(I-\frac{1}{2}M)^{-1}(I+\frac{1}{2}M)=(C-\frac{h}{2}A)^{-1}(C+\frac{h}{2}A)$, which is \eqref{eq:gauss_step}.
\end{proof}

The scheme \eqref{eq:gauss_step} of Proposition~\ref{prop:rational} is therefore not a new linear integrator but a pH representation of the classical \emph{fourth-order two-stage Gauss method} \cite{HairerLubichWanner2006}, whose dissipation sign the proposition explains.
$C$ is invertible for every step-size $h>0$ if $A$ is Hurwitz, whereas for a mode $\pm\iu\omega$ it is singular at $h\omega=\sqrt{12}$.

\section{Model Problems}\label{sec:models}

\subsection{Damped Hamiltonian oscillator}\label{ssec:model_damped}

For Example~\ref{ex:damped}, $H$ is a diagonal quadratic form, so the AVF, Gonzalez and Itoh--Abe gradients coincide with the midpoint gradient (Lemma~\ref{lem:quadratic}), and the \pHNSFD{} scheme with constant $J$, $R$, $G$ reads
\begin{equation}\label{eq:pq_scheme}
    \frac{q_{n+1}-q_n}{\phi(h)} = \frac{p_{n+1}+p_n}{2},
    \qquad
    \frac{p_{n+1}-p_n}{\phi(h)}
    = -\omega^2 \,\frac{q_{n+1}+q_n}{2} - \gamma \,\frac{p_{n+1}+p_n}{2}
       + u_{\nph}.
\end{equation}
Its amplification matrix is $\mathcal{C}_{\phi(h)}(A)$ from \eqref{eq:cayley} with $A=\begin{psmallmatrix}0&1\\-\omega^2&-\gamma\end{psmallmatrix}$, and for $u\equiv0$ the energy satisfies $H_{n+1}-H_n = -\phi(h)\,\gamma \bigl(\frac{p_{n+1}+p_n}{2}\bigr)^2\le 0$ for every $h>0$.

\paragraph{Choice of $\phi$}
The eigenvalues are $\lambda_{\pm}=-\frac{\gamma}{2}\pm\iu\omega_d$ with
$\omega_d=\sqrt{\omega^2-\gamma^2/4}$ in the underdamped case.  
Natural options are $\phi(h)=h$, $\frac{2}{\omega_d}\tan\bigl(\frac{\omega_dh}{2}\bigr)$ (phase matched only for $\gamma=0$), the phase-matching denominator \eqref{eq:phasematching} with $a=\frac{\gamma}{2}$, $b=\omega_d$, and the exact matrix denominator $\Phi_{\mathrm{M}}(h)=2A^{-1}\tanh\bigl(\frac{h}{2}A\bigr)$ (Theorem~\ref{thm:locallyexactisph}, Section~\ref{ssec:exp3}).

\subsection{Nonlinear pendulum}\label{ssec:model_pendulum}
For Example~\ref{ex:pendulum} the AVF discrete gradient \eqref{eq:avf} can be evaluated in closed form,
\begin{equation}\label{eq:pendulum_dgrad}
    \dgrad H(x_n,x_{n+1})
    =\begin{pmatrix}
       \frac{\cos q_n - \cos q_{n+1}}{q_{n+1}-q_n} \\[1mm]
       \frac{p_n+p_{n+1}}{2}
     \end{pmatrix},
\end{equation}
with the removable singularity at $q_{n+1}=q_n$ replaced by $\sin q_n$.
To avoid cancellation, it is evaluated as $\sin\bigl(\frac{q_n+q_{n+1}}{2}\bigr) \operatorname{sinc} \bigl(\frac{q_{n+1}-q_n}{2}\bigr)$, 
$\operatorname{sinc}(a) = \sin(a)/a\approx1-\frac{a^2}{6}+\frac{a^4}{120}$ for small $|a|$. 
For Table~\ref{tab:exp4}, however, the implementation used the above quotient directly, replacing it by $\sin q_n$ whenever $|q_{n+1}-q_n|<10^{-12}$.
The factor $\operatorname{sinc}$ distinguishes it from the implicit midpoint rule, which does not conserve $H$ exactly.
%
The \emph{\pHNSFD{} scheme} reads
\begin{equation}\label{eq:pendulum_scheme}
    \frac{q_{n+1}-q_n}{\phi(h)}=\frac{p_n+p_{n+1}}{2},
    \qquad
    \frac{p_{n+1}-p_n}{\phi(h)}=-\frac{\cos q_n-\cos q_{n+1}}{q_{n+1}-q_n} .
\end{equation}
It conserves $H$ exactly for every solved step and every $\phi>0$, so the numerical trajectory stays on the exact level set of $H$.
Energy preservation alone, however, does not guarantee that a librating solution remains in its potential well: 
since the potential is periodic, the same energy level is attained in every well $\bigl((2k-1)\pi,(2k+1)\pi\bigr)$ when $q$ is tracked on $\R$, and the implicit step may select a solution in a different well.
For instance, $q_n=2\pi n$, $p_n=1$ satisfies \eqref{eq:pendulum_scheme} with $\phi=2\pi$ and preserves $H_n=\frac{1}{2}<2$, while the unwrapped angle grows without bound.
The following result gives a sufficient condition that excludes this behavior.

\begin{theorem}[Unique AVF step and preservation of the libration well]\label{thm:pendulum}
Let $s=\phi(h)$ with $0<s<2$. Then:
\begin{enumerate}[label=(\alph*),leftmargin=2.2em]
  \item for every $(q_n,p_n) \in\R^2$, the step \eqref{eq:pendulum_scheme} has a unique solution $(q_{n+1},p_{n+1}) \in\R^2$, which depends continuously on $(q_n,p_n,s)$;
  \item if $H(q_n,p_n)<2$ and $q_n \in \bigl( (2k-1)\pi, (2k+1)\pi \bigr)$ for some $k\in\Z$, then $q_{n+1}$ lies in the same interval.
\end{enumerate}
Consequently, a librating initial condition remains in its well for all $n$.
\end{theorem}
\begin{proof}
Let $d=q_{n+1}-q_n$. Eliminating $p_{n+1}=2d/s-p_n$ from \eqref{eq:pendulum_scheme}, with the AVF form $\frac{\cos q_n-\cos q_{n+1}}{q_{n+1}-q_n}=\int_0^1\sin(q_n+\theta d)\,\dd\theta$, gives the scalar equation
\begin{equation}\label{eq:pendulum_scalar}
    F(d):=d-sp_n+\frac{s^2}{2}\int_0^1\sin(q_n+\theta d)\,\dd\theta=0 .
\end{equation}
Its derivative satisfies
\begin{equation*}
    F'(d)=1+\frac{s^2}{2}\int_0^1\theta\cos(q_n+\theta d)\,\dd\theta\ \ge\ 1-\frac{s^2}{4}>0 .
\end{equation*}
Since the integral term in \eqref{eq:pendulum_scalar} is bounded, $F(d)\to\pm\infty$ as $d\to\pm\infty$, so $F$ has exactly one root, and the implicit function theorem gives continuous dependence on $(q_n,p_n,s)$, this proves (a).
For (b), we consider the solution branch $\sigma\mapsto d(\sigma)$ for $\sigma\in(0,s]$.
It depends continuously on $\sigma$, and $d(\sigma)\to0$ as $\sigma\to0$, because $(1-\sigma^2/4)\,|d(\sigma)| \le |F(0)| \le \sigma|p_n| + \sigma^2/2$ by the mean value theorem.
Along the branch, the energy is preserved: $H\bigl( q_n+d(\sigma),p_{n+1}(\sigma) \bigr) = H(q_n,p_n)<2$.
If $q_n+d(s)$ were outside the well, by continuity some $\sigma$ would give $q_n+d(\sigma) = (2k\pm1)\pi$, where $1-\cos q=2$ and hence $H\ge2$, a contradiction.
Repea\-ting the argument at each step proves the final statement.
\end{proof}

For the locally exact choice $\phi(h) = 2\tan(h/2)$, the sufficient (not necessarily sharp) condition $0<\phi(h)<2$ is equivalent to $0<h<\pi/2$.
A bounded denominator such as \eqref{eq:phi_family_bounded} with $\tau^{-1/(2m)}<2$ gives the result for every $h>0$, at the price of the exact frequency fitting at the origin.
The monotone scalar equation \eqref{eq:pendulum_scalar} can be solved by a bracketing method.

\subsection{Forced port-Hamiltonian system and the power-balance residual}
For \eqref{eq:ph_nsfd} with $u\neq0$ we monitor the residual
\begin{equation}\label{eq:residual}
    \mathcal{R}_{H,n}
    :=    \frac{H_{n+1}-H_n}{\phi(h)} +\mathcal{D}_n
    -    y_{\nph}^\top u_{\nph},
\end{equation}
which by Theorem~\ref{thm:powerbalance} vanishes for the exact solution of the step equation and equals $\dgrad H_n^\top r/\phi(h)$ for a solver residual $r$.
It is therefore a useful code-verification diagnostic.
If the discrete gradient is only approximated by $\tilde{g}$ (e.g.\ by inexact quadrature), with $\mathcal{D}_n$, $y_{\nph}$ computed from $\tilde{g}$, the residual equals $(\tilde{g}^\top r+\varepsilon_{\mathrm{DG}})/\phi(h)$ with the chain-rule defect $\varepsilon_{\mathrm{DG}}=H_{n+1}-H_n-\tilde{g}^\top(x_{n+1}-x_n)$, which is not reduced by the solver tolerance. 

\section{Numerical Experiments}\label{sec:numerics}
All experiments were carried out in double precision using NumPy and SciPy. 
We solved the nonlinear pendulum steps using \texttt{scipy.optimize.fsolve}, which employs a MINPACK hybrid method with a finite-difference approximation of the Jacobian when the Jacobian is not supplied. 
We set \texttt{xtol} to $10^{-13}$ and initialized the solver with the explicit Euler predictor $x_n+\varphi(h)J\nabla H(x_n)$.
The unforced linear experiments were computed using direct linear solves or scalar amplification factors. 
For the forced linear experiment in Section~\ref{ssec:exp5}, we also used \texttt{fsolve} with \texttt{xtol} set to $10^{-13}$, initialized at $x_n$.
For the pendulum and forced linear experiments, each returned state was independently checked against the residual criterion $\|F_n(x_{n+1})\|_2<10^{-13}$, where $F_n=0$ denotes the corresponding step equation. 
All returned states satisfied this criterion, including those for which \texttt{fsolve} reported a non-success termination status.
The reference solutions for the homogeneous linear problems are $x(t) = \mathrm{e}^{tA} x_0$. 
The exact pendulum period is $T(q_0)=4K(\sin^2(q_0/2))$, where $K$ denotes the complete elliptic integral of the first kind in the parameter convention.
Unless otherwise stated, errors are relative final-time errors $E(h)=\|x_N-x(T)\|_2/\|x(T)\|_2$, observed orders are $\log_2\bigl(E(h)/E(h/2)\bigr)$, and the relative energy drift is $\max_{0\le n\le N}|H_n-H_0|/|H_0|$. 
Quantities at the round-off level may vary with the platform and numerical library versions.

\subsection{Conservative oscillator: energy versus phase} \label{ssec:exp1}
We integrate Example~\ref{ex:damped} with $\gamma=0$, $\omega=2$, $x_0=(1,0)^\top$, comparing the implicit midpoint / discrete-gradient scheme ($\phi=h$) with the \pHNSFD{} scheme using the denominator function $\phi(h) = \frac{2}{\omega}\tan(\frac{\omega h}{2})$.

\begin{table}[h]
\centering
\caption{Conservative oscillator, $\omega=2$, $T=20$. 
Both methods conserve the energy to machine precision.
Only the \pHNSFD{} scheme is free of phase error. 
``ord'' is the observed order of the midpoint error.}
\label{tab:exp1}
\footnotesize
\setlength{\tabcolsep}{7pt}
\begin{tabular}{@{}lccccc@{}}
\toprule
& \multicolumn{3}{c}{error at $T$} & \multicolumn{2}{c}{$\max_n|H_n-H_0|/H_0$}\\
\cmidrule(lr){2-4}\cmidrule(l){5-6}
$h$ & midpoint & ord & PH-NSFD & midpoint & PH-NSFD \\
\midrule
$0.4$   & $1.099\cdot10^{0}$  & ---   & $3.69\cdot10^{-14}$ & $4.00\cdot10^{-15}$ & $3.78\cdot10^{-15}$\\
$0.2$   & $4.008\cdot10^{-1}$ & $1.46$& $3.27\cdot10^{-14}$ & $4.89\cdot10^{-15}$ & $1.06\cdot10^{-14}$\\
$0.1$   & $1.186\cdot10^{-1}$ & $1.76$& $4.85\cdot10^{-14}$ & $2.67\cdot10^{-15}$ & $2.13\cdot10^{-14}$\\
$0.05$  & $3.081\cdot10^{-2}$ & $1.94$& $4.76\cdot10^{-14}$ & $2.02\cdot10^{-14}$ & $1.95\cdot10^{-14}$\\
$0.025$ & $7.775\cdot10^{-3}$ & $1.99$& $3.02\cdot10^{-14}$ & $8.24\cdot10^{-14}$ & $1.71\cdot10^{-14}$\\
\bottomrule
\end{tabular}
\end{table}

Table~\ref{tab:exp1} confirms Corollary~\ref{cor:conservation} (both schemes conserve $H$ to round-off, independently of $h$) and Proposition~\ref{prop:exactscalar} (the \pHNSFD{} scheme is exact).
The midpoint rule converges with order two, and its unwrapped phase lag increases linearly  with $T$. 
However, the final-time state error does not necessarily grow monotonically because the trajectories are periodic.

Over a long horizon, for $h=0.1$, $T=1000$, the final-time relative state errors are
\begin{equation*}
    \mathrm{err}_{\mathrm{mid}}=1.91\cdot10^{-1},
    \qquad
    \mathrm{err}_{\pHNSFD{}}=1.89\cdot10^{-11} ,
\end{equation*}
a reduction by ten orders of magnitude at identical cost per step.
The moderate $\mathrm{err}_{\mathrm{mid}}$ hides a midpoint phase lag $\omega T-\frac{2T}{h}\arctan\frac{\omega h}{2}\approx6.63$ rad, more than one revolution (the \pHNSFD{} scheme has none), so a final-state norm is not a phase diagnostic.
Energy conservation and dynamic consistency are different properties.
The denominator provides the second one for one scalar function evaluation.

\subsection{Dissipative mode: $\tanh$ versus $\exp$ denominators}\label{ssec:exp2}
We integrate the scalar pH system with
$H(x) = \tfrac{1}{2}Q x^2$, $Q>0$, $J=0$ and $R>0$, namely
$\dot{x} = -RQx$, with $\alpha = RQ = 3$
to $T=1$, comparing $\varphi=h$, $\phi_{\mathrm{E}} = (1-\mathrm{e}^{-\alpha h})/\alpha$ and $\phi_{\tanh} = \frac{2}{\alpha} \tanh(\frac{\alpha h}{2})$.

\begin{table}[h]
\centering
\caption{Scalar dissipative pH system, $\alpha=3$, $T=1$: relative error and observed order.
The Mickens exponential denominator degrades the scheme to first order, exactly as predicted by Proposition~\ref{prop:tanhvsexp}. The $\tanh$ denominator is exact.}
\label{tab:exp2}
\small
\begin{tabular}{@{}lccccc@{}}
\toprule
$h$ & $\phi=h$ & ord & $\phi=\phi_{\mathrm{E}}$ & ord & $\phi=\phi_{\tanh}$\\
\midrule
$0.25$      & $1.426\cdot10^{-1}$ & ---    & $1.313\cdot10^{0}$  & ---    & $1.39\cdot10^{-16}$\\
$0.125$     & $3.528\cdot10^{-2}$ & $2.01$ & $6.122\cdot10^{-1}$ & $1.10$ & $8.36\cdot10^{-16}$\\
$0.0625$    & $8.797\cdot10^{-3}$ & $2.00$ & $2.940\cdot10^{-1}$ & $1.06$ & $2.65\cdot10^{-15}$\\
$0.03125$   & $2.198\cdot10^{-3}$ & $2.00$ & $1.439\cdot10^{-1}$ & $1.03$ & $2.79\cdot10^{-15}$\\
$0.015625$  & $5.493\cdot10^{-4}$ & $2.00$ & $7.113\cdot10^{-2}$ & $1.02$ & $5.16\cdot10^{-15}$\\
\bottomrule
\end{tabular}
\end{table}

Table~\ref{tab:exp2} is the numerical counterpart of Proposition~\ref{prop:tanhvsexp}. 
On the tested step sizes, the exponential denominator gives
larger errors than the standard choice $\varphi(h)=h$.
%
For this midpoint-type discretization of the scalar decay equation, the amplification factor is reproduced exactly by $\varphi_{\tanh}$, whereas $\varphi_E$ reduces the convergence order to one.

\subsection{Damped oscillator with a matrix denominator}\label{ssec:exp3}
We consider Example~\ref{ex:damped} with $\omega=2$, $\gamma=0.5$, $x_0=(1,0)^\top$, $T=20$, comparing $\phi=h$ with the exact matrix denominator $\Phi_{\mathrm{M}}(h)=2A^{-1}\tanh\bigl(\frac{h}{2}A\bigr)$.

For $h=0.4$, $0.2$, $0.1$ and $0.05$, the relative error at $T$ of the scheme with $\phi=h$ is $1.386$, $4.211\cdot10^{-1}$, $1.148\cdot10^{-1}$ and $2.915\cdot10^{-2}$ (observed orders $1.72$, $1.87$, $1.98$), whereas the matrix denominator gives errors between $3.8\cdot10^{-14}$ and $5.1\cdot10^{-14}$.
In both cases the numerical 
(discrete) energy decreases monotonically for all tested step sizes (Corollary~\ref{cor:dissipation}).

As a check of Theorem~\ref{thm:locallyexactisph}, at $h=0.4$ the eigenvalues of $R_h=-\frac{1}{h}Q^{-1}\sym(QM_h)Q^{-1}$ are $\{0.0076,\,0.557\}$, compared with $\{0,\,0.5\}$ for $R$.

\subsection{Nonlinear pendulum: period error}\label{ssec:exp4}
We integrate Example~\ref{ex:pendulum} with $p_0=0$ using the AVF discrete gradient \eqref{eq:pendulum_dgrad} and compare the standard AVF method, $\phi=h$, with the AVF-based \pHNSFD{} method with $\phi(h)=2\tan(h/2)$, built from the linearized frequency $\omega_0=1$ and locally exact at the origin. 
Both step sizes satisfy $h<\pi/2$, so by Theorem~\ref{thm:pendulum} each step is uniquely solvable and the solution remains in its libration well.
%
The numerical period is estimated by averaging 39 consecutive periods. 
We set $\tau_0=0$ and record the next 39 downward zero crossings of $p$, locating each crossing by linear interpolation between consecutive time steps. 
The reported estimate is
\begin{equation*}
   T_{\mathrm{num}}=\frac{\tau_{39}-\tau_0}{39}.
\end{equation*}
The maximal relative energy drift is measured in a separate run
with $N_E=\operatorname{round}(200/h)$ time steps:
\begin{equation*}
    \max_{0\le n\le N_E}\frac{|H_n-H_0|}{|H_0|}.
\end{equation*}
The corresponding final observation times are $200.10$ for $h=0.30$ and $199.95$ for $h=0.15$.

\begin{table}[h]
\centering
\caption{Nonlinear pendulum: AVF method ($\varphi(h)=h$) versus AVF-based
PH-NSFD method.
Relative period error and maximal relative energy drift over $N_E=\operatorname{round}(200/h)$ time steps. 
Energy is conserved to round-off in all reported runs.
The denominator built from the linearization at the origin gives an approximately 40-fold reduction in period error for small amplitudes and no improvement for large amplitudes, consistent with Remark~\ref{rem:localaccuracy}.}\label{tab:exp4}
\small
\begin{tabular}{@{}llcccc@{}}
\toprule
$q_0$ & $h$ & scheme & $T_{\mathrm{num}}$ & rel.\ period error &
$\max_n|H_n-H_0|/H_0$\\
\midrule
\multirow{4}{*}{$0.3$}
 & $0.30$ & AVF & $6.3647851$ & $7.29\cdot10^{-3}$ & $1.74\cdot10^{-14}$\\
 & $0.30$ & \pHNSFD{}& $6.3176735$ & $1.64\cdot10^{-4}$ & $1.24\cdot10^{-14}$\\
 & $0.15$ & AVF & $6.3302763$ & $1.83\cdot10^{-3}$ & $2.49\cdot10^{-14}$\\
 & $0.15$ & \pHNSFD{}& $6.3184458$ & $4.21\cdot10^{-5}$ & $9.94\cdot10^{-15}$\\
\cmidrule(l){2-6}
 & \multicolumn{5}{l}{exact period $T(0.3)=6.3187115$}\\
\midrule
\multirow{4}{*}{$1.5$}
 & $0.30$ & AVF & $7.3282427$ & $3.75\cdot10^{-3}$ & $2.75\cdot10^{-15}$\\
 & $0.30$ & \pHNSFD{}& $7.2736121$ & $3.73\cdot10^{-3}$ & $3.70\cdot10^{-15}$\\
 & $0.15$ & AVF & $7.3077020$ & $9.36\cdot10^{-4}$ & $3.11\cdot10^{-15}$\\
 & $0.15$ & \pHNSFD{}& $7.2940206$ & $9.37\cdot10^{-4}$ & $2.75\cdot10^{-15}$\\
\cmidrule(l){2-6}
 & \multicolumn{5}{l}{exact period $T(1.5)=7.3008649$}\\
\bottomrule
\end{tabular}
\end{table}

Table~\ref{tab:exp4} is instructive in both directions.  
%
For $q_0=0.3$ the pendulum is close to its linear regime, and the denominator reduces the period error by factors of approximately 44.5 at $h=0.3$ and 43.5 at $h=0.15$, at the cost of one scalar evaluation per step size.
For $q_0=1.5$ the amplitude-dependent frequency deviates substantially from $\omega_0=1$ and a denominator built from the linearization at the origin brings no benefit: the two errors have opposite sign but nearly equal magnitude.  

This observation motivates further investigation of amplitude-dependent or locally adapted denominator functions.
Their effect on the period error at larger amplitudes remains to be quantified in a future work.

\subsection{Forced system: power-balance residual}\label{ssec:exp5}
For the forced pH system of Example~\ref{ex:damped} with $\omega=2$, $\gamma=0.3$, $G=(0,1)^\top$, $u(t)=\sin(2t)$, $x_0={(1,\,0.5)^\top}$, $\phi(h)={h}$, $h=0.2$, $T=20$, the residual
\eqref{eq:residual} satisfies
\begin{equation*}
    \max_{0\le n<N}|\mathcal{R}_{H,n}| = 5.19\cdot10^{-15},
\end{equation*}
%
with $y_{\nph}=\Gb_n^\top\dgrad H_n$ and $u_{\nph}=u(t_n+\frac{h}{2})$. 
This residual provides a numerical consistency check of the
discrete power-balance \eqref{eq:discrete_power_balance}. 

\subsection{Additional structural tests}\label{ssec:exp6}
\paragraph{Forced equilibrium}
For $\dot{x}=-x+u$, $u=1$, $x^\ast=1$ (Example~\ref{ex:congruence_shift}), one step from $x^\ast$ leaves $x^\ast$ unchanged (to round-off) for the scalar scheme with $\varphi(h)=2\tanh(h/2)$ and for left multiplication with $L_h = (1+h^2)^2$, whereas the congruence form with $P=(1+h^2)^2$ gives the drifts $9.6\cdot10^{-4}$, $1.12\cdot10^{-1}$, $6.67\cdot10^{-1}$ and $1.54$ for $h=0.1$, $0.5$, $1$ and $2$. 
The left-multiplied and congruence schemes use $\varphi(h)=h$.

\paragraph{Casimir preservation}
For the rigid body of Example~\ref{ex:rigid_discrete} with $(I_1,I_2,I_3)=(1,2,3)$, $x_0=(1,0.5,-0.3)^\top$, $\varphi(h)=h=0.1$ and $1000$ steps, the relative Casimir drift is $5.3\cdot10^{-14}$ for $\rho=0$ and $3.3\cdot10^{-14}$ for $\rho=0.1$. 
The relative final-time energy deviation
$|H_N-H_0|/|H_0|$ is $2.6\cdot10^{-14}$ for $\rho=0$  and the energy decreases monotonically by 61\% for $\rho=0.1$.

%
For the rigid-body experiment, we used Newton iteration with a forward finite-difference Jacobian with increment $10^{-7}$, initialized at $x_n$. 
Each step was accepted only when the Euclidean norm of the unscaled step residual was below the tolerance $10^{-14}$, up to a maximum of 50 Newton updates.

\paragraph{Pendulum step beyond the cap}
For $(q_n,p_n) = (3,0)$, close to the separatrix, the numerical root search for \eqref{eq:pendulum_scalar}  found one solution with $q_{n+1}=q_n+d\in(-\pi,\pi)$ for $\varphi\in\{1,1.9,1.99\}$, and three solutions, two of which give $q_{n+1}\notin(-\pi,\pi)$, for $\varphi\in\{2.5,4,2\pi\}$.
For example, when $\varphi=2.5$, the three values are $q_{n+1}\approx0.629,\,3.856,\,4.980$.
These results demonstrate that, that uniqueness and preservation of the libration well can fail beyond the sufficient bound $0<\varphi(h)<2$ in Theorem~\ref{thm:pendulum}.


\subsection{Qualitative comparison with existing integrators}\label{ssec:comparison}
Table~\ref{tab:comparison} compares the proposed scheme qualitatively with classical, discrete-gradient \cite{Gonzalez1996, McLachlanQuispelRobidoux1999, QuispelMcLaren2008, CelledoniEtAl2012}, exponential \cite{HochbruckOstermann2010} and locally exact \cite{Cieslinski2013camwa, Cieslinski2013pla} integrators.
For existing structure-preserving pH integrators see the references in Section~\ref{Sec1}.

\begin{table}[h]
{\centering
\caption{Qualitative comparison of guarantees for the specified constructions under the stated assumptions; ``no'' means ``not guaranteed in general''. 
An em dash ('---') indicates that no general assertion is assigned to that property in this comparison.
``PB'' = exact discrete power balance; ``passive'' = passivity for all solved steps without step-size restriction; ``elem.'' = preservation of the equilibria and of their linear stability type (Definition~\ref{def:elemstab}), under the conditions stated in the text; ``loc.\ ex.'' = reproduces the linearized flow at an equilibrium; ``DG'' = discrete gradient; the last two rows refer to this work.}
\label{tab:comparison}
\footnotesize
\setlength{\tabcolsep}{4pt}
\begin{tabular}{@{}lcccccc@{}}
\toprule
method & order & PB & passive & Casimirs & elem.\ & loc.\ ex.\\
\midrule
explicit Euler            & 1 & no  & no  & affine  & no  & no\\
implicit Euler            & 1 & no  & conv.\ $H$\textsuperscript{$\dagger$} & affine  & no\textsuperscript{$\ddagger$} & no\\
implicit midpoint         & 2 & quadr.\ $H$ & quadr.\ $H$ & lin./quadr. & w/ cond.\textsuperscript{$\P$} & no\\
Gauss collocation, order $2s$ & $2s$ & quadr.\ $H$ & quadr.\ $H$ & lin./quadr. & w/ cond.\textsuperscript{$\P$} & no\\
symmetric DG (AVF, Gonzalez) & 2\textsuperscript{$\S$} & yes & yes & w/ cond. & yes\textsuperscript{$\ast$} & no\\
Itoh--Abe DG & 1 & yes & yes & w/ cond. & --- & no\\
exponential integrator    & $\ge1$ & no & no & no & --- & method-dep.\\
locally exact DG \cite{Cieslinski2013pla} & 2 & canonical & --- & --- & w/ cond.$^{\diamond}$ & yes\\
PH-NSFD, scalar & 2\textsuperscript{$\S$} & yes & yes &
   w/ cond. & yes\textsuperscript{$\ast$} & fitted mode\\
PH-NSFD, matrix & method-dep.\textsuperscript{a} & yes & yes &
   --- & w/ cond. & linear\textsuperscript{$\parallel$}\\
\bottomrule
\end{tabular}
\par}

  \vspace{0.5em}
  \footnotesize
\textsuperscript{$\dagger$}\,Since $H(z)-H(x)\le\nabla H(z)^\top(z-x)=-h\mathcal{D}(z)+hy(z)^\top u$ for convex $H$; with excess numerical dissipation.\quad
\textsuperscript{$\ddagger$}\,Turns conservative centers into spurious sinks.\quad
\textsuperscript{$\P$}\,On step-size ranges where the selected step map is well defined, retains the equilibrium branches and has no spurious fixed points (cf.\ the resonant Gauss step after Definition~\ref{def:elemstab}).
For implicit midpoint, linearized solvability requires $I-\frac h2 A_*$ to be invertible.
For $s$-stage Gauss collocation, the corresponding condition is the invertibility of $I_s\otimes I_n-h\mathcal A_s\otimes A_*$,
where $\mathcal A_s$ is the Butcher matrix. \quad
$^{\diamond}$ Subject to solvability, selection of an equilibrium-retaining branch, and exclusion of step sizes producing spurious fixed points.\quad
\textsuperscript{$\S$}\,Under (DG2), e.g.\ $H\in C^3$ (Lemma~\ref{lem:dgradorder}), and, for PH-NSFD, $\phi(h)=h+\mathcal{O}(h^3)$.\quad
\textsuperscript{$\ast$}\,at nondegenerate energy minima, on step-size ranges on which the step map is well defined and retains the diagonal solutions (Theorem~\ref{thm:elemstab}).\quad 
\textsuperscript{a}\,The order depends on the construction: it is two under the second-order consistency assumptions of Section~\ref{Sec8.3}, exact for the homogeneous linear constructions of Theorems~\ref{thm:locallyexactisph} and~\ref{thm:williamson}, and four on linear problems for the rational construction of Proposition~\ref{prop:rational}, on the stated admissible
step-size ranges.\quad
\textsuperscript{$\parallel$}\,For the constructions of Theorems~\ref{thm:locallyexactisph} and~\ref{thm:williamson}, on their stated step-size ranges.
\end{table}
\section*{Conclusion}
We have proposed and analyzed a structure-preserving NSFD method for pH systems, which separates the preservation of the algebraic pH structure, achieved through discrete gradients and skew-symmetric and positive semidefinite structure operators, from the approximation of the time scale, controlled by a Mickens-type denominator function $\phi(h)$.

The resulting guarantees hold at three levels.
First, the exact discrete power balance, passivity, energy conservation in the lossless case and, under explicit algebraic conditions, Casimir preservation hold for every solved step, independently of the step size and the denominator.
Second, all equilibria are exact solutions of the step equation of the scalar class (and fixed points of every branch that retains them), and the Lyapunov, LaSalle and linear stability results hold under their stated assumptions.
At nondegenerate energy minima the linearization is a Cayley transform that preserves the linear stability type for every admissible positive denominator.
Third, exact matrix fitting for linear systems is part of a step-dependent extension. This extension involves agreement between the fixed-point and equilibrium sets within the stated nonresonant step-size ranges.
For other structured extensions, such as the congruence construction, equilibrium compatibility must be checked separately, including at unforced equilibria that are not critical points of $H$.
%

Because the discrete-gradient scheme is of midpoint type, the appropriate denominators are $\frac{2}{\lambda}\tanh\bigl(\frac{\lambda h}{2}\bigr)$ rather than the exponential family of Euler-type NSFD schemes, which generically reduces the order from two to one.
For damped modes, a scalar phase-matching denominator reproduces the rotation angle of the exact flow.
The exact matrix denominator and its $[1/2]$ Pad\'e approximation, which coincides with the two-stage Gauss method, admit step-dependent pH representations with positive semidefinite dissipation; and a modal construction is exactly energy-preserving and frequency-matching in every mode for conservative linear systems.
The numerical experiments confirm the theory, with exactness for linear pH systems, a ten-order-of-magnitude reduction of the final-time error after $1000$ time units for the conservative oscillator and a forty-fold reduction of the period error for the near-linear pendulum, at the cost of one scalar function evaluation.

Future work will investigate the use of NSFD denominators within existing structure-preserving discretizations of pH differential-algebraic equations \cite{KinonMorandinSchulze2026, MehrmannMorandin2019}, as well as their potential application in structure-preserving, data-driven identification of such systems \cite{Zwerschke2026}.
Additionally, extensions to spatially semi-discretized pH partial differential equations \cite{Jaeschke2022} will be explored.


\section*{Acknowledgements} 
\noindent
M.~Ehrhardt acknowledges funding by the Deutsche Forschungsgemeinschaft (DFG, German Research Foundation) -- Project-ID 531152215 -- CRC 1701.


\begin{appendix}
\section{Assessment of the numerator function \texorpdfstring{$\psi(h)$}{psi(h)}}\label{app:numerator}
\renewcommand{\thedefinition}{A.\arabic{definition}} 

Mickens' general NSFD ansatz replaces $\dot{x}(t_n)$ by $(x_{n+1}-\psi(h)x_n)/\varphi(h)$.
For smooth $x$, $\varphi(h)=h+\mathcal{O}(h^2)$ and $\psi(h)-1=\mathcal{O}(h)$,
\begin{equation}\label{eq:psi_consistency}
    \dfrac{x(t+h)-\psi(h)x(t)}{\varphi(h)}
    = \dot{x}(t)-\frac{\psi(h)-1}{h}\,x(t)+o(1)
    \qquad (h\to0^+),
\end{equation}
so consistency requires $\psi(h)=1+o(h)$, a first-order approximation of the derivative requires $\psi(h)=1+O(h^2)$, and, for an otherwise second-order midpoint-type scheme
with $\varphi(h) = h+O(h^3)$, retaining second order requires $\psi(h) =1+O(h^3)$ in general.
Inserting the ansatz into the \pHNSFD{} scheme gives the generalized scheme
\begin{equation}\label{eq:psi_scheme}
    \dfrac{x_{n+1}-\psi(h)x_n}{\varphi(h)} = \bigl(\Jb_n-\Rb_n\bigr)\dgrad H_n + \Gb_n u_{\nph},
    \qquad
    y_{\nph} = \Gb_n^{\top}\dgrad H_n .
\end{equation}

\begin{proposition}[Rigidity of the scalar numerator]\label{prop:rigidity}
Consider \eqref{eq:psi_scheme} with $\phi(h)>0$ and discrete operators satisfying the diagonal consistency conditions \eqref{eq:consistency_gradient} and \eqref{eq:discrete_consistency}.
If the scheme preserves a prescribed nonzero equilibrium $(x^\ast,u^\ast)$ of \eqref{eq:ph_system} at each step, then $\psi(h)=1$ for every step size $h$ considered.
\end{proposition}
\begin{proof}
For $x_n=x_{n+1}=x^\ast$ and $u_{\nph}=u^\ast$, the right-hand side of \eqref{eq:psi_scheme} vanishes by diagonal consistency and \eqref{eq:equilibrium}.
Hence $(1-\psi(h))x^\ast/\phi(h)=0$ and $\psi(h)=1$.
\end{proof}

The following further obstructions explain the structural reason behind this result.

\paragraph{Loss of the exact power balance}
The increment of \eqref{eq:psi_scheme} 
has the additional term $(\psi(h)-1)x_n$, so the discrete chain rule \eqref{eq:discrete_gradient} yields
\begin{equation}\label{eq:modified-power-balance}
    H(x_{n+1})-H(x_n)
    = -\phi(h)\,\mathcal{D}_n
    + \phi(h)\, y_{\nph}^{\top}u_{\nph}
    + \mathcal{R}_{\psi,n},
    \quad
    \mathcal{R}_{\psi,n} := \bigl(\psi(h)-1\bigr)\dgrad H_n^{\top}x_n .
\end{equation}
The residual $\mathcal{R}_{\psi,n}$ is an artificial power term without pH interpretation, so Theorem~\ref{thm:powerbalance} is lost.
%
Similarly, even under the Casimir conditions \eqref{eq:casimirconditions} one obtains $C(x_{n+1})-C(x_n)=\bigl(\psi(h)-1\bigr)\dgrad C_n^{\top}x_n$, so exact Casimir preservation is destroyed as well.

\paragraph{Linear analysis}
For $\dot{\xi} = A_\ast\xi$ and the midpoint discrete gradient, the generalized scheme has, for an eigenvalue $\lambda$ of $A_\ast$, the amplification factor $r_{\lambda}(h)=\bigl(\psi(h)+\frac{\phi(h)}{2}\lambda\bigr)/\bigl(1-\frac{\phi(h)}{2}\lambda\bigr)$, so $\psi$ changes the numerator of the Cayley-type stability function.
For a conservative mode $\lambda=\iu\omega\neq0$ and real $\phi$, $\psi$, one has $|r_{\iu\omega}(h)|^2=\bigl(\psi(h)^2+\frac{1}{4}\phi(h)^2\omega^2\bigr)/\bigl(1+\frac{1}{4}\phi(h)^2\omega^2\bigr)$, so exact conservation of the quadratic Hamiltonian requires $\psi(h)^2=1$, and continuity with $\psi(0)=1$ selects $\psi\equiv1$.

\end{appendix}


\begin{thebibliography}{99}

\bibitem{Alalhareth2025}
F.~K.~Alalhareth, M.~Gupta, S.~Roy, H.~V.~Kojouharov, 
\textit{Second-order modified positive and elementary stable nonstandard numerical methods for $n$-dimensional autonomous differential equations}, 
Math.\ Meth.\ Appl.\ Sci.\ \textbf{48} (2025), 8037--8057.
 \href{https://doi.org/10.1002/mma.9560} {\path{doi:10.1002/mma.9560}}.


\bibitem{AnguelovEtAl2022}
R.~Anguelov, J.~M.-S.~Lubuma,
\textit{Second-order nonstandard finite difference schemes for a class of models in bioscience},
Preprint arXiv:2207.11618, 2022.
\href{https://doi.org/10.48550/arXiv.2207.11618}{\path{doi:10.48550/arXiv.2207.11618}}.

\bibitem{AnguelovLubuma2001}
R.~Anguelov, J.~M.-S.~Lubuma,
\textit{Contributions to the mathematics of the nonstandard finite difference method and applications},
Numer.\ Meth.\ Part.\ Diff.\ Eqs. \textbf{17}(5) (2001), 518--543.
 \href{https://doi.org/10.1002/num.1025}{\path{doi:10.1002/num.1025}}.
 
\bibitem{AnguelovLubuma2003}
R.~Anguelov, J.~M.-S.~Lubuma,
\textit{Nonstandard finite difference method by nonlocal approximation},
Math.\ Comput.\ Simul.\ \textbf{61} (2003), 465--475.
\href{https://doi.org/10.1016/S0378-4754(02)00106-4}{\path{doi:10.1016/S0378-4754(02)00106-4}}.

\bibitem{AouesEtAl2017}
S.~Aoues, M.~Di~Loreto, D.~Eberard, W.~Marquis-Favre,
\textit{Hamiltonian systems discrete-time approximation: Losslessness, passivity and composability},
Syst.\ Contr.\ Lett.\ \textbf{110} (2017), 9--14.
\href{https://doi.org/10.1016/j.sysconle.2017.10.003}{\path{doi:10.1016/j.sysconle.2017.10.003}}.

\bibitem{BeattieEtAl2018}
C.~Beattie, V.~Mehrmann, H.~Xu, H.~Zwart,
\textit{Linear port-Hamiltonian descriptor systems},
Math.\ Contr. Sign. Syst. \textbf{30} (2018), Art.~17.
\href{https://doi.org/10.1007/s00498-018-0223-3}{\path{doi:10.1007/s00498-018-0223-3}}.


\bibitem{CardosoRibeiroMatignonLefevre2021}
F.~L.~Cardoso-Ribeiro, D.~Matignon, L.~Lef\`evre,
\textit{A partitioned finite element method for power-preserving discretization of open systems of conservation laws},
IMA J.\ Math.\ Contr.\ Inform.\ \textbf{38} (2021), 493--533.
\href{https://doi.org/10.1093/imamci/dnaa038}{\path{doi:10.1093/imamci/dnaa038}}.



\bibitem{CelledoniEtAl2012}
E.~Celledoni, V.~Grimm, R.~I.~McLachlan, D.~I.~McLaren, D.~O'Neale, B.~Owren, G.~R.~W.~Quispel,
\textit{Preserving energy resp.\ dissipation in numerical PDEs using the ``average vector field'' method},
J.\ Comput.\ Phys.\ \textbf{231} (2012), 6770--6789.
\href{https://doi.org/10.1016/j.jcp.2012.06.022}{\path{doi:10.1016/j.jcp.2012.06.022}}.



\bibitem{CelledoniHoiseth2017}
E.~Celledoni, E.~H.~H{\o}iseth,
\textit{Energy-preserving and passivity-consistent numerical discretization of port-Hamiltonian systems},
Preprint arXiv:1706.08621, 2017.
\href{https://doi.org/10.48550/arXiv.1706.08621}{\path{doi:10.48550/arXiv.1706.08621}}.

\bibitem{Cieslinski2013camwa}
J.~L.~Cie\'sli\'nski,
\textit{Locally exact modifications of numerical schemes},
Comput.\ Math.\ Appl.\ \textbf{65} (2013), 1920--1938.
\href{https://doi.org/10.1016/j.camwa.2013.04.015}{\path{doi:10.1016/j.camwa.2013.04.015}}.

\bibitem{Cieslinski2013pla}
J.~L.~Cie\'sli\'nski,
\textit{Locally exact modifications of discrete gradient schemes},
Phys.\ Lett.\ A \textbf{377} (2013), 592--597.
\href{https://doi.org/10.1016/j.physleta.2013.01.005}{\path{doi:10.1016/j.physleta.2013.01.005}}.




\bibitem{Conte2026}
D.~Conte, G.~Pagano, T.~Rold\'an,
\textit{High-order nonstandard methods preserving mass and positivity for production-destruction systems},
ESAIM: Math. Model. Numer. Anal., 60(6) (2026) 2741--2767.
\href{https://doi.org/10.1051/m2an/2026080}{\path{doi:10.1051/m2an/2026080}}.



\bibitem{DahlbyOwrenYaguchi2011}
M.~Dahlby, B.~Owren, T.~Yaguchi,
\textit{Preserving multiple first integrals by discrete gradients},
J.\ Phys.\ A: Math.\ Theor.\ \textbf{44} (2011), 305205.
\href{https://doi.org/10.1088/1751-8113/44/30/305205}{\path{doi:10.1088/1751-8113/44/30/305205}}.


\bibitem{Dambrosio2026}
  R.~D’Ambrosio, S.~Di~Donato, 
  \textit{Discrete gradient $\theta$-methods for port-Hamiltonian systems}, 
  In: O.~Gervasi, et al.\ Computational Science and Its Applications -- ICCSA 2025 Workshops. ICCSA 2025. Lecture Notes in Computer Science 15887, Springer, Cham, 2026, pp.~225--236. \href{https://doi.org/10.1007/978-3-031-97589-9_16} {\path{doi:10.1007/978-3-031-97589-9_16}}.

\bibitem{DangHoang2018}
Q.~A.~Dang, M.~T.~Hoang, 
\textit{Exact Finite Difference Schemes for Three-Dimensional Linear Systems with Constant Coefficients}, 
Vietnam J.\  Math.\ \textbf{46} (2018), 471--492.
\href{https://doi.org/10.1007/s10013-017-0249-8}{\path{doi:10.1007/s10013-017-0249-8}}.


\bibitem{DangHoang2020}
Q.~A.~Dang, M.~T.~Hoang,
\textit{Positive and elementary stable explicit nonstandard Runge-Kutta methods for a class of autonomous dynamical systems},
Int.\ J.\ Comput.\ Math.\ \textbf{97} (2020), 2036--2054. 
\href{https://doi.org/10.1080/00207160.2019.1677895}{\path{doi:10.1080/00207160.2019.1677895}}.




\bibitem{DuindamEtAl2009}
V.~Duindam, A.~Macchelli, S.~Stramigioli, H.~Bruyninckx (eds.),
\textit{Modeling and Control of Complex Physical Systems: The Port-Hamiltonian Approach},
Springer, Berlin, 2009.
\href{https://doi.org/10.1007/978-3-642-03196-0}{\path{doi:10.1007/978-3-642-03196-0}}.

\bibitem{Egger2019}
H.~Egger,
\textit{Structure preserving approximation of dissipative evolution problems},
Numer.\ Math.\ \textbf{143} (2019), 85--106.
\href{https://doi.org/10.1007/s00211-019-01050-w}{\path{doi:10.1007/s00211-019-01050-w}}.


\bibitem{EhrhardtMickens2013}
M.~Ehrhardt, R.~E.~Mickens,
\textit{A nonstandard finite difference scheme for convection-diffusion equations having constant coefficients},
Appl.\ Math.\ Comput.\ \textbf{219} (2013), 6591--6604.
\href{https://doi.org/10.1016/j.amc.2012.12.068}{\path{doi:10.1016/j.amc.2012.12.068}}.



 \bibitem{Ehrhardt2026}
 M.~Ehrhardt, M.~Günther, D.~Sevcovic, 
 Structure-Preserving Coupling and Decoupling of Port-Hamiltonian Systems, 
 Appl.\ Math.\ Lett.\ \textbf{177} (2026), 109894. 
 \href{https://doi.org/10.1016/j.aml.2026.109894}{\path{doi:10.1016/j.aml.2026.109894}}.


\bibitem{FalaizeHelie2016}
A.~Falaize, T.~H\'elie,
\textit{Passive guaranteed simulation of analog audio circuits: A port-Hamiltonian approach},
Appl.\ Sci.\ \textbf{6} (2016), 273.
\href{https://doi.org/10.3390/app6100273}{\path{doi:10.3390/app6100273}}.

\bibitem{Gautschi1961}
W.~Gautschi,
\textit{Numerical integration of ordinary differential equations based on trigonometric polynomials},
Numer.\ Math.\ \textbf{3} (1961), 381--397.
\href{https://doi.org/10.1007/BF01386037}{\path{doi:10.1007/BF01386037}}.


\bibitem{GiesselmannKarsaiTscherpel2025}
J.~Giesselmann, A.~Karsai, T.~Tscherpel,
\textit{Energy-consistent Petrov-Galerkin time discretization of port-Hamiltonian systems},
SMAI J.\ Comput.\ Math.\ \textbf{11} (2025), 335--367.
\href{https://doi.org/10.5802/smai-jcm.127}{\path{doi:10.5802/smai-jcm.127}}.


\bibitem{Gonzalez1996}
O.~Gonzalez,
\textit{Time integration and discrete Hamiltonian systems},
J.\ Nonlin. Sci.\ \textbf{6} (1996), 449--467.
\href{https://doi.org/10.1007/BF02440162}{\path{doi:10.1007/BF02440162}}.



\bibitem{GrimmQuispel2005}
V.~Grimm, G.~R.~W.~Quispel,
\textit{Geometric integration methods that preserve Lyapunov functions},
BIT Numer.\ Math.\ \textbf{45} (2005), 709--723.
\href{https://doi.org/10.1007/s10543-005-0034-z}{\path{doi:10.1007/s10543-005-0034-z}}.



\bibitem{HairerLubichWanner2006}
E.~Hairer, C.~Lubich, G.~Wanner,
\textit{Geometric Numerical Integration: Structure-Preserving Algorithms for Ordinary Differential Equations}, 2nd ed.,
Springer, Berlin, 2006.
\href{https://doi.org/10.1007/3-540-30666-8}{\path{doi:10.1007/3-540-30666-8}}.





\bibitem{Hoang2023}
M.~T.~Hoang, 
\textit{An exact finite difference scheme for a three-layer beam model}, 
J.\ Differ.\ Eqs. Appl.\ \textbf{29} (2023), 419--431.
\href{https://doi.org/10.1080/10236198.2023.2209204}{\path{doi:10.1080/10236198.2023.2209204}}.

\bibitem{HoangMatthias2024}
M.~T.~Hoang, M.~Ehrhardt, 
\textit{A second-order nonstandard finite difference method for a general Rosenzweig--MacArthur predator--prey model},  
J.\ Comput.\ Appl.\ Math.\ \textbf{444} (2024), 115752.
\href{https://doi.org/10.1016/j.cam.2024.115752}{\path{doi:10.1016/j.cam.2024.115752}}.


\bibitem{HoangMatthias2026}
M.~T.~Hoang, M.~Ehrhardt,
\textit{Differential equation models for infectious diseases: Mathematical modeling, qualitative analysis, numerical methods and applications}, 
SeMa Journal \textbf{83} (2026), 435--470.
\href{https://doi.org/10.1007/s40324-025-00404-9}{\path{doi:10.1007/s40324-025-00404-9}}.


\bibitem{HoangMatthias2026NSFD}
M.~T.~Hoang, M.~Ehrhardt,  
\textit{A generalized second-order positivity-preserving numerical method for non-autonomous dynamical systems with applications}, 
Appl.\ Math.\ Comput.\ \textbf{524} (2026), 130029.
\href{https://doi.org/10.1016/j.amc.2026.130029}{\path{doi:10.1016/j.amc.2026.130029}}.




\bibitem{HochbruckLubich1999}
M.~Hochbruck, C.~Lubich,
\textit{A Gautschi-type method for oscillatory second-order differential equations},
Numer.\ Math.\ \textbf{83} (1999), 403--426.
\href{https://doi.org/10.1007/s002110050456}{\path{doi:10.1007/s002110050456}}.


\bibitem{HochbruckOstermann2010}
M.~Hochbruck, A.~Ostermann,
\textit{Exponential integrators},
Acta Numer.\ \textbf{19} (2010), 209--286.
\href{https://doi.org/10.1017/S0962492910000048}{\path{doi:10.1017/S0962492910000048}}.



\bibitem{ItohAbe1988}
T.~Itoh, K.~Abe,
\textit{Hamiltonian-conserving discrete canonical equations based on variational difference quotients},
J.\ Comput.\ Phys.\ \textbf{76} (1988), 85--102.
\href{https://doi.org/10.1016/0021-9991(88)90132-5}{\path{doi:10.1016/0021-9991(88)90132-5}}.


\bibitem{IxaruVandenBerghe2004}
L.~Gr.~Ixaru, G.~Vanden~Berghe,
\textit{Exponential Fitting},
Kluwer, Dordrecht, 2004.
\href{https://doi.org/10.1007/978-1-4020-2100-8}
{\path{doi:10.1007/978-1-4020-2100-8}}.

\bibitem{Jaeschke2022}
J.~Jäschke, M.~Ehrhardt, M.~Günther, B.~Jacob, 
\textit{Discrete port-Hamiltonian Coupled Heat Transfer},
in: M. Ehrhardt, M. Günther (eds.), Progress in Industrial Mathematics at ECMI 2021, The European Consortium for Mathematics in Industry, Springer, 2022, pp. 439-445.
\href{https://doi.org/10.1007/978-3-031-11818-0_57}
{\path{doi:10.1007/978-3-031-11818-0_57}}.



\bibitem{KinonEtAl2023}
P.~L.~Kinon, T.~Thoma, P.~Betsch, P.~Kotyczka,
\textit{Discrete nonlinear elastodynamics in a port-Hamiltonian framework}, Proc. Appl. Math. Mech. \textbf{23}(2023), e202300144.
\href{https://doi.org/10.1002/pamm.202300144}{\path{doi:10.1002/pamm.202300144}}.

\bibitem{KinonMorandinSchulze2026}
P.~L.~Kinon, R.~Morandin, P.~Schulze,
\textit{Discrete gradient methods for port-Hamiltonian differential-algebraic equations},
Appl.\ Numer.\ Math.\ \textbf{223} (2026), 45--75.
\href{https://doi.org/10.1016/j.apnum.2025.12.006}{\path{doi:10.1016/j.apnum.2025.12.006}}.



\bibitem{KojouharovEtAl2021}
H.~V.~Kojouharov, S.~Roy, M.~Gupta, F.~Alalhareth, J.~M.~Slezak,
\textit{A second-order modified nonstandard theta method for one-dimensional autonomous differential equations},
Appl.\ Math.\ Lett.\ \textbf{112} (2021), 106775.
\href{https://doi.org/10.1016/j.aml.2020.106775}{\path{doi:10.1016/j.aml.2020.106775}}.




\bibitem{Kotyczka2019}
P.~Kotyczka,
\textit{Numerical Methods for Distributed Parameter Port-Hamiltonian Systems},
TUM University Press, Munich, 2019.
\href{https://doi.org/10.14459/2019md1510230}
{\path{doi:10.14459/2019md1510230}}.

\bibitem{KotyczkaLefevre2018}
P.~Kotyczka, L.~Lef\`evre,
\textit{Discrete-time port-Hamiltonian systems based on Gauss-Legendre collocation},
IFAC-PapersOnLine \textbf{51} (2018), 125--130.
\href{https://doi.org/10.1016/j.ifacol.2018.06.035}
{\path{doi:10.1016/j.ifacol.2018.06.035}}.

\bibitem{KotyczkaLefevre2019}
P.~Kotyczka, L.~Lef\`evre,
\textit{Discrete-time port-Hamiltonian systems: A definition based on symplectic integration},
Syst.\ Contr.\ Lett.\ \textbf{133} (2019), 104530.
\href{https://doi.org/10.1016/j.sysconle.2019.104530}
{\path{doi:10.1016/j.sysconle.2019.104530}}.

\bibitem{KotyczkaMaschkeLefevre2018}
P.~Kotyczka, B.~Maschke, L.~Lef\`evre,
\textit{Weak form of Stokes-Dirac structures and geometric discretization of port-Hamiltonian systems},
J.\ Comput.\ Phys.\ \textbf{361} (2018), 442--476.
\href{https://doi.org/10.1016/j.jcp.2018.02.006}
{\path{doi:10.1016/j.jcp.2018.02.006}}.


\bibitem{LaSalle1976}
J.~P.~LaSalle,
\textit{The Stability of Dynamical Systems},
SIAM, Philadelphia, 1976.
\href{https://doi.org/10.1137/1.9781611970432}
{\path{doi:10.1137/1.9781611970432}}.

%





\bibitem{LeimkuhlerReich2004}
B.~Leimkuhler, S.~Reich,
\textit{Simulating Hamiltonian Dynamics},
Cambridge University Press, Cambridge, 2004.
\href{https://doi.org/10.1017/CBO9780511614118}
{\path{doi:10.1017/CBO9780511614118}}.


 \bibitem{Maamar2024}
 M.H.~Maamar, M.~Ehrhardt, L.~Tabharit, 
 \textit{A Nonstandard Finite Difference Scheme for a Time-Fractional Model of Zika Virus Transmission}, 
 Math. Biosci. Engrg.\ \textbf{21}(1) (2024), 924--962.
 \href{https://doi.org/10.3934/mbe.2024039}{\path{doi:10.3934/mbe.2024039}}.
%




\bibitem{McLachlanQuispelRobidoux1999}
R.~I.~McLachlan, G.~R.~W.~Quispel, N.~Robidoux,
\textit{Geometric integration using discrete gradients},
Philos.\ Trans.\ R.\ Soc.\ Lond.\ A \textbf{357} (1999), 1021--1045.
\href{https://doi.org/10.1098/rsta.1999.0363}
{\path{doi:10.1098/rsta.1999.0363}}.



\bibitem{McLachlanQuispelTse2009}
R.~I.~McLachlan, G.~R.~W.~Quispel, P.~S.~P.~Tse,
\textit{Linearization-preserving self-adjoint and symplectic integrators},
BIT Numer.\ Math.\ \textbf{49} (2009), 177--197.
\href{https://doi.org/10.1007/s10543-009-0214-3}
{\path{doi:10.1007/s10543-009-0214-3}}.


\bibitem{McLachlanQuispel2014}
R.~I.~McLachlan, G.~R.~W.~Quispel,
\textit{Discrete gradient methods have an energy conservation law},
Discr.\ Contin.\ Dyn.\ Syst.\ \textbf{34} (2014), 1099--1104.
\href{https://doi.org/10.3934/dcds.2014.34.1099}
{\path{doi:10.3934/dcds.2014.34.1099}}.


\bibitem{MehrmannMorandin2019}
V.~Mehrmann, R.~Morandin,
\textit{Structure-preserving discretization for port-Hamiltonian descriptor systems},
in: Proc.\ 58th IEEE Conf.\ Decision and Control (CDC), Nice, 2019, pp.~6863--6868.
\href{https://doi.org/10.1109/CDC40024.2019.9030180}
{\path{doi:10.1109/CDC40024.2019.9030180}}.




\bibitem{MehrmannUnger2023}
V.~Mehrmann, B.~Unger,
\textit{Control of port-Hamiltonian differential-algebraic systems and applications},
Acta Numer.\ \textbf{32} (2023), 395--515.
\href{https://doi.org/10.1017/S0962492922000083}
{\path{doi:10.1017/S0962492922000083}}.



\bibitem{Mickens1994}
R.~E.~Mickens,
\textit{Nonstandard Finite Difference Models of Differential Equations},
World Scientific, Singapore, 1994.
\href{https://doi.org/10.1142/2081}
{\path{doi:10.1142/2081}}.

\bibitem{MickensRamadhani1994}
R.~E.~Mickens,  I.~Ramadhani, 
\textit{Finite-difference schemes having the correct linear stability properties for all finite step-sizes III}, 
Comput.\ Math.\ Appl. \textbf{27} (1994), 77--84.
\href{https://doi.org/10.1016/0898-1221(94)90056-6}
{\path{doi:10.1016/0898-1221(94)90056-6}}.


\bibitem{Mickens1996energy}
R.~E.~Mickens,
\textit{Construction of finite difference schemes for coupled nonlinear oscillators derived from a discrete energy function},
J.\ Differ.\ Eq.\ Appl.\ \textbf{2} (1996), 185--193.
\href{https://doi.org/10.1080/10236199608808053}
{\path{doi:10.1080/10236199608808053}}.



\bibitem{Mickens2002}
R.~E.~Mickens,
\textit{Nonstandard finite difference schemes for differential equations},
J.\ Differ. Eqs.\ Appl.\ \textbf{8} (2002), 823--847.
\href{https://doi.org/10.1080/1023619021000000807}
{\path{doi:10.1080/1023619021000000807}}.


\bibitem{Mickens2005dc}
R.~E.~Mickens,
\textit{Dynamic consistency: a fundamental principle for constructing nonstandard finite difference schemes for differential equations},
J.\ Differ. Eqs.\ Appl.\ \textbf{11} (2005), 645--653.
\href{https://doi.org/10.1080/10236190412331334527}{\path{doi:10.1080/10236190412331334527}}.



\bibitem{Mickens2005hamilton}
R.~E.~Mickens,
\textit{A numerical integration technique for conservative oscillators combining
nonstandard finite-difference methods with a Hamilton's principle},
J.\ Sound Vib.\ \textbf{285} (2005), 477--482.
\href{https://doi.org/10.1016/j.jsv.2004.09.027}{\path{doi:10.1016/j.jsv.2004.09.027}}.
%

\bibitem{Mickens2007denominator}
R.~E.~Mickens,
\textit{Calculation of denominator functions for nonstandard finite difference schemes for differential equations satisfying a positivity condition},
Numer.\ Meth. Part. Diff. Eqs \textbf{23} (2007), 672--691.
\href{https://doi.org/10.1002/num.20198}{\path{doi:10.1002/num.20198}}.

\bibitem{Mickens2013}
R.~E.~Mickens, T.~M.~Washington, 
\textit{NSFD discretizations of interacting population models satisfying conservation laws},
Comput.\ Math.\ Appl. \textbf{66} (2013), 2307--2316.
\href{https://doi.org/10.1016/j.camwa.2013.06.011}
{\path{doi:10.1016/j.camwa.2013.06.011}}.



\bibitem{Mickens2020}
R.~E.~Mickens,
\textit{Nonstandard Finite Difference Schemes: Methodology and Applications},
World Scientific, Singapore, 2020.
\href{https://doi.org/10.1142/11891}{\path{doi:10.1142/11891}}.

\bibitem{MickensWashington2020}
R.~E.~Mickens, T.~M.~Washington,
\textit{A note on a positivity preserving nonstandard finite difference scheme for a modified parabolic reaction-advection-diffusion PDE},
J.\ Differ.\ Eqs.\ Appl.\ \textbf{26} (2020), 1423--1427.
\href{https://doi.org/10.1080/10236198.2020.1841755}{\path{doi:10.1080/10236198.2020.1841755}}.

\bibitem{Moench2024}
M.~M\"onch, N.~Marheineke,
\textit{Fourth-order force-gradient splitting for linear port-Hamiltonian systems},
Proc. Appl. Math. Mech. \textbf{24}(3) (2024), e202400132.
\href{https://doi.org/10.1002/pamm.202400132}{\path{doi:10.1002/pamm.202400132}}.

\bibitem{Morais2024}
G.~Morais Rodrigues Costa, M.~Lobosco, M.~Ehrhardt, R.~F.~Reis, \textit{Mathematical Analysis and a Nonstandard Scheme for a Model of the Immune Response against COVID-19}, 
in: A.~Gumel (ed.), \textit{Mathematical and Computational Modeling of Phenomena Arising in Population Biology and Nonlinear Oscillations}, 
Contemporary Mathematics \textbf{793},
Amer. Math. Soc., Providence, RI, 2024, pp. 251--270. 
\href{https://doi.org/10.1090/conm/793/15881}{\path{doi:10.1090/conm/793/15881}}.


\bibitem{MoullaLefevreMaschke2012}
R.~Moulla, L.~Lef\`evre, B.~Maschke,
\textit{Pseudo-spectral methods for the spatial symplectic reduction of open systems of conservation laws},
J.\ Comput.\ Phys.\ \textbf{231} (2012), 1272--1292.
\href{https://doi.org/10.1016/j.jcp.2011.10.008}{\path{doi:10.1016/j.jcp.2011.10.008}}.



\bibitem{NortonQuispel2014}
R.~A.~Norton, G.~R.~W.~Quispel,
\textit{Discrete gradient methods for preserving a first integral of an ordinary differential equation},
Discr.\ Contin.\ Dyn.\ Syst.\ \textbf{34} (2014), 1147--1170.
\href{https://doi.org/10.3934/dcds.2014.34.1147}{\path{doi:10.3934/dcds.2014.34.1147}}.


\bibitem{Oettinger2018}
H.~C.~\"Ottinger,
\textit{GENERIC integrators: Structure preserving time integration for thermodynamic systems},
J.\ Non-Equilib.\ Thermodyn.\ \textbf{43} (2018), 89--100.
\href{https://doi.org/10.1515/jnet-2017-0034}{\path{doi:10.1515/jnet-2017-0034}}.



%






\bibitem{Patidar2016}
K.~C.~Patidar,
\textit{Nonstandard finite difference methods: recent trends and further developments},
J.\ Differ.\  Eqs.\ Appl.\ \textbf{22} (2016), 817--849.
\href{https://doi.org/10.1080/10236198.2016.1144748}{\path{doi:10.1080/10236198.2016.1144748}}.




\bibitem{QuispelMcLaren2008}
G.~R.~W.~Quispel, D.~I.~McLaren,
\textit{A new class of energy-preserving numerical integration methods},
J.\ Phys.\ A: Math.\ Theor.\ \textbf{41} (2008), 045206.
\href{https://doi.org/10.1088/1751-8113/41/4/045206}{\path{doi:10.1088/1751-8113/41/4/045206}}.



\bibitem{Roeger2008}
L.-I.~W.~Roeger,
\textit{Exact difference schemes for two-dimensional linear systems with constant coefficients},
J.\ Comput.\ Appl.\ Math.\ \textbf{219} (2008), 102--109.
\href{https://doi.org/10.1016/j.cam.2007.07.015}{\path{doi:10.1016/j.cam.2007.07.015}}.


\bibitem{SanzSernaCalvo1994}
J.~M.~Sanz-Serna, M.~P.~Calvo,
\textit{Numerical Hamiltonian Problems},
Chapman \& Hall, London, 1994.


\bibitem{Schulze2023}
P.~Schulze,
\textit{Structure-preserving time discretization of port-Hamiltonian systems via discrete gradient pairs},
In: K.~Burnecki, J.~Szwabiński, M.~Teuerle (eds.),
Progress in Industrial Mathematics at ECMI 2023. 
Mathematics in Industry \textbf{41}. Springer, Cham, 2026, pp. 501--511.
\href{https://doi.org/10.1007/978-3-032-20404-2_47}{\path{doi:10.1007/978-3-032-20404-2_47}}.


\bibitem{SongoloBidegaray2024}
M.~E.~Songolo, B.~Bid\'egaray-Fesquet,
\textit{Extending nonstandard finite difference scheme rules to systems of nonlinear ODEs with constant coefficients},
J.\ Differ. Eqs.\ Appl.\ \textbf{30} (2024), 545--570.
\href{https://doi.org/10.1080/10236198.2023.2300726}{\path{doi:10.1080/10236198.2023.2300726}}.

\bibitem{Strogatz}
S. H.~Strogatz, 
\textit{Nonlinear dynamics and chaos}, 2nd edition,
Westview press, 2015; reprinted CRC press 2018.
\href{https://doi.org/10.1201/9780429492563}{\path{doi:10.1201/9780429492563}}.

\bibitem{TalasilaEtAl2006}
V.~Talasila, J.~Clemente-Gallardo, A.~J.~van~der~Schaft,
\textit{Discrete port-Hamiltonian systems},
Syst. Contr. Lett.\ \textbf{55} (2006), 478--486.
\href{https://doi.org/10.1016/j.sysconle.2005.10.001}{\path{doi:10.1016/j.sysconle.2005.10.001}}.



\bibitem{Treibert2022}
S.~Treibert, H.~Brunner, M.~Ehrhardt, 
\textit{A Nonstandard Finite Difference Scheme for the SVICDR model to predict COVID-19 dynamics}, 
Math. Biosci. Engrg.\ \textbf{ 19}(2) (2022), 1213--1238. 
\href{https://doi.org/10.3934/mbe.2022056}{\path{doi:10.3934/mbe.2022056}}.







\bibitem{TrenchantEtAl2018}
V.~Trenchant, H.~Ramirez, Y.~Le~Gorrec, P.~Kotyczka,
\textit{Finite differences on staggered grids preserving the port-Hamiltonian structure with application to an acoustic duct},
J.\ Comput.\ Phys.\ \textbf{373} (2018), 673--697.
\href{https://doi.org/10.1016/j.jcp.2018.06.051}{\path{doi:10.1016/j.jcp.2018.06.051}}.



\bibitem{vanderSchaft2017}
A.~J.~van~der~Schaft,
\textit{$L_2$-Gain and Passivity Techniques in Nonlinear Control}, 3rd ed.,
Springer, Cham, 2017. \href{https://doi.org/10.1007/978-3-319-49992-5}{\path{doi:10.1007/978-3-319-49992-5}}.



\bibitem{vanderSchaftJeltsema2014}
A.~J.~van~der~Schaft, D.~Jeltsema,
\textit{Port-Hamiltonian systems theory: An introductory overview},
Found.\ Trends Syst.\ Contr.\ \textbf{1} (2014), 173--378.
\href{https://doi.org/10.1561/2600000002}{\path{doi:10.1561/2600000002}}.


\bibitem{Tanhc}
E.~W.~Weisstein, 
\textit{Tanhc Function},
From MathWorld--A Wolfram Resource,
\url{https://mathworld.wolfram.com/TanhcFunction.html}
(last access Nov. 25, 2026)


\bibitem{Willems1972}
J.~C.~Willems,
\textit{Dissipative dynamical systems. Part I: General theory},
Arch.\ Ration.\ Mech.\ Anal.\ \textbf{45} (1972), 321--351.
\href{https://doi.org/10.1007/BF00276493}{\path{doi:10.1007/BF00276493}}.


\bibitem{Williamson1936}
J.~Williamson,
\textit{On the algebraic problem concerning the normal forms of linear dynamical systems},
Amer.\ J.\ Math.\ \textbf{58} (1936), 141--163.
\href{https://doi.org/10.2307/2371062}{\path{doi:10.2307/2371062}}.




\bibitem{ZinihiEhrhardtSidiAmmi2026}
A.~Zinihi, M.~Ehrhardt, M.~R.~Sidi~Ammi,
\textit{A nonstandard finite difference scheme for an SEIQR epidemiological PDE model},
Appl.\ Math.\ Comput.\ \textbf{520} (2026), 129953.
\href{https://doi.org/10.1016/j.amc.2026.129953}{\path{doi:10.1016/j.amc.2026.129953}}.

\bibitem{Zwerschke2026}
T. Zwerschke, M. Ehrhardt, M. Günther, C. Totzeck, 
\textit{Structure-Preserving Data-Driven Identification of Port-Hamiltonian Differential-Algebraic Systems},
Preprint arXiv:2609.07187, 2026.
\href{https://doi.org/10.48550/arXiv.2609.07187}{\path{doi:10.48550/arXiv.2609.07187}}.

\end{thebibliography}
\end{document}